\pdfoutput=1
\RequirePackage[l2tabu, orthodox]{nag} %
\documentclass{diss}

\usepackage{subcaption}
\usepackage[utf8]{inputenc}   %
\usepackage{theorems}         %
\usepackage{Baskervaldx}

\usepackage[Baskervaldx]{newtxmath}
\usepackage[scr=rsfs,bb=ams]{mathalpha}

\usepackage{amsfonts}

\DeclareMathAlphabet{\mathbbe}{U}{bbold}{m}{n}

\usepackage[activate={true,nocompatibility},final,tracking=true,kerning=true,spacing=true]{microtype}

\usepackage{amsmath}
\usepackage[warning,all,nodollardollar]{onlyamsmath}      %
\usepackage[mathic]{mathtools}
\usepackage[low-sup]{subdepth}

\usepackage{xspace}
\usepackage{notation}
\usepackage{mydiagrams}
\usepackage[inline]{enumitem}
\usepackage{csquotes}
\usepackage[style=alphabetic,
 maxnames=4,
 maxbibnames=99,
 uniquename=false,doi=false,url=false,isbn=false,]{biblatex}
\addglobalbib{references.bib}

\DeclareSourcemap{
  \maps[datatype=bibtex]{
    \map[overwrite]{ %
      \step[fieldsource=fjournal,fieldtarget=journaltitle]
    }
  }  
}
\usepackage[unicode=true,colorlinks=true,hypertexnames=false]{hyperref}
\usepackage{bookmark}
\usepackage{xcolor}
\hypersetup{allcolors={DodgerBlue4}}%

\usepackage{needspace}

\title{Pasting in Simplicial Categories}
\author{Tobias Columbus}

\begin{document}
\frontmatter
\maketitle
\tableofcontents

\mainmatter

\chapter*{Overview}

This text is concerned with a pasting theorem for categories enriched over
quasi-categories. The notion of pasting goes back to the theory of
$2$-categories, where one frequently encounters diagrams such as 
  \begin{equation}\label{dia:pasting}
    \begin{tikzpicture}[diagram]
      \matrix[objects] {
        |(0)| A \& |(1)| B \& |(2)| C \& \\[-1ex]
                \& |(3)| X \& |(4)| Y \& |(5)| Z\smash{.} \\
      };
      \path[maps,->]
        (0) edge[out=-90,in=180] node[below left] (1c) {$u$} (3)
        (0) edge node[above] {$a $} (1)
        (1) edge node[above] {$b$} (2)
        (2) edge node[right] {$g$} (4)
        (1) edge node[left] {$f$} (3)
        (3) edge node[below] {$v$} (4)
        (2) edge[out=0,in=90] node[above right] (1b) {$c$} (5)
        (4) edge node[below] {$w$} (5)
      ;
      \node[anchor=center, rotate=-135] at (barycentric cs:1=0.5,2=0.5,3=0.5,4=0.5) (arrow) {$\Rightarrow$};
      \node[left,maps] at (arrow) {$\phi$};
      \node[anchor=center,rotate=-135] at (barycentric cs:0=0.5,3=0.5,1=0.5,1c=0.5) (eps) {$\Rightarrow$};
      \node[anchor=center,rotate=-135] at (barycentric cs:2=0.5,4=0.5,5=0.5,1b=0.5) (eta) {$\Rightarrow$};
      \node[left,maps] at (eps) {$\psi$};
      \node[right,maps] at (eta) {$\xi$};
    \end{tikzpicture}
  \end{equation}
  The axioms of a $2$-category ensure that one can compose two cells like
  e.\,g.\ $\psi$ and $\phi$ above.  In this way, one obtains a composite
  \begin{equation*}
    cba \xto{\xi \cdot ba} wgba \xto{  w \cdot \phi \cdot a} wvfa \xto{wv \cdot \psi} wvu
  \end{equation*}
  of the whole diagram. Pasting refers to this cobbling of cells in
  $2$-categories and a diagram such as \eqref{dia:pasting} is consequently called a pasting
  diagram.
  
  In the case of \eqref{dia:pasting}, the order in which one
  forms the composites of two such cells so as to obtain a composite of the
  whole diagram is actually unique. However, there are other diagrams such as
  e.\,g.\
  \begin{equation*}
    \begin{tikzpicture}[diagram]
      \matrix[objects] {
        |(a)| A \&[+1em] |(b)| B \&[+1em] |(c)| C, \\
      };
      \path[maps,->] 
        (a) edge[bend left] node[above] (u) {$u$} (b)
        (a) edge[bend right] node[below] (v) {$v$} (b)
        (b) edge[bend left] node[above] (f) {$f$} (c)
        (b) edge[bend right] node[below] (g) {$g$} (c)
      ;
      \node[anchor=center,rotate=-90] (phi) at  ($ (u) ! 0.5 ! (v) $) {$\Rightarrow$}; 
      \node[anchor=center,rotate=-90] (psi) at  ($ (f) ! 0.5 ! (g) $) {$\Rightarrow$}; 
      \node[left,font=\scriptsize] at (phi) {$\phi$};
      \node[left,font=\scriptsize] at (psi) {$\psi$};
    \end{tikzpicture}
  \end{equation*}
  where we have the two possibilities 
  \begin{equation*}
    fu \xto{\psi \cdot u} gu \xto{g \cdot \phi} gv  \qquad\text{and}\qquad
    fu \xto{f \cdot \phi} fv \xto{\psi \cdot v} gv 
  \end{equation*}
  to obtain a composite of the whole diagram. These two compositions actually
  do coincide by the axioms of a $2$-category. A \emph{pasting theorem} simply
  asserts that any well-formed diagram admits a composite that is independent
  from the order in which we form compositions of individual cells in the diagram.
  Such a theorem was conjectured for $2$-categories in
  \cite{kelly-street;review-of-the-elements-of-2-categories} and proven by John
  Power in \cite{power;a-2-categorical-pasting-theorem}. 

  The work of Power in \cite{power;a-2-categorical-pasting-theorem} has two
  facets. First of all, the formulation of a pasting theorem requires a formal
  definition of what exactly a well-formed diagram actually is. Power chooses
  certain plane graphs as his model of a diagram in a $2$-category. Building on
  this notion of a diagram and its basic combinatorial properties, Power then
  goes on to show that any labeling of such a diagram in a $2$-category admits
  a  uniquely determined composite.
  
  We prove a pasting theorem for a specific model of
  $(\infty,2)$-categories, namely categories enriched over quasi-categories.
  We reuse the model of diagrams introduced by Power and extend his work from
  $2$-categories to $(\infty,2)$-categories. However, we have to cope with quite some technical difficulties due to the fact that composition in
  $(\infty,2)$-categories is only unital and associative up to higher
  dimensional invertible cells. Moreover, we cannot even expect a diagram to have a
  unique composite but have to settle with composites that are unique up to
  higher dimensional invertible cells. More precisely, we show that the space
  of compositions of a given diagram is nonempty and contractible: 
  \begin{restatable*}{mytheorem}{pastingthm}\label{thm:pasting}
    Consider a globular graph\/ $G$ and a category\/ $\dA$ enriched over
    quasi-categories.  The space\/ $C(\Lambda)$ of compositions of a given
    labeling\/ $\Lambda$ of\/ $G$ in\/ $\dA$ is a nonempty contractible Kan complex.
  \end{restatable*}
  A labeling as in the statement of the above theorem or the work of Power
  simply refers to a sensible assignment of cells in the category $\dA$ to any
  cell in the graph $G$.

  Let us say a few words about our model of diagrams before we comment on the
  proof of \autoref{thm:pasting}. We model our diagrams on certain plane graphs
  that we call globular graphs for they typically look as follows:
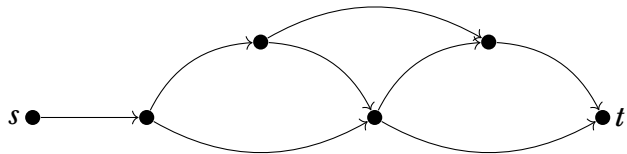
\begin{equation*}
\begin{tikzpicture}
            \node[fill,circle,inner sep=2pt] (1) at (0,0) {};
            \node[fill,circle,inner sep=2pt] (2) at (1.5,0) {};
            \node[fill,circle,inner sep=2pt] (3) at (4.5,0) {}; 
            \node[fill,circle,inner sep=2pt] (4) at (7.5,0) {}; 
            \node[fill,circle,inner sep=2pt] (5) at (3,1) {}; 
            \node[fill,circle,inner sep=2pt] (6) at (6,1) {}; 
            \node[left] at (1) {$s$};
            \node[right] at (4) {$t$};
            \draw (1) edge[->] (2);
            \draw (2) edge[->,bend right]  (3);
            \draw (2) edge[->,bend left]  (5);
            \draw (5) edge[->,bend left]  (3);
            \draw (3) edge[->,bend left]  (6);
            \draw (6) edge[->,bend left]  (4);
            \draw (3) edge[->,bend right](4);
            \draw (5) edge[->,bend left] (6);
        \end{tikzpicture}
\end{equation*}
Although the concept of globular graphs is due to
\cite{power;a-2-categorical-pasting-theorem}, we have to extend and adapt it to
our model of $(\infty,2)$-categories. This means in particular that we have to
associate with any such graph a family $X_{a,b}$ of simplicial sets indexed by
pairs $(a,b)$ of vertices of the graph. To this end, we introduce the nerve of
such a graph and develop a pictorial calculus to describe its simplices. A simplex in $\nerve(G)$ essentially corresponds to a picture such as 
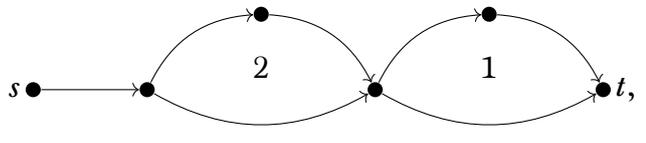
\begin{equation*}
\begin{tikzpicture}
            \node[fill,circle,inner sep=2pt] (1) at (0,0) {};
            \node[fill,circle,inner sep=2pt] (2) at (1.5,0) {};
            \node[fill,circle,inner sep=2pt] (3) at (4.5,0) {}; 
            \node[fill,circle,inner sep=2pt] (4) at (7.5,0) {}; 
            \node[fill,circle,inner sep=2pt] (5) at (3,1) {}; 
            \node[fill,circle,inner sep=2pt] (6) at (6,1) {}; 
            \node[left] at (1) {$s$};
            \node[right] at (4) {$t$,};
            \draw (1) edge[->] (2);
            \draw (2) edge[->,bend right]  (3);
            \draw (2) edge[->,bend left]  (5);
            \draw (5) edge[->,bend left]  (3);
            \draw (3) edge[->,bend left]  (6);
            \draw (6) edge[->,bend left]  (4);
            \draw (3) edge[->,bend right](4);
            \node at (3,0.3) {$2$};
            \node at (6,0.3) {$1$};
        \end{tikzpicture},
\end{equation*}
where the numbers in the faces should be thought of as specifying the order of
composition. In order to describe partial composites in terms of these
simplicial sets, we then go on to restrict the type of simplices that may
appear in this nerve. This line of thought eventually leads to our notion of
pasting diagram.

  Let us now comment on the proof of \autoref{thm:pasting}.  There are three
  main aspects to the proof and each aspect can be handled individually. In
  fact, the three main chapters of this text are more or less independent
  from each other and deal with one issue at a time. 
  \subsection*{Global aspect}
    If a pasting diagram $\Sigma$ satisfies a certain technical assumption,
    then we are able to associate with $\Sigma$ a simplicial category
    $\dC[\Sigma]$. For the maximal pasting diagram with underlying graph $G$
    one obtains for example the free $2$-category on this graph in the sense of
    Power.
    We use these simplicial categories $\dC[\Sigma]$ to model partial composites
    of diagrams in some simplicial category $\dA$ by simplicial functors
    $\dC[\Sigma] \to \dA$.  Given such a partially composed diagram $u\colon
    \dC[\Sigma] \to \dA$, we may add all missing composite cells to the
    abstract diagram $\Sigma$ and thus obtain a new diagram $\Pi$.  The problem
    of finding a composition of the concrete diagram $u\colon \dC[\Sigma] \to
    \dA$ then amounts to the problem of finding an extension of $u$ as in the
    diagram 
      \begin{equation}\label{eq:ext}
        \begin{tikzpicture}[diagram]
          \matrix[objects] {
            |(s)| \dC[\Sigma] \& |(a)| \dA \\
            |(p)| \dC[\Pi]\smash{.} \\
          };
          \path[maps,->] 
            (s) edge node[above] {$u$} (a)
            (s) edge node[left] {inclusion} (p)
            (p) edge[dashed] (a)
          ;
        \end{tikzpicture}
      \end{equation}
      We attack the problem of constructing such an extension by
      considering the more general problem of finding a simplicial functor
      $\dC[\Pi] \to \dA$ that renders a diagram such as
      \begin{equation}\label{eq:llp}
        \begin{tikzpicture}[diagram]
          \matrix[objects] {
            |(s)| \dC[\Sigma] \& |(a)| \dA \\
            |(p)| \dC[\Pi] \& |(b)| \dB \\
          };
          \path[maps,->] 
            (s) edge node[above] {$u$} (a)
            (s) edge node[left] {inclusion} (p)
            (p) edge[dashed] (a)
            (p) edge (b)
            (a) edge node[right] {$p$} (b)
          ;
        \end{tikzpicture}
      \end{equation}
      commutative. This solves the original problem of finding a composition of
      $u \colon \dC[\Sigma] \to \dA$ as long as we allow the unique
      simplicial functor $\dA \to \ast$ to feature as the functor $p$.

      Consider a class $\sR$ of maps of simplicial sets that contains all
      isomorphisms and let $\sL = {}^\lifts \sR$ be the class of maps that have
      the left lifting property against all maps in $\sR$. We solve the
      problem~\eqref{eq:llp} for those diagrams of simplicial categories, where
      $p$ is a local $\sR$-functor, that is, $p \colon \dA(a,a') \to
      \dB(pa,pa')$ is an $\sR$-map for all $a,a' \in \dA$. Our solution to
      \eqref{eq:llp} then takes the form of the following theorem:
      \begin{restatable*}{mytheorem}{globalthm}\label{thm:global}
      The functor\/ $\dC[\Sigma] \to \dC[\Pi]$ induced by an inclusion of complete
      pasting diagrams has the left lifting property against all local\/
      $\sR$-functors if and only if the map
      \begin{equation*}
        \nerve( \Sigma_{x,y} \hc \Pi_{x,y} ) \to \nerve(\Pi_{x,y})
      \end{equation*}
      is an\/ $\sL$-map for all vertices\/ $x,y \in \Sigma$.
    \end{restatable*}
    The diagrams $\Sigma_{x,y} \hc \Pi_{x,y}$ appearing in the statement of
    \autoref{thm:global} are certain intermediate diagrams between $\Sigma$ and
    $\Pi$ that are introduced at the very end of
    \autoref{sub:nerves_of_pasting_diagrams}.

    \subsection*{Local aspect}
    Now suppose that the classes $\sL$ and $\sR$ from above are the classes of
    mid anodyne maps and mid fibrations, respectively. In this case, the unique
    functor $\dA \to \pt$ from a simplicial category $\dA$ to the terminal
    simplicial category is a local $\sR$-functor if and only if $\dA$ is
    enriched over quasi-categories. \autoref{thm:global} hence reduces the
    problem \eqref{eq:ext} of finding a composition of a diagram $u\colon
    \dC[\Sigma] \to \dA$ in some category enriched over quasi-categories to the
    problem of showing that certain maps
    \begin{equation*}
      \nerve(\Sigma_{x,y} \hc \Pi_{x,y} ) \to \nerve(\Pi_{x,y})
    \end{equation*}
    are mid anodyne. We identify a certain class of inclusions $\Sigma \to \Pi$ of
    pasting diagrams for which these maps are indeed mid anodyne. More
    precisely, we deduce that the maps in question are mid anodyne from the
    following theorem:
    \begin{restatable*}{mytheorem}{localthm}\label{thm:local}
      Let\/ $\Sigma \to \Pi$ be an inclusion of complete pasting diagrams such that
      both\/ $\Sigma$ and $\Pi$ contain all the interior faces of the underlying
      graph and are closed under taking subdivisions. Then
      \begin{equation*}
        \nerve(\Sigma) \to \nerve(\Pi)
      \end{equation*}
      is mid anodyne.
    \end{restatable*}
    
    \subsection*{Bootstrapping}
    There is still one problem left that did not yet appear in the above
    discussion. The simplicial categories $\dC[\Sigma]$ featuring as the domain
    of diagrams $\dC[\Sigma] \to \dA$ exist only under certain technical
    assumptions. In order to complete the proof of \autoref{thm:pasting} that
    we sketched so far, we have to make sure that we can associate with any
    labeling $\Lambda$ of a diagram a simplicial functor $\dC[\Sigma] \to \dA$.
    This is achieved by 
\begin{restatable*}{mytheorem}{labelingthm}\label{thm:labeling}
  Suppose that\/ $\Sigma$ is the minimal complete pasting diagram on some
  globular graph\/ $G$. The map
  \begin{equation*}
    \scat(\dC[\Sigma],\dA) \to L(G,\dA),\quad u \mapsto \Lambda_u,
  \end{equation*}
  that sends a simplicial functor\/ $u$ to its associated labeling\/ $\Lambda_u$ is a
  bijection.
\end{restatable*}

\newpage
The individual chapters correspond to these three aspects of our main
theorem. The mutual interdependence of the individual chapters can be seen in
the following diagram:
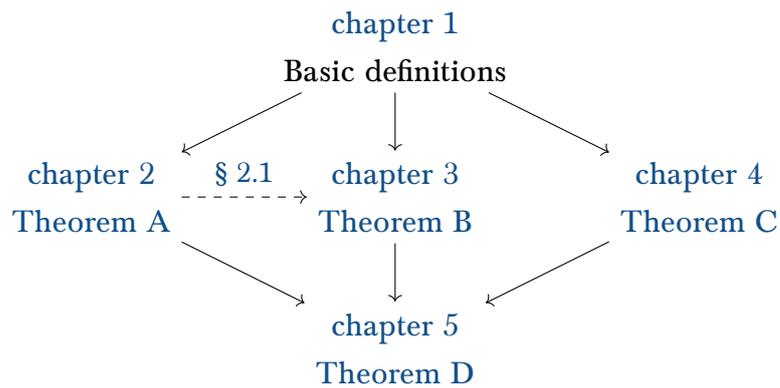
\begin{center}
  \begin{tikzpicture}
    \node[align=center] (pd) at (0,-4) {\autoref{chap:pasting-diagrams}\\Basic definitions};
    \node[align=center] (labeling) at (-4,-6) {\autoref{chap:labeling}\\\autoref{thm:labeling}};
    \node[align=center] (global) at (0,-6) {\autoref{chap:global}\\\autoref{thm:global}};
    \node[align=center] (local) at (4,-6) {\autoref{chap:local}\\\autoref{thm:local}};
    \node[align=center] (thm) at (0,-8) {\autoref{chap:thm}\\\autoref{thm:pasting}};
    \path[->]
      (pd) edge (labeling)
      (pd) edge (local)
      (pd) edge (global)
      (labeling) edge (thm)
      (local) edge (thm)
      (global) edge (thm)
      (labeling) edge[dashed] node[above,font=\small] {\autoref{sec:scat}} (global)
      ;
  \end{tikzpicture}
\end{center}

\subsection*{Related work}%
\label{sub:related_work_}

This text is a slightly abridged version of the author's unpublished 2017
thesis. The decision to finally release some version of this work into the wild
was sparked by the recent preprint
\cite{hackney-ozornova-riehl-rovelli;pasting} and private communication with
its authors Hackney, Ozornova, Riehl and Rovelli.\footnote{The author
appreciates both the pleasant interaction and the impetus to finally publish
this thesis in some form.}

In \cite{hackney-ozornova-riehl-rovelli;pasting}, the authors also prove a
pasting theorem for $(\infty,2)$-categories based on the notion of a pasting
diagram from \cite{power;a-2-categorical-pasting-theorem} and there is some
overlap between the work of Hackney, Ozorova, Riehl and Rovelli and this text.
Their Theorem~2.3.5 is similar to the combinatorial and pictorial description
of the nerve of a globular graph that we give in
\autoref{sub:nerves_of_globular_graphs} and their Proposition~4.1.3 is our
\autoref{rem:simplicial-computad}. 

Despite these similarities, there are some clear differences between
\cite{hackney-ozornova-riehl-rovelli;pasting} and the approach taken in this
text: Even though \cite[Theorem~4.1.10]{hackney-ozornova-riehl-rovelli;pasting}
more or less recovers our \autoref{thm:global} and \autoref{thm:local}, the
method of proof is different.  Whereas our argument in \autoref{chap:global} and
\ref{chap:local} is mostly combinatorial, the proof in
\cite{hackney-ozornova-riehl-rovelli;pasting} employs a new variant of a result
of Thomason on pushouts along Dwyer maps and their compatibility with the Joyal
model structure.

Moreover, instead of producing a space of composites via some concrete
simplicial enrichment of simplicial categories and proving its contractibility
as we do in \autoref{thm:pasting}, Hackney, Ozornova, Riehl and Rovelli deduce
a model-independent result by making good use of a model structure on
simplicial categories and the axiomatic theory of $(\infty,2)$-categories due
to \cite{unicity}. 

Thus, even though both \autoref{thm:pasting} and
\cite[Corollary~4.4.2]{hackney-ozornova-riehl-rovelli;pasting} assert that we
have contractible (and hence equivalent) spaces of composites, the spaces are
\emph{a priori} different and it is not clear to the author how one might
compare both results on a conceptual level.

\chapter{Globular Graphs and Pasting Diagrams}\label{chap:pasting-diagrams}

This chapter has essentially two parts. The first part comprises
\autoref{sub:globgraphs}--\ref{sec:globular_graphs_as_canonical_arities} and
is concerned with globular graphs. We review the definition of globular graphs such as 
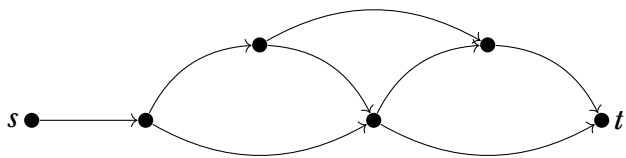
\begin{equation*}
\begin{tikzpicture}
            \node[fill,circle,inner sep=2pt] (1) at (0,0) {};
            \node[fill,circle,inner sep=2pt] (2) at (1.5,0) {};
            \node[fill,circle,inner sep=2pt] (3) at (4.5,0) {}; 
            \node[fill,circle,inner sep=2pt] (4) at (7.5,0) {}; 
            \node[fill,circle,inner sep=2pt] (5) at (3,1) {}; 
            \node[fill,circle,inner sep=2pt] (6) at (6,1) {}; 
            \node[left] at (1) {$s$};
            \node[right] at (4) {$t$};
            \draw (1) edge[->] (2);
            \draw (2) edge[->,bend right]  (3);
            \draw (2) edge[->,bend left]  (5);
            \draw (5) edge[->,bend left]  (3);
            \draw (3) edge[->,bend left]  (6);
            \draw (6) edge[->,bend left]  (4);
            \draw (3) edge[->,bend right](4);
            \draw (5) edge[->,bend left] (6);
        \end{tikzpicture}
\end{equation*}
from \cite{power;a-2-categorical-pasting-theorem} and discuss some technical
observations due to Power.  

As we eventually want to interpret these globular graphs as diagrams in categories enriched over
simplicial sets or quasi-categories, we then go on to define the nerve of a
globular graph and establish some of its technical properties. It turns out,
that the nerve of a globular graph appears as the nerve of the category
$F_2G(s,t)$, where $F_2G$ is the free $2$-category on the graph $G$. This means
in particular that $\nerve(G)$ contains the composite of any of its
$1$-simplices and is therefore inadequate as a model of the input data of a
vertical composition in a category enriched over quasi-categories.  In the
second part of this chapter we thus introduce pasting diagrams in
\autoref{sub:pasting_diagrams} and their nerves in
\autoref{sub:nerves_of_pasting_diagrams}. Our notion of pasting diagram is
built on the notion of a globular graph but allows for a specification of which
vertical composites are present in its nerve.

\section{Globular Graphs}\label{sub:globgraphs}

In this section, we review relevant parts of the work
\cite{power;a-2-categorical-pasting-theorem} of Power and extend it by some
definitions and mostly trivial observations of our own.

\subsection*{Globular graphs}
Let $G$ be a graph.  A vertex $s \in G$ is called a \emph{source} if there is a
directed path from~$s$ to any vertex $u \neq s$ of $G$. A \emph{target} in $G$
is a source in $G\op$.  An \emph{$st$-graph} is a nontrivial plane graph $G$
with unique source and target that are both incident to the exterior face.

A face $\phi$ of a plane graph is \emph{globular} if its clockwise oriented
boundary $\del \phi$ decomposes as $\del \phi = p \cdot q\op$ for two
nontrivial directed paths $p$ and $q$. If $\phi$ is an interior face, we call
$\dom \phi = p$ the \emph{domain} and $\cod \phi = q$ the \emph{codomain} of $\phi$.
However, we use the exact opposite convention for the exterior face $\epsilon$,
i.\,e.  $\dom \varepsilon = q$ and $\cod \varepsilon = p$ if $\del \varepsilon
= p \cdot q\op$.

As $s(p) = s(q)$
and $t(p) = t(q)$, we simply denote these vertices by $s(\phi)$ and $t(\phi)$,
respectively.  A \emph{globular graph} is an $st$-graph in which all faces are
globular.  

\begin{example}
  Consider the graph shown in \autoref{fig:globular-graph}. It certainly is an
  $st$-graph for the vertices marked $s$ and $t$ are its source and target,
  respectively. Moreover, the clockwise boundary $\del \phi$ of any face $\phi$
  of $G$ decomposes as $\del \phi = p \cdot q\op$ for two directed paths $p$
  and $q$, that is, all faces of $G$ are globular. Let us check this for the
  interior face $\phi$ marked in \autoref{fig:globular-graph} and the exterior
  face $\varepsilon$ of $G$. The clockwise boundary of $\phi$ is $e_2 \cdot e_5
  \cdot e_7\op = p \cdot q\op$ with $p = e_2 \cdot e_5$ and $q = e_7$. The
  face $\phi$ thus has $\dom \phi = e_2 \cdot e_5$ and $\cod \phi = e_7$.  The
  exterior face $\varepsilon$ has clockwise boundary $\del \varepsilon = e_4\op
  \cdot e_3\op \cdot e_2\op \cdot e_1\op \cdot e_1 \cdot e_7 \cdot e_8$ and
  this clockwise boundary decomposes up to rotation as $\del \varepsilon = p
  \cdot q\op$ with $p = e_1 \cdot e_7 \cdot e_8$ and $q = e_1  \cdot e_2 \cdot
  e_3 \cdot e_4$. We thus have $\dom \varepsilon= e_1 \cdot e_2 \cdot e_3 \cdot
  e_4$ and $\cod \varepsilon = e_1 \cdot e_7 \cdot e_8$.
  \begin{figure}
    \centering
        \begin{tikzpicture}
            \node[fill,circle,inner sep=2pt] (1) at (0,0) {};
            \node[fill,circle,inner sep=2pt] (2) at (1.5,0) {};
            \node[fill,circle,inner sep=2pt] (3) at (4.5,0) {}; 
            \node[fill,circle,inner sep=2pt] (4) at (7.5,0) {}; 
            \node[fill,circle,inner sep=2pt] (5) at (3,1) {}; 
            \node[fill,circle,inner sep=2pt] (6) at (6,1) {}; 
            \node[left] at (1) {$s$};
            \node[right] at (4) {$t$};
            \draw (1) edge[->] node[fill=white,inner sep=1.5pt,font=\scriptsize] {$e_1$} (2);
            \draw (2) edge[->,bend right] node[fill=white,inner sep=1pt,font=\scriptsize] {$e_7$} (3);
            \draw (2) edge[->,bend left] node[fill=white,inner sep=1pt,font=\scriptsize] {$e_2$} (5);
            \draw (5) edge[->,bend left] node[fill=white,inner sep=1pt,font=\scriptsize] {$e_5$} (3);
            \draw (3) edge[->,bend left] node[fill=white,inner sep=1pt,font=\scriptsize] {$e_6$} (6);
            \draw (6) edge[->,bend left] node[fill=white,inner sep=1pt,font=\scriptsize] {$e_4$} (4);
            \draw (3) edge[->,bend right] node[fill=white,inner sep=1pt,font=\scriptsize] {$e_8$} (4);
            \draw (5) edge[->,bend left] node[fill=white,inner sep=1pt,font=\scriptsize] {$e_3$} (6);
            \node[font=\scriptsize] at (barycentric cs:2=0.5,3=0.5,5=0.5) {$\phi$};
        \end{tikzpicture}
    \caption{An example of a globular graph.}\label{fig:globular-graph}
  \end{figure}
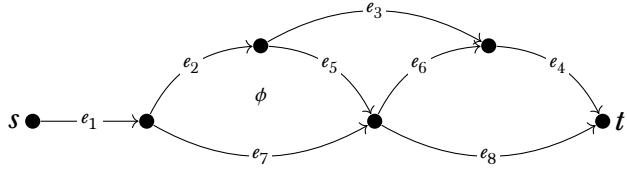
\end{example}

The following characterisation of globular graphs is taken from
\cite[Proposition~2.6]{power;a-2-categorical-pasting-theorem}. 
\begin{proposition}\label{prop:ps:power}
  A nontrivial $st$-graph\/ $G$ is globular if and only if it is acyclic, i.\,e.\ contains no directed cycles.
\end{proposition}

The boundary, domain and codomain of a globular graph are the boundary, domain
and codomain of its exterior face. Note that our conventions concerning the
domain and codomain of the exterior face now ensure that the domain of a face
does not depend upon whether we consider it as a subgraph or as a face.

\begin{definition}\label{def:wide}
  A globular subgraph $H \subseteq G$ is \emph{wide} if $s(H) = s(G)$ and $t(H) = t(G)$.
\end{definition}

\begin{definition}\label{def:subdivision}
  If $H \subseteq G$ is a globular subgraph of $G$, we call a subgraph $K \subseteq G$ a \emph{subdivision} of $H$ if $H \subseteq K$ and $\del H = \del K$. 
\end{definition}

\begin{definition}\label{def:label}
  A \emph{glob} $\gamma$ in a globular graph $G$ is a globular subgraph of $G$
  with the property that any of its edges is incident with the exterior face of
  $\gamma$.  A glob is \emph{nondegenerate} if it has at least one interior face
  and \emph{degenerate} otherwise. A glob is \emph{proper} if it is $2$-connected.
\end{definition}
\begin{example}
\strut
  \begin{enumerate}[label=(\alph*)]
    \item Any path $p$ in a globular graph $G$ is a degenerate glob.
    \item Any face $\phi$ of a globular graph $G$ is a nondegenerate, proper glob.
    \item Let $H \subseteq G$ be a globular subgraph of some globular graph
      $G$. The boundary $\del H$ of $H$ is a glob in $G$ that is degenerate if
      and only $H$ has no interior faces and proper if and only if $H$ is
      $2$-connected.
    \item Consider the graph $G$ shown in \autoref{subfig:glob-example-graph}
      and its subgraphs $\gamma_1$, $\gamma_2$ and $H$ shown in
      \autoref{subfig:glob-example-proper-glob},
      \ref{subfig:glob-example-nonproper-glob} and
      \ref{subfig:glob-example-no-glob}, respectively.  Both $\gamma_1$ and
      $\gamma_2$ are globs, but $H$ is not since it contains an edge that is
      not incident with the exterior face of $H$.
      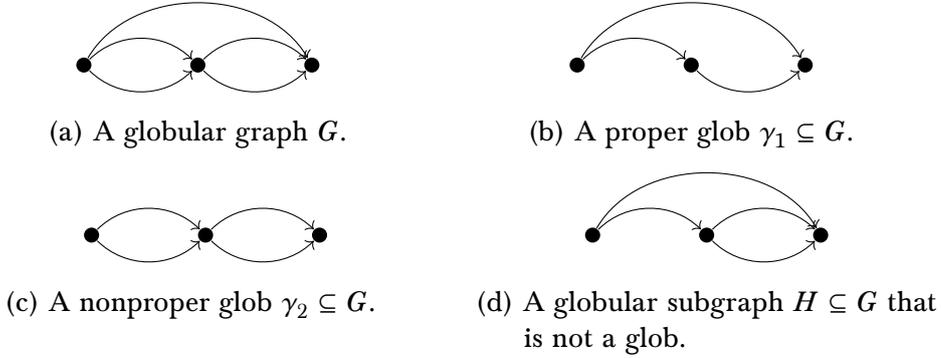
\begin{figure}
        \centering
        \begin{subfigure}[b]{0.45\textwidth}
          \centering
          \begin{tikzpicture}
            \node[fill,circle,inner sep=2pt] (s) at (0,0) {};
            \node[fill,circle,inner sep=2pt] (m) at (1.5,0) {};
            \node[fill,circle,inner sep=2pt] (t) at (3,0) {}; 
            \draw (s) edge[->,out=45,in=135] (m);
            \draw (s) edge[->,out=-45,in=-135] (m);
            \draw (m) edge[->,out=45,in=135] (t);
            \draw (m) edge[->,out=-45,in=-135] (t);
            \draw (s) edge[->,out=60,in=120] (t);
        \end{tikzpicture}
          \caption{A globular graph $G$.}\label{subfig:glob-example-graph}
        \end{subfigure}
        \hfil
        \begin{subfigure}[b]{0.45\textwidth}
          \centering
          \begin{tikzpicture}
            \node[fill,circle,inner sep=2pt] (s) at (0,0) {};
            \node[fill,circle,inner sep=2pt] (m) at (1.5,0) {};
            \node[fill,circle,inner sep=2pt] (t) at (3,0) {}; 
            \draw (s) edge[->,out=45,in=135] (m);
            \draw[white] (m) edge[->,out=45,in=135] (t);
            \draw (m) edge[->,out=-45,in=-135] (t);
            \draw (s) edge[->,out=60,in=120] (t);
        \end{tikzpicture}
          \caption{A proper glob $\gamma_1 \subseteq G$.}\label{subfig:glob-example-proper-glob}
        \end{subfigure}
        \newline
        \begin{subfigure}[b]{0.45\textwidth}
          \centering
          \begin{tikzpicture}
            \node[fill,circle,inner sep=2pt] (s) at (0,0) {};
            \node[fill,circle,inner sep=2pt] (m) at (1.5,0) {};
            \node[fill,circle,inner sep=2pt] (t) at (3,0) {}; 
            \draw (s) edge[->,out=45,in=135] (m);
            \draw (s) edge[->,out=-45,in=-135] (m);
            \draw (m) edge[->,out=45,in=135] (t);
            \draw (m) edge[->,out=-45,in=-135] (t);
        \end{tikzpicture}
          \caption{A nonproper glob $\gamma_2 \subseteq G$.\newline\strut}\label{subfig:glob-example-nonproper-glob}
        \end{subfigure}
        \hfil
        \begin{subfigure}[b]{0.45\textwidth}
          \centering
          \begin{tikzpicture}
            \node[fill,circle,inner sep=2pt] (s) at (0,0) {};
            \node[fill,circle,inner sep=2pt] (m) at (1.5,0) {};
            \node[fill,circle,inner sep=2pt] (t) at (3,0) {}; 
            \draw (s) edge[->,out=45,in=135] (m);
            \draw (m) edge[->,out=45,in=135] (t);
            \draw (m) edge[->,out=-45,in=-135] (t);
            \draw (s) edge[->,out=60,in=120] (t);
        \end{tikzpicture}
          \caption{A globular subgraph $H \subseteq G$ that is not a
          glob.}\label{subfig:glob-example-no-glob}
        \end{subfigure}
        \caption{Example of a globular graph $G$ and two of its globs $\gamma_1$ and $\gamma_2$.}
      \end{figure}
  \end{enumerate}
\end{example}
\begin{lemma}\label{lem:degenerate-glob}
  The following are equivalent for a glob $\gamma$:
  \begin{enumerate}[label=(\roman*),noitemsep]
    \item $\gamma$ is degenerate.
    \item $\gamma$ is a directed path in $G$.
    \item $\del\gamma\op$ is oriented clockwise.
  \end{enumerate}
\end{lemma}
\begin{proof}
  Let $\gamma$ be a degenerate glob with source $s$ and target $t$ in $G$.  If
  $\gamma$ is no directed path, then there are two $st$-paths $p \neq q$ in
  $\gamma$. The cycle $p \cdot q\op$ encloses at least one interior face of
  $\gamma$ --- a contradiction. As directed paths are certainly degenerate
  globs, this proves the equivalence of (i) and (ii). 

  The oriented boundary of a directed path $p$ is simply $\del p = p
  \cdot p\op$ and this cycle clearly satisfies $\del p\op = \del p$, which
  proves the implication \enquote{(ii) $\Rightarrow$ (iii)}.

  Now assume that both $\del\gamma$ and $\del\gamma\op$ are clockwise
  orientations of $\del\gamma$. Write $\del \gamma = p \cdot q\op$ and
  $\del\gamma\op = q\cdot p\op$. As clockwise orientations are unique up to
  cyclic rotation, we conclude that $p = q$, i.\,e.\ that $\gamma$ is a simple
  directed path in $G$.
\end{proof}

\begin{remark}
  If $u,v$ are two vertices of a globular graph $G$, the subgraph $G_{u,v}$
  consisting of all paths with source $u$ and target $v$ is either empty or
  itself a globular graph. If $\{u,v\} \neq \{s,t\}$, then $G_{u,v}$ is
  strictly smaller 
  than $G$ in terms of arrows and vertices.
\end{remark}

The following two lemmata are due to Power \cite{power;a-2-categorical-pasting-theorem}.
\begin{lemma}\label{lem:cyclic-order}
  Let $v$ be a vertex of a globular graph $G$. The clockwise cyclic order of
  edges around $v$ is $e_1 \prec \dots \prec e_r \prec d_1 \prec \dots \prec
  d_s$, where the edges $e_i$ are the edges with source $v$ and the edges $d_i$
  are the edges with target $v$. 
\end{lemma}

\begin{lemma}\label{lemma:glob-graph:face}
  Let $G$ be a globular graph with exterior face $\varepsilon$ and at least one interior
  face. There then exists an interior face $\phi$ of $G$ with $\dom(\phi)$ lying
  entirely in $\dom(G)$.
\end{lemma}

\subsection*{Joins of globular graphs}
The operation of glueing two globular graphs at their respective source and
target vertices very much looks like horizontal composition in $2$-categories.
In fact, this is precisely the role that this operation is going to play in the
later chapters of this work. In this paragraph, we just define this operation
and prove one technical lemma. 
\begin{definition}\label{def:join}
  The \emph{join} $G\join H$ of two globular graphs is the globular graph
  obtained by gluing $t(G)$ to $s(H)$. This is a well-defined globular graph as we consider
  topological equivalence classes of actual embeddings.
\end{definition}

\begin{remark}
  Note that the domain and codomain of $G \join H$ are given by $\dom(G \join
  H) = \dom (G) \cdot \dom (H)$ and $\cod(G \join H) = \cod (G) \cdot \cod
  (H)$. 
\end{remark}

\begin{remark}
  Consider a globular graph $G$ and some vertex $x \in G$.  Recall that
  $G_{s,x}$ and $G_{x,t}$ are the subgraphs of $G$ consisting of all the paths
  from $s$ to $x$ and from $x$ to $t$, respectively.  As $t(G_{s,x}) = x = 
  s(G_{x,t})$, we can form their join $G_{s,x} \join G_{x,t}$ and obtain a wide
  globular subgraph of $G$.
\end{remark}

We end this paragraph with a technical lemma on the interplay between joins and
intersections of globular subgraphs of the form $G_{u,v}$. 
\begin{lemma}\label{lem:intersection-of-xy-globs}
    Let $G$ be a globular graph with source $s$ and target $t$. Further let
    $x,y \in G$ be two vertices.  
    The intersection
    \begin{equation*}
      \bigl( G_{s,x} \join G_{x,t} \bigr) \cap \bigl( G_{s,y} \join G_{y,t} \bigr)
    \end{equation*}
    contains a directed path from $s$ to $t$ if and only if\/ $G$ contains a
    directed path from $x$ to $y$ or from $y$ to $x$.  Moreover, if\/ $G$
    contains a directed path, say, from $x$ to $y$, then 
    \begin{equation*}
      \bigl( G_{s,x} \join G_{x,t} \bigr) \cap \bigl( G_{s,y} \join G_{y,t} \bigr)
      = G_{s,x} \join G_{x,y} \join G_{y,t}.
    \end{equation*}
  \end{lemma}
  \begin{proof}
    Let $H_x = G_{s,x} \join G_{x,t}$ and $H_y = G_{s,y} \join G_{y,t}$.  If
    $H_x \cap H_y$ contains a directed path $p$ from $s$ to $t$, this path $p$
    passes through $x$ as $p \subseteq H_x$ and it passes through $y$ as $p
    \subseteq H_y$. Taking the subpath of $p$ between $x$ and $y$, we have
    found a directed path between $x$ and $y$. 
    
    Now suppose that there exists some directed path $p$, say, from $x$ to $y$ in
    $G$.  Then $G_{x,y} \neq \emptyset$ and we have $G_{s,x} \join G_{x,y}
    \subseteq G_{s,y}$ and $G_{x,y} \join G_{y,t} \subseteq G_{x,t}$.  Thus
    \begin{equation*}
      G_{s,x} \join G_{x,y} \join G_{y,t} \subseteq H_x \cap H_y.
    \end{equation*}
    In order to show the converse inclusion, let us consider some path $p$ in
    $H_x \cap H_y$ that consists of at most one edge. We then find
    directed paths $q_x$ and $q_y$ from $s$ to $t$ such that $x \in q_x$,  $p
    \subseteq q_x$, $y \in q_y$ and $p \subseteq q_y$.  We distinguish the
    following cases:
    \begin{enumerate}
      \item The vertex $x$ does not precede $t(p)$ on $q_x$. Note that $p$ then lies on a path from $s$ to $x$, i.\,e.\ in $G_{s,x} \subseteq G_{s,x} \join G_{x,y} \join G_{y,t}$.
      \item The vertex $s(p)$ does not precede $y$ on $q_y$. The path $p$ then lies on
        a path from $y$ to $t$ and hence in $G_{y,t} \subseteq G_{s,x} \join
        G_{x,y} \join G_{y,t}$.
      \item The vertex $x$ precedes $t(p)$ on $q_x$ and the vertex $s(p)$ precedes $y$ on~$q_y$.
        Write $q_x = a_x \cdot b_x \cdot p \cdot c_x$ with $a_x$ a path from $s$ to $x$ and $b_x$ a possibly empty path from $x$ to $s(p)$. This is possible as $p$ consists of at most one edge. Similarly, write $q_y = a_y \cdot p \cdot b_y \cdot c_y$ with $a_y$ a path from $s$ to $s(p)$ and $b_y$ a path from $t(p)$ to $y$. 
        The concatenation $b_x \cdot p \cdot b_y$ is a path from $x$ to $y$ that contains $p$ and we conclude that $p \subseteq G_{x,y} \subseteq G_{s,x} \join G_{x,y} \join G_{y,t}$.
    \end{enumerate}
  \end{proof}

\section{Nerves of globular graphs}
\label{sub:nerves_of_globular_graphs}

In this section, we first give a formal definition of the nerve $\nerve(G)$ of
a globular graph as the nerve of a certain partially ordered set. After some
basic examples, we then show that the simplices of $\nerve(G)$ are in bijection
with certain marked subgraphs of $G$. This ultimately leads to a quite
intuitive description of the nerve of a globular graph and to a pictorial
calculus for the action of simplicial operators on $\nerve(G)$.

\subsection*{Definition and examples}
Let $G$ be a globular graph and consider two $st$-paths $p$ and $q$.  Let
us write $p \leq q$ if there exists a glob $\gamma \subseteq G$ and
possibly trivial paths $a$ and $b$ in $G$ such that $p = a \cdot
\dom\gamma \cdot b$ and $q = a \cdot \cod\gamma \cdot b$. We call any
such glob $\gamma$ a \emph{witness} for the relation $p \leq q$.  Observe
that minimal witnesses for a relation $p < q$ are unique.
\begin{lemma}\label{lem:path-partial-order}
 The relation\/ \enquote{\:$\leq$} defines a partial order on the set\/
  $PG$ of\/ $st$-paths in\/ $G$. 
\end{lemma}
\begin{proof}
  We have $p \leq p$ as $p$ itself is a glob.  Suppose $p \leq q$ and $q \leq
  r$ for $st$-paths $p$, $q$ and $r$ in $G$. Choose globs $\gamma,\delta
  \subseteq G$ and decompositions $p = a \cdot \dom \gamma \cdot  b $, $q = a
  \cdot \cod\gamma \cdot b = a' \cdot \dom
  \delta \cdot b'$ and $r = a' \cdot \cod \delta \cdot b'$.  We may assume
  without loss of generality that $a = a'$ and $b = b'$, for we could otherwise
  decompose $a = a_0 \cdot a_1$, $a' = a_0 \cdot a_1'$, $b = b_1 \cdot b_0$ and
  $b' = b_1' \cdot b_0$ and choose $\gamma = a_1 \cup \gamma \cup b_1$ and
  $\delta = a_1' \cup \delta \cup b_1'$ as witnesses for $p \leq q$ and $q \leq
  r$. But if $a = a'$ and $b = b'$, then $\del(\gamma \cup \delta)$ is a glob
  witnessing $p \leq r$.

  Now suppose $p \leq q$ and $q \leq p$. As above, we find witnesses $\gamma$
  and $\delta$ such that $p = a \cdot \dom\gamma \cdot b = a \cdot \cod\delta
  \cdot b$ and $q = a \cdot \cod \gamma \cdot b  = a \cdot \dom \gamma \cdot
  b$. We thus have 
  \begin{equation*}
    \del \gamma = \dom\gamma \cdot \cod\gamma\op = \cod \delta \cdot \dom\delta\op = (\dom\delta \cdot \cod\delta)\op = \del\delta\op
  \end{equation*}
  for the clockwise directed boundaries of $\gamma$ and $\delta$. This implies
  that both $\gamma$ and $\delta$ are degenerate.
\end{proof}
\begin{definition}\label{def:nerve-globular-graph}
  The\/ \emph{nerve $\nerve(G)$} of a globular graph\/ $G$ is the nerve of the
  partially ordered set\/ $(PG,\leq)$. An $n$-simplex $\sigma \in \nerve(G)_n$
  is thus an $(n+1)$-chain 
  \begin{equation*}
   \sigma = (p_0 \leq \dots \leq p_n) 
  \end{equation*}
  of $st$-paths $p_i \in PG$ and the action of a simplicial operator
  $\alpha\colon [m] \to [n]$ on this simplex is determined by $\sigma \alpha =
  (q_0 \leq \dots \leq q_m)$ with $q_i = p_{\alpha(i)}$.
\end{definition}
The reader might object that the purely combinatorial
\autoref{def:nerve-globular-graph} is not satisfactory, for any nerve
$\nerve\colon \cB \to \cAhat$ taking values in a category $\cBhat$ of
presheaves should arise as a composition $\nerve = i^* \comp y$, where $y
\colon \cB \to \cBhat$ denotes the Yoneda embedding and $i^*\colon \cBhat \to
\cAhat$ is restriction along a preferably dense functor $i \colon \cA \to \cB$.
It is in fact easy to show, though, that the nerve of globular graphs arises in
exactly this manner and we will come back to this topic in
\autoref{sec:globular_graphs_as_canonical_arities}.

\begin{example}\label{ex:nerve-bn}
  Consider the graph $B_n$ with two vertices $s$ and $t$ and~$n+1$ distinct edges
  $e_0,\dots,e_n$ between them. The set $PB_n$ of $st$-paths in $B_n$ is the
  set $\{e_0,\dots,e_n\}$. We embed $B_n$ such that for all $i \in
  \{1,\dots,n\}$ there is an interior face $\phi_i$ with oriented boundary
  $\del \phi_i = e_{i-1} \cdot e_i\op$. A glob $\gamma$ in $B_n$ is then
  nothing but a pair $(e_i,e_j)$ of edges with $i \leq j$. The boundary of such
  a $\gamma = (e_i,e_j)$ is given by $\del\gamma = e_i \cdot e_j\op$ and we
  thus have $\dom\gamma = e_i$ and $\cod\gamma = e_j$.  This implies that $e_i
  \leq e_j$ if and only if $i \leq j$. Altogether, we see that $PB_n$ is
  isomorphic to the ordinal $[n]$ and $\nerve(B_n)$ is isomorphic
  to~$\Delta^n$. Note that a nondegenerate $1$-simplex $(p_0 < p_1)$ of
  $\nerve(B_n)$ is contained in the spine of $\nerve(B_n)$ if and only if there
  is a single face of~$B_n$ witnessing the relation $p_0 < p_1$. This
  observation along with our computation of $\nerve(B_n)$ is also
  illustrated in \autoref{fig:b2} for the case $n = 2$. 
\begin{figure}
  \centering
  \begin{subfigure}[b]{0.3\textwidth}
    \centering
    \begin{tikzpicture}
      \node[fill,circle,inner sep=2pt] (s) at (0,0) {};
      \node[fill,circle,inner sep=2pt] (t) at (3,0) {}; 
      \node[left] at (s) {$s$};
      \node[right] at (t) {$t$};
      \draw (s) edge[->,out=45,in=135] node[fill=white,font=\scriptsize] {$e_0$} (t);
      \draw (s) edge[->,out=0,in=180] node[fill=white,font=\scriptsize] {$e_1$} (t);
      \draw (s) edge[->,out=-45,in=-135] node[fill=white,font=\scriptsize] {$e_2$} (t);
  \end{tikzpicture}
    \caption{$B_2$}
  \end{subfigure}
  \hfil
  \begin{subfigure}[b]{0.3\textwidth}
    \centering
    \begin{tikzpicture}
      \node[fill,circle,inner sep=2pt] (e0) at (0,0) {};
      \node[fill,circle,inner sep=2pt] (e1) at (0,1) {}; 
      \node[fill,circle,inner sep=2pt] (e2) at (0,2) {};
      \node[left] at (e0) {$e_0$};
      \node[left] at (e1) {$e_1$};
      \node[left] at (e2) {$e_2$};
      \fill (e0) circle [radius=2pt];
      \fill (e1) circle [radius=2pt];
      \fill (e2) circle [radius=2pt];
      \path[->] (e0) edge coordinate (e0e1)  (e1)
        (e1) edge coordinate (e1e2) (e2)
        (e0) edge[out=45,in=-45] coordinate (e0e2) (e2);
      \node[anchor=east] at (e1e2) {%
        \begin{tikzpicture}[scale=0.3]
          \coordinate (s) at (0,0);
          \coordinate (t) at (3,0);
          \draw (s) edge[-,out=-45,in=-135]  (t);
          \draw (s) edge[-,out=0,in=180]  (t);
        \end{tikzpicture}};
      \node[anchor=west] at (e0e2) {%
        \begin{tikzpicture}[scale=0.3]
          \coordinate (s) at (0,0);
          \coordinate (t) at (3,0);
          \draw (s) edge[-,out=-45,in=-135] (t);
          \draw (s) edge[-,out=45,in=135] (t);
        \end{tikzpicture}};
      \node[anchor=east] at (e0e1) {
        \begin{tikzpicture}[scale=0.3]
          \coordinate (s) at (0,0);
          \coordinate (t) at (3,0);
          \draw (s) edge[-,out=45,in=135]  (t);
          \draw (s) edge[-,out=0,in=180]  (t);
        \end{tikzpicture}};
    \end{tikzpicture}
    \caption{$PB_2$}
  \end{subfigure}
  \hfil
  \begin{subfigure}[b]{0.3\textwidth}
    \centering
    \begin{tikzpicture}
      \node[fill,circle,inner sep=2pt] (e0) at (0,0) {};
      \node[fill,circle,inner sep=2pt] (e2) at (3,0) {}; 
      \node[fill,circle,inner sep=2pt] (e1) at (1.5,1.5) {};
      \node[below] at (e0) {$e_0$};
      \node[below] at (e2) {$e_2$};
      \node[above] at (e1) {$e_1$};
      \path[->] (e0) edge coordinate (e0e2) (e2);
      \path[->] (e0) edge coordinate (e0e1) (e1);
      \path[->] (e1) edge coordinate (e1e2) (e2);
      \node[rectangle,anchor=south east] at (e0e1) {%
        \begin{tikzpicture}[scale=0.3]
          \coordinate (s) at (0,0);
          \coordinate (t) at (3,0);
          \draw (s) edge[-,out=0,in=180] (t);
          \draw (s) edge[-,out=45,in=135] (t);
        \end{tikzpicture}
      };
      \node[rectangle,anchor = south west] at (e1e2) {%
        \begin{tikzpicture}[scale=0.3]
          \coordinate (s) at (0,0);
          \coordinate (t) at (3,0);
          \draw (s) edge[-,out=0,in=180] (t);
          \draw (s) edge[-,out=-45,in=-135] (t);
        \end{tikzpicture}
        };
      \node[rectangle,anchor=north] at (e0e2) {%
        \begin{tikzpicture}[scale=0.3]
          \coordinate (s) at (0,0);
          \coordinate (t) at (3,0);
          \draw (s) edge[-,out=-45,in=-135] (t);
          \draw (s) edge[-,out=45,in=135] (t);
        \end{tikzpicture}
        };
      \node at (barycentric cs:e0=0.5,e1=0.5,e2=0.5) {%
        \begin{tikzpicture}[scale=0.3]
          \coordinate (s) at (0,0);
          \coordinate (t) at (3,0);
          \draw (s) edge[-,out=0,in=180] (t);
          \draw (s) edge[-,out=-45,in=-135] (t);
          \draw (s) edge[-,out=45,in=135] (t);
        \end{tikzpicture}
        };
    \end{tikzpicture}
    \caption{$\nerve(B_2) \iso \Delta^2$}
  \end{subfigure}
  \caption{Computation of $\nerve(B_2)$.}\label{fig:b2}
\end{figure}
\end{example}

\begin{example}\label{ex:nerve-join-bn-bm}
Let us next consider the join $G = B_{n} \join B_{m}$ of two such graphs $B_n$
and $B_m$. We denote the edges of $B_n$ and $B_m$ by $e_0,\dots,e_n$ and
$d_0,\dots,d_m$, respectively.  The set of $st$-paths in $G$ is then given by
\begin{equation*}
  PG = \{ e_{i} \cdot d_k \mid
  0 \leq i \leq n \text{ and } 0 \leq k \leq m\}.
\end{equation*}
Any glob $\gamma$ in $G$ is the join of two globs $(e_i,e_j)$ in $B_n$ and
$(d_k,d_l)$ in $B_m$. We thus have $e_i \cdot d_k \leq e_j \cdot d_l$  if and
only if $i \leq j$ and $k \leq l$. Summing up, we have $PG$ isomorphic to the
product $[n] \times [m]$ and hence $\nerve(G) \iso \Delta^n \times \Delta^m$.
This observation on joins and products actually holds true for all globular
graphs as we discuss in a more general setting in \autoref{lem:nerve-of-join}
below. An explicit computation of $\nerve(B_1 \join B_1)$ along the lines of
the computation of $\nerve(B_2)$ in \autoref{fig:b2} can be seen in
\autoref{fig:b1-times-b1}.

\begin{figure}
  \centering
  \begin{subfigure}[b]{0.8\textwidth}
    \centering
    \begin{tikzpicture}
      \node[fill,circle,inner sep=2pt] (s) at (0,0) {};
      \node[fill,circle,inner sep=2pt] (m) at (3,0) {}; 
      \node[fill,circle,inner sep=2pt] (t) at (6,0) {}; 
      \node[left] at (s) {$s$};
      \node[right] at (t) {$t$};
      \draw[->] (s) edge[out=45,in=135] node[fill=white,font=\scriptsize] {$e_0$} (m);
      \draw[->] (s) edge[out=-45,in=-135] node[fill=white,font=\scriptsize] {$e_1$} (m);
      \draw[->] (m) edge[out=45,in=135] node[fill=white,font=\scriptsize] {$d_0$} (t);
      \draw[->] (m) edge[out=-45,in=-135] node[fill=white,font=\scriptsize] {$d_1$} (t);
  \end{tikzpicture}
  \caption{$ B_1 \join B_1$}
  \end{subfigure}
  \vskip\baselineskip
  \begin{subfigure}[b]{0.49\textwidth}
    \centering
    \begin{tikzpicture}
      \node[fill,circle,inner sep=2pt] (d0e0) at (0,0) {};
      \node[fill,circle,inner sep=2pt] (d1e0) at (-1.5,1.5) {};
      \node[fill,circle,inner sep=2pt] (d0e1) at (1.5,1.5) {};
      \node[fill,circle,inner sep=2pt] (d1e1) at (0,3) {};
      \node[below] at (d0e0) {$e_0 \cdot d_0$};
      \node[left] at (d1e0) {$e_0 \cdot d_1$};

      \node[right] at (d0e1) {$e_1 \cdot d_0$};
      \node[above] at (d1e1) {$e_1 \cdot d_1$};
      \path[->] 
        (d0e0) edge coordinate (d0e0-d0e1) (d0e1)
        (d0e0) edge coordinate (d0e0-d1e0) (d1e0)
        (d1e0) edge coordinate (d1e0-d1e1) (d1e1)
        (d0e1) edge coordinate (d0e1-d1e1) (d1e1)
        (d0e0) edge coordinate (d0e0-d1e1) (d1e1);
      \node[anchor=north east] at (d0e0-d1e0) {%
        \begin{tikzpicture}[scale=0.2]
          \coordinate (s) at (0,0);
          \coordinate (m) at (3,0); 
          \coordinate (t) at (6,0); 
          \draw (s) edge[out=45,in=135] (m);
          \draw (m) edge[out=45,in=135] (t);
          \draw (m) edge[out=-45,in=-135] (t);
        \end{tikzpicture}
      };
      \node[anchor=north west] at (d0e0-d0e1) {%
        \begin{tikzpicture}[scale=0.2]
          \coordinate (s) at (0,0);
          \coordinate (m) at (3,0); 
          \coordinate (t) at (6,0); 
          \draw[] (s) edge[out=45,in=135] (m);
          \draw[] (s) edge[out=-45,in=-135] (m);
          \draw[] (m) edge[out=45,in=135] (t);
        \end{tikzpicture}
      };
      \node[anchor=south west] at (d0e1-d1e1) {%
        \begin{tikzpicture}[scale=0.2]
          \coordinate (s) at (0,0);
          \coordinate (m) at (3,0); 
          \coordinate (t) at (6,0); 
          \draw[] (s) edge[out=-45,in=-135] (m);
          \draw[] (m) edge[out=45,in=135] (t);
          \draw[] (m) edge[out=-45,in=-135] (t);
        \end{tikzpicture}
      };
      \node[anchor=south east] at (d1e0-d1e1) {%
        \begin{tikzpicture}[scale=0.2]
          \coordinate (s) at (0,0);
          \coordinate (m) at (3,0); 
          \coordinate (t) at (6,0); 
          \draw[] (s) edge[out=45,in=135] (m);
          \draw[] (s) edge[out=-45,in=-135] (m);
          \draw[] (m) edge[out=-45,in=-135] (t);
        \end{tikzpicture}
      };
      \node[fill=white,inner sep=1pt] at (d0e0-d1e1) {%
        \begin{tikzpicture}[scale=0.2]
          \coordinate (s) at (0,0);
          \coordinate (m) at (3,0); 
          \coordinate (t) at (6,0); 
          \draw[] (s) edge[out=45,in=135] (m);
          \draw[] (s) edge[out=-45,in=-135] (m);
          \draw[] (m) edge[out=-45,in=-135] (t);
          \draw[] (m) edge[out=45,in=135] (t);
        \end{tikzpicture}
      };
    \end{tikzpicture}
    \caption{$P(B_1 \join B_1)$}
  \end{subfigure}
  \hfil
  \begin{subfigure}[b]{0.49\textwidth}
    \centering
    \begin{tikzpicture}
      \node[fill,circle,inner sep=2pt] (d0e0) at (0,0) {};
      \node[fill,circle,inner sep=2pt] (d1e0) at (-1.5,1.5) {};
      \node[fill,circle,inner sep=2pt] (d0e1) at (1.5,1.5) {};
      \node[fill,circle,inner sep=2pt] (d1e1) at (0,3) {};
      \node[below] at (d0e0) {$e_0 \cdot d_0$};
      \node[left] at (d1e0) {$e_0 \cdot d_1$};

      \node[right] at (d0e1) {$e_1 \cdot d_0$};
      \node[above] at (d1e1) {$e_1 \cdot d_1$};
      \path[->] 
        (d0e0) edge coordinate (d0e0-d0e1) (d0e1)
        (d0e0) edge coordinate (d0e0-d1e0) (d1e0)
        (d1e0) edge coordinate (d1e0-d1e1) (d1e1)
        (d0e1) edge coordinate (d0e1-d1e1) (d1e1)
        (d0e0) edge coordinate[pos=0.75] (d0e0-d1e1) (d1e1);
      \node[anchor=north east] at (d0e0-d1e0) {%
        \begin{tikzpicture}[scale=0.2]
          \coordinate (s) at (0,0);
          \coordinate (m) at (3,0); 
          \coordinate (t) at (6,0); 
          \draw (s) edge[out=45,in=135] (m);
          \draw (m) edge[out=45,in=135] (t);
          \draw (m) edge[out=-45,in=-135] (t);
        \end{tikzpicture}
      };
      \node[anchor=north west] at (d0e0-d0e1) {%
        \begin{tikzpicture}[scale=0.2]
          \coordinate (s) at (0,0);
          \coordinate (m) at (3,0); 
          \coordinate (t) at (6,0); 
          \draw[] (s) edge[out=45,in=135] (m);
          \draw[] (s) edge[out=-45,in=-135] (m);
          \draw[] (m) edge[out=45,in=135] (t);
        \end{tikzpicture}
      };
      \node[anchor=south west] at (d0e1-d1e1) {%
        \begin{tikzpicture}[scale=0.2]
          \coordinate (s) at (0,0);
          \coordinate (m) at (3,0); 
          \coordinate (t) at (6,0); 
          \draw[] (s) edge[out=-45,in=-135] (m);
          \draw[] (m) edge[out=45,in=135] (t);
          \draw[] (m) edge[out=-45,in=-135] (t);
        \end{tikzpicture}
      };
      \node[anchor=south east] at (d1e0-d1e1) {%
        \begin{tikzpicture}[scale=0.2]
          \coordinate (s) at (0,0);
          \coordinate (m) at (3,0); 
          \coordinate (t) at (6,0); 
          \draw[] (s) edge[out=45,in=135] (m);
          \draw[] (s) edge[out=-45,in=-135] (m);
          \draw[] (m) edge[out=-45,in=-135] (t);
        \end{tikzpicture}
      };
      \node[fill=white,inner sep=1pt] at (d0e0-d1e1) {%
        \begin{tikzpicture}[scale=0.175]
          \coordinate (s) at (0,0);
          \coordinate (m) at (3,0); 
          \coordinate (t) at (6,0); 
          \draw[] (s) edge[out=45,in=135] (m);
          \draw[] (s) edge[out=-45,in=-135] (m);
          \draw[] (m) edge[out=-45,in=-135] (t);
          \draw[] (m) edge[out=45,in=135] (t);
        \end{tikzpicture}
      };
      \node[] at (barycentric cs:d0e0=0.5,d1e0=0.8,d1e1=0.5) {%
        \begin{tikzpicture}[scale=0.175]
          \coordinate (s) at (0,0);
          \coordinate (m) at (3,0); 
          \coordinate (t) at (6,0); 
          \draw[] (s) edge[out=45,in=135] (m);
          \draw[] (s) edge[out=-45,in=-135] (m);
          \draw[] (m) edge[out=-45,in=-135] (t);
          \draw[] (m) edge[out=45,in=135] (t);
          \node[font=\tiny] at ($(s) ! 0.5 ! (m)$) {$2$};
          \node[font=\tiny] at ($(m) ! 0.5 ! (t)$) {$1$};
        \end{tikzpicture}
      };
      \node[] at (barycentric cs:d0e0=0.5,d0e1=0.8,d1e1=0.5) {%
        \begin{tikzpicture}[scale=0.175]
          \coordinate (s) at (0,0);
          \coordinate (m) at (3,0); 
          \coordinate (t) at (6,0); 
          \draw[] (s) edge[out=45,in=135] (m);
          \draw[] (s) edge[out=-45,in=-135] (m);
          \draw[] (m) edge[out=-45,in=-135] (t);
          \draw[] (m) edge[out=45,in=135] (t);
          \node[font=\tiny] at ($(s) ! 0.5 ! (m)$) {$1$};
          \node[font=\tiny] at ($(m) ! 0.5 ! (t)$) {$2$};
        \end{tikzpicture}
      };
    \end{tikzpicture}
    \caption{$\nerve(B_1 \join B_1) \iso \Delta^1 \times \Delta^1$}
  \end{subfigure}
  \caption{Computation of $\nerve(B_1 \join B_1)$.}\label{fig:b1-times-b1}
\end{figure}
\end{example}

\subsection*{A pictorial representation of $\nerve(G)$}
Even though the definition of $\nerve(G)$ in \autoref{def:nerve-globular-graph}
is concise and sufficient for most technical purposes, it leaves a lot to be
desired from an intuitional point of view. Therefore, in this paragraph, we
give another description of $\nerve(G)$ in terms of certain marked globular
subgraphs of $G$ such as those that we already used in the computation of $\nerve(B_1 \join B_1)$ in \autoref{fig:b1-times-b1}. 
This description also plays a role in our proof of \autoref{thm:labeling} in \autoref{chap:labeling}.

Let us fix a globular graph $G$ with source $s$ and target $t$ and consider an
$n$-simplex $\sigma = (p_0 \leq \dots \leq p_n) \in \nerve(G)$. Observe that
\autoref{prop:ps:power} implies that $P = \bigcup_{i = 0}^n p_i$ is a wide
globular subgraph of $G$, for it is acyclic and contains both a source and
target.  The following lemma is crucial to the description of $\nerve(G)$ that
we are about to give: 
\begin{lemma}\label{lem:exists-unique-witness}
  Let $\sigma = (p_0 \leq \dots \leq p_n)$ be an $n$-simplex in $\nerve(G)$ and
  let $\phi$ be an interior face of $P = \bigcup_{i = 0}^n p_i$. There then
  exists some $i \in \{1,\dots,n\}$ such that $\cod \phi \subseteq p_i$.
  Moreover, $\dom \phi \subseteq p_{i-1}$ if and only if\/ $i$ is chosen minimal.
\end{lemma}
\begin{proof}
  The proof is by induction on $n$ and the claim is obvious for $n=0$ and
  $n=1$.  Now let $n \geq 2$. The interior faces of $P$ are the bounded
  connected components of $(\dR^2 \smallsetminus P') \smallsetminus p_n$, where
  $P' = \bigcup_{i = 0}^{n-1} p_i$. Observe that $p_n$ and $P'$ intersect only
  in $p_{n-1} \subseteq \del P'$ since $p_{n-1} \leq p_n$. This implies in
  particular that $p_n$ does not intersect any interior face of $P'$ and that
  any interior face of $P$ that is not already an interior face of $P'$ is
  necessarily bounded by subpaths of $p_{n-1}$ and $p_n$.  The claim follows.
  
\end{proof}

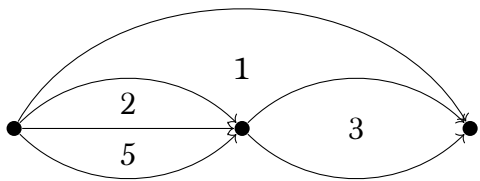
\begin{figure}
   \centering
          \begin{tikzpicture}
            \node[fill,circle,inner sep=2pt] (s) at (0,0) {};
            \node[fill,circle,inner sep=2pt] (m) at (3,0) {};
            \node[fill,circle,inner sep=2pt] (t) at (6,0) {}; 
            \draw (s) edge[->,out=45,in=135] coordinate(1) (m);
            \draw (s) edge[->,out=-45,in=-135] coordinate(3) (m);
            \draw (s) edge[->] coordinate(2) (m);
            \draw (m) edge[->,out=45,in=135] (t);
            \draw (m) edge[->,out=-45,in=-135] (t);
            \draw (s) edge[->,out=60,in=120] coordinate[pos=0.5] (oben) (t);
            \node at ($ (1) ! 0.5 ! (2) $) {$2$};
            \node at ($ (2) ! 0.5 ! (3) $) {$5$};
            \node at ($ (m) ! 0.5 ! (oben) $) {$1$};
            \node at ($ (m) ! 0.5 ! (t) $) {$3$};
        \end{tikzpicture}
          \caption{A $5$-marked subgraph $(P,\lambda)$.}\label{fig:marked-subgraph}
\end{figure}

\autoref{lem:exists-unique-witness} implies in particular that we may associate
with a given $n$-simplex $\sigma =(p_0 \leq \dots \leq p_n)$ of $\nerve(G)$ a
tuple $(P_\sigma,\lambda_\sigma)$, where $P_\sigma = \bigcup_{i = 0}^n p_i$ is
a wide globular subgraph of $G$ and where $\lambda_\sigma \colon \Phi(P_\sigma)
\to \{1,\dots,n\}$ is the map that associates with any interior face $\phi \in
\Phi(P_\sigma)$ the unique $i \in \{1,\dots,n\}$ such that $\phi \subseteq
p_{i-1} \cup p_i$.  This observation motivates the following definition:
\begin{definition}\label{def:marked}
  An \emph{$n$-marked subgraph} $(P,\lambda)$ of a globular graph $G$
  consists of a wide globular subgraph $P \subseteq G$ and a map
  $\lambda\colon \Phi(P) \to \{1,\dots,n\}$ from the set of its interior faces
  into $\{1,\dots,n\}$.
\end{definition}
We picture $n$-marked subgraphs by drawing the graph $P$ and labelling
the interior faces $\phi$ of $P$ with $\lambda(\phi)$, see \autoref{fig:marked-subgraph}. The
reader should note that this pictorial representation of $(P,\lambda)$ does not
determine the number $n$, though.

It is clear that not all $n$-marked globular subgraphs $(P,\lambda)$ arise from
$n$-simplices $\sigma$ of $\nerve(G)$. One condition
that all the maps $\lambda_\sigma$ satisfy, though, is the following:
\begin{lemma}
  Let $(P_\sigma,\lambda_\sigma)$ be the $n$-marked subgraph associated with
  some $n$-simplex $\sigma \in \nerve(G)$. Suppose that $\phi$ and $\psi$ are
  two faces of\/ $P_\sigma$ such that there exists some edge $e$ of\/ $P_\sigma$
  such that $e \subseteq \cod \phi$ and $e \subseteq \dom \psi$. Then
  $\lambda_\sigma(\phi) < \lambda_\sigma(\psi)$.
\end{lemma}
\begin{proof}
  We know that $j = \lambda_\sigma(\psi)$ is minimal with $\cod \psi \subseteq
  p_j$. Any $p_k$ that contains the edge $e$ cannot contain $\cod \psi$ simply
  because $e \subseteq \dom \psi$. We thus have $j > k$ for all $k$ with $e
  \subseteq p_k$ and this implies in particular $j > \lambda_\sigma(\phi)$ as $e \subseteq \cod(\phi)$ by assumption.
\end{proof}
\begin{definition}\label{def:admissible}
  An $n$-marked subgraph\/ $(P,\lambda)$ is \emph{admissible} if\/
  $\lambda(\phi) < \lambda(\psi)$ whenever $\cod \phi \cap \dom \psi$ contains
  an edge.
\end{definition}
\begin{remark}\label{rem:left-right-of-path}
  Let $(P,\lambda)$ be admissible and let $p$ be some $st$-path in $P$. Any face $\phi$ on the left of $p$ then satisfies 
  \begin{equation*}
    \lambda(\phi) \leq \max\{ \lambda(\xi) \mid \cod\xi \cap p \text{ contains an edge} \}
  \end{equation*}
  and any face $\psi$ on the right of $p$ satisfies
  \begin{equation*}
    \lambda(\psi) \geq \min\{ \lambda(\xi) \mid \dom\xi \cap p \text{ contains an edge} \}.
  \end{equation*}
  Indeed, if $\phi = \phi_0$ is an interior face on the left of $p$ such that
  $\cod \phi_0 \cap p$ does not contain an edge, we find an interior face
  $\phi_1$ to the left of $p$ such that $\cod \phi_0 \cap \dom \phi_1$ contains
  an edge, i.\,e.\ $\lambda(\phi_0) < \lambda(\phi_1)$. If $\cod \phi_1 \cap p$
  contains an edge, we are done. Otherwise, we can continue in this manner to
  obtain a sequence of faces $\phi_0,\dots,\phi_k$ that are on the left of $p$
  and satisfy $\lambda(\phi_0) < \dots < \lambda(\phi_k)$. This procedure
  terminates after a finite number of steps, for $P$ is a globular graph.  We
  thus find a face $\phi_k$ such that $\lambda(\phi_0) < \lambda(\phi_k)$ and
  $\cod \phi_k \cap p$ contains an edge.
\end{remark}
\needspace{4\baselineskip}
\begin{lemma}\label{lem:paths-exist-for-admissible}
  Let\/ $\lambda \colon P \to \{1,\dots,n\}$ be admissible.  For each\/ $k \in
  \{0,\dots,n\}$ there exists a unique $st$-path $p_k \subseteq P$ such that
  all the interior faces $\phi$ on the left of\/ $p_k$ satisfy $\lambda(\phi) \leq  k$ and all the interior faces $\phi$ on the right of\/ $p_k$ satisfy
  $\lambda(\phi) > k$.  
\end{lemma}
\begin{proof}
  We obviously have $p_0 = \dom P$.  Given such a path $p_k$, we can construct
  $p_{k+1}$ by replacing each subpath of $p_k$ of the form $\dom \phi$ for some
  interior face $\phi$ with $\lambda(\phi) = k+1$ by $\cod \phi$. This is
  indeed the sought-for path $p_k$ by \autoref{rem:left-right-of-path}.  Let us
  now suppose that we have two such paths $p_k \neq q_k$.  As both $p_k$ and
  $q_k$ start in $s$, there is then some vertex $u$ and edges $e \neq d$ with
  source $u$ such that $e \subseteq p_k$ and $d \subseteq q_k$.  We know from
  \autoref{lem:cyclic-order} that the edges with source $u$ appear
  consecutively in the cyclic order of edges around $u$ and after possibly changing the role of $p_k$ and $q_k$ we therefore find
  interior faces $\phi_1,\dots,\phi_r$ such that $e \subseteq \dom \phi_1$, $d
  \subseteq \cod \phi_r$ and such that $\cod \phi \cap \dom \phi_{i+1}$
  contains an edge for all $1 \leq i < r$, see \autoref{fig:path-exists}. We then have $k <
  \lambda(\phi_1) \leq \lambda(\phi_r) \leq k$ by admissibility of $\lambda$
  and the definitions of $p_k$ and $q_k$. This is a contradiction and we conclude $p_k = q_k$.

  \begin{figure}
          \centering
        \begin{tikzpicture}
            \node[fill,circle,inner sep=2pt] (u) at (0,0) {};
            \node[fill,circle,inner sep=2pt] (te) at (90:2) {};
            \node[fill,circle,inner sep=2pt] (dummy1) at (50:2) {};
            \node[fill,circle,inner sep=2pt] (dummy2) at (10:2) {};
            \node[fill,circle,inner sep=2pt] (dummy3) at (-20:2) {};
            \node[fill,circle,inner sep=2pt] (td) at (-60:2) {};
            \node[left] at (u) {$u$};
            \path (u) edge[->] node[left] {$e$} coordinate[pos=0.8] (e) (te) 
                  (u) edge[->] node[below left] {$d$} coordinate[pos=0.8] (d) (td)
                  (u) edge[->] coordinate[pos=0.8] (e2) (dummy1)
                  (u) edge[->] coordinate[pos=0.8] (e3) (dummy2)
                  (u) edge[->] coordinate[pos=0.8] (e4) (dummy3)
                  ;
            \node at ($ (e) ! 0.5 ! (e2) $) {$\phi_1$};
            \node at ($ (e2) ! 0.5 ! (e3) $) {$\phi_2$};
            \node at ($ (e4) ! 0.5 ! (d) $) {$\phi_r$};
            \node at ($ (e3) ! 0.5 ! (e4) $) {$\vdots$};
        \end{tikzpicture}
          \caption{An illustration of an argument in the proof of \autoref{lem:paths-exist-for-admissible}.}\label{fig:path-exists}
  \end{figure}
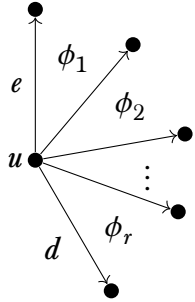

\end{proof}

\begin{remark}\label{rem:pi-characterisation}
  One consequence of \autoref{lem:paths-exist-for-admissible} is that the path
  $p_i$ occuring in some simplex $\sigma = (p_0 \leq \dots \leq p_n)$ of
  $\nerve(G)$ can be characterised as the unique $st$-path $p$ in $P_\sigma$
  with the property that all interior faces $\phi$ of $P_\sigma$ with
  $\lambda_\sigma(\phi) \leq i$ are on the left hand side of $p$.
\end{remark}

The following proposition gives us the promised description of the simplices of
$\nerve(G)$ in terms of $n$-marked graphs:
\begin{proposition}\label{prop:nmarked-is-simplex}
  Let $G$ be a globular graph.  The assignment $\sigma \mapsto
  (P_\sigma,\lambda_\sigma)$ defines a bijection between the set of
  $n$-simplices $\sigma$ of\/ $\nerve(G)$ and the set of $n$-marked subgraphs of\/
  $G$.
\end{proposition}
\begin{proof}
Given some admissible $n$-marked subgraph $(P,\lambda)$ of $G$, we obtain by
  \autoref{lem:paths-exist-for-admissible} a set of $st$-paths $p_0, \dots,
  p_n$ in the wide subgraph $P$ of $G$ such that $p_i$ has precisely the
  interior faces of $P$ with $\lambda(p_i) \leq i$ on its left. The graph
  $p_{i-1} \cup p_i$ is a glob witnessing $p_{i-1} \leq p_i$ and we thus obtain
  an $n$-simplex $\sigma = \sigma(P,\lambda) = (p_0 \leq \dots \leq p_n)$ of
  $\nerve(G)$. We obviously have $P = P_\sigma$. Furthermore, the interior
  faces $\phi$ of $P$ with $\lambda(\phi) = i$ are precisely those interior
  faces of $P$ that lie on the right of $p_{i-1}$ and on the left of $p_i$,
  i.\,e.\ in $p_{i-1} \cup p_i$. It now follows that $\lambda =
  \lambda_\sigma$. 

  Conversely, given a simplex $\sigma$ in $\nerve(G)$, it follows immediately
  from \autoref{rem:pi-characterisation} that the simplex
  $\sigma(P_\sigma,\lambda_\sigma)$ is $\sigma$ again. 
\end{proof}

It remains to desribe the action of simplicial operators in terms of $n$-marked graphs. This is achieved in the following proposition.
\begin{proposition}\label{prop:marked-simplicial}
  Let $\alpha \colon [m] \to [n]$ be a simplicial operator and let
  $(P_\sigma,\lambda_\sigma)$ be the tuple associated with some $n$-simplex
  $\sigma \in \nerve(G)$. The tuple $(P_{\sigma\alpha},\lambda_{\sigma\alpha})$
  associated with $\sigma\alpha$ can then be obtained from
  $(P_\sigma,\lambda_\sigma)$ by the following steps:
  \begin{enumerate}
    \item Remove the edges and interior vertices of all paths $\dom(\phi)$, where $\phi$ is some interior face of\/ $P_\sigma$ with $\lambda_\sigma(\phi) \leq \alpha(0)$.
    \item Remove the edges and interior vertices of all paths $\cod(\phi)$, where $\phi$ is some interior face of\/ $P_\sigma$ with $\lambda_\sigma(\phi) > \alpha(m)$.
    \item Define $\what{\lambda}$ on the remaining graph by
      \begin{equation*}
        \what{\lambda}(\phi) = \min\bigl\{\, k \in \{1,\dots,m\} \mid \alpha(k) \geq \lambda_\sigma(\phi) \bigr\}.
      \end{equation*}
    \item In each maximal globular subgraph $Q \subseteq P$ on which
      $\what{\lambda}$ takes some constant value $c$, remove all edges and
      vertices of\/ $Q$ that are not incident with the exterior face of\/ $Q$ and
      define $\lambda_{\sigma\alpha}(\phi) = c$ on the resulting interior
      faces.   
  \end{enumerate}
\end{proposition}
\begin{proof}
  Let us write $\sigma = (p_0 \leq \dots \leq p_n)$.  The simplex
  $\sigma\alpha$ is then given by $\sigma\alpha = (q_0 \leq \dots \leq q_m)$
  with $q_j = p_{\alpha(j)}$.  Observe that $Q = \bigcup q_j$ is a subgraph of
  $P = \bigcup p_i$. We fix some embedding of $P$ in $\dR^2$ and consider the
  induced embedding of $Q$.  The closure $\overline{\psi}$ of any face $\psi$
  of $Q$ can be written as a union of the closures of a certain set of faces of
  $P$ and each face of $P$ occurs in at most one such union.  Now consider a
  path $q_j = p_{\alpha(j)}$ occuring in $\sigma\alpha$. By
  \autoref{rem:pi-characterisation}, it is the unique $st$-path in $P$ such
  that an interior face $\phi \in \Phi(P)$ lies to the left of $p_{\alpha(j)}$ if and only if $\lambda_\sigma(\phi) \leq \alpha(j)$. It follows that the interior faces $\phi$ in $P$ between
  $q_{j-1}$ and $q_j$ are precisely those interior faces of $P$ with
  $\alpha(j-1) < \lambda_\sigma(\phi) \leq \alpha(j)$. This implies that all
  interior faces $\phi$ of $P$ with $\lambda_\sigma(\phi) \leq \alpha(0)$ or
  $\lambda_\sigma(\phi) > \alpha(m)$ are not contained in any interior face
  $\psi$ of $Q$ and justifies the correctness of the first two steps. Moreover,
  the closure of any interior face $\psi$ of $q_{j-1} \cup q_j$ is the union of
  the closure of interior faces $\phi$ of $P$ with $\alpha(j-1) <
  \lambda_\sigma(\phi) \leq \alpha(j)$. But for any of these faces $\phi$ we certainly have
  \begin{equation*}
    \lambda_{\sigma\alpha}(\psi) = j = \min\bigl\{\, k \in \{1,\dots,m\} \mid \alpha(k) \geq \lambda_\sigma(\phi) \bigr\}
  \end{equation*}
  and this justifies the remaining two steps.

\end{proof}

\begin{example}
  Let us compute the action of $d_1, d_2 \colon [2] \to [3]$ on the admissible
  $3$-marked globular subgraph $(P,\lambda)$ shown in
  \autoref{subfig:face-computation-1}.  Let us first consider the action of
  $d_1$.  As $d_1(0) = 0$ and $d_1(2) = 3$, the first two steps of the
  procedure given in \autoref{prop:marked-simplicial} do not change the graph
  $P$. The values of the map $\what{\lambda}$ of step 3 are easily computed and
  can be seen in \autoref{subfig:d1-whatlambda}. There is then one edge between
  two faces $\phi$ and $\psi$ of $P$ with $1 = \what{\lambda}(\phi) =
  \what{\lambda}(\psi)$. Removing this edge then results in $d_1(P,\lambda)$ as
  can be seen in \autoref{subfig:d1-removal} and \ref{subfig:d1-plambda}.
  
  Let us now consider the action of $d_2$. Again, there are no edges to be
  removed in the first two steps and the map $\what{\lambda}$ can easily seen
  to be the one in \autoref{subfig:d2-whatlambda}. It turns out that
  $(P,\what{\lambda})$ is already the result of the algorithm given in
  \autoref{prop:marked-simplicial}.
  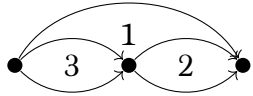
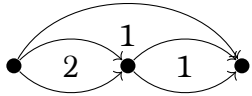
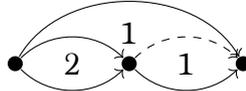
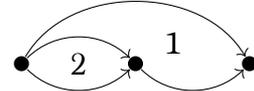
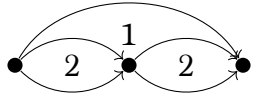
\begin{figure}[p]
        \centering
        \begin{subfigure}[b]{0.8\textwidth}
          \centering
          \begin{tikzpicture}
            \node[fill,circle,inner sep=2pt] (s) at (0,0) {};
            \node[fill,circle,inner sep=2pt] (m) at (1.5,0) {};
            \node[fill,circle,inner sep=2pt] (t) at (3,0) {}; 
            \draw (s) edge[->,out=45,in=135] (m);
            \draw (s) edge[->,out=-45,in=-135]  (m);
            \draw (m) edge[->,out=45,in=135] (t);
            \draw (m) edge[->,out=-45,in=-135] (t);
            \draw (s) edge[->,out=60,in=120] coordinate[pos=0.5] (oben) (t);
            \node at ($ (s) ! 0.5 ! (m) $) {$3$};
            \node at ($ (m) ! 0.5 ! (t) $) {$2$};
            \node at ($ (m) ! 0.5 ! (oben) $) {$1$};
        \end{tikzpicture}
          \caption{An admissible $3$-marked subgraph $(P,\lambda)$.}\label{subfig:face-computation-1}
        \end{subfigure}
        \vskip2ex
        \hfil
        \begin{subfigure}[b]{0.3\textwidth}
          \centering
          \begin{tikzpicture}
            \node[fill,circle,inner sep=2pt] (s) at (0,0) {};
            \node[fill,circle,inner sep=2pt] (m) at (1.5,0) {};
            \node[fill,circle,inner sep=2pt] (t) at (3,0) {}; 
            \draw (s) edge[->,out=45,in=135] (m);
            \draw (s) edge[->,out=-45,in=-135]  (m);
            \draw (m) edge[->,out=45,in=135] (t);
            \draw (m) edge[->,out=-45,in=-135] (t);
            \draw (s) edge[->,out=60,in=120] coordinate[pos=0.5] (oben) (t);
            \node at ($ (s) ! 0.5 ! (m) $) {$2$};
            \node at ($ (m) ! 0.5 ! (t) $) {$1$};
            \node at ($ (m) ! 0.5 ! (oben) $) {$1$};
        \end{tikzpicture}
          \caption{The map $\what{\lambda}$ in the computation of $d_1 (P,\lambda)$.}\label{subfig:d1-whatlambda}
        \end{subfigure}
        \hfil
        \begin{subfigure}[b]{0.3\textwidth}
          \centering
          \begin{tikzpicture}
            \node[fill,circle,inner sep=2pt] (s) at (0,0) {};
            \node[fill,circle,inner sep=2pt] (m) at (1.5,0) {};
            \node[fill,circle,inner sep=2pt] (t) at (3,0) {}; 
            \draw (s) edge[->,out=45,in=135] (m);
            \draw (s) edge[->,out=-45,in=-135]  (m);
            \draw (m) edge[dashed,->,out=45,in=135] (t);
            \draw (m) edge[->,out=-45,in=-135] (t);
            \draw (s) edge[->,out=60,in=120] coordinate[pos=0.5] (oben) (t);
            \node at ($ (s) ! 0.5 ! (m) $) {$2$};
            \node at ($ (m) ! 0.5 ! (t) $) {$1$};
            \node at ($ (m) ! 0.5 ! (oben) $) {$1$};
        \end{tikzpicture}
          \caption{Removal of the edges of $P$ that are not present in $d_1(P,\lambda)$.}\label{subfig:d1-removal}
        \end{subfigure}
        \hfil
        \begin{subfigure}[b]{0.3\textwidth}
          \centering
          \begin{tikzpicture}
            \node[fill,circle,inner sep=2pt] (s) at (0,0) {};
            \node[fill,circle,inner sep=2pt] (m) at (1.5,0) {};
            \node[fill,circle,inner sep=2pt] (t) at (3,0) {}; 
            \draw (s) edge[->,out=45,in=135] (m);
            \draw (s) edge[->,out=-45,in=-135]  (m);
            \draw (m) edge[->,out=-45,in=-135] (t);
            \draw (s) edge[->,out=60,in=120] coordinate[pos=0.5] (oben) (t);
            \node at ($ (s) ! 0.5 ! (m) $) {$2$};
            \node at (barycentric cs:m=0.5,t=0.5,oben=0.5) {$1$};
        \end{tikzpicture}
          \caption{The $2$-marked subgraph $d_1(P,\lambda)$.\newline\strut}\label{subfig:d1-plambda}
        \end{subfigure}
        \vskip2ex
        \hfil
        \begin{subfigure}[b]{0.8\textwidth}
          \centering
          \begin{tikzpicture}
            \node[fill,circle,inner sep=2pt] (s) at (0,0) {};
            \node[fill,circle,inner sep=2pt] (m) at (1.5,0) {};
            \node[fill,circle,inner sep=2pt] (t) at (3,0) {}; 
            \draw (s) edge[->,out=45,in=135] (m);
            \draw (s) edge[->,out=-45,in=-135]  (m);
            \draw (m) edge[->,out=45,in=135] (t);
            \draw (m) edge[->,out=-45,in=-135] (t);
            \draw (s) edge[->,out=60,in=120] coordinate[pos=0.5] (oben) (t);
            \node at ($ (s) ! 0.5 ! (m) $) {$2$};
            \node at ($ (m) ! 0.5 ! (t) $) {$2$};
            \node at ($ (m) ! 0.5 ! (oben) $) {$1$};
        \end{tikzpicture}
          \caption{The map $\what{\lambda}$ in the computation of $d_2 (P,\lambda)$.}\label{subfig:d2-whatlambda}
        \end{subfigure}
       \caption{Examples of the procedure given in
       \autoref{prop:marked-simplicial} for the computation of the action of
       simplicial operators on $n$-marked subgraphs.}
      \end{figure}
\end{example}

\section{Globular graphs and \texorpdfstring{$ 2 $}{2}-computads}
\label{sec:globular_graphs_as_canonical_arities}

It is well-known that the category $\cat$ of small categories is monadic over
the category $\graph$ of graphs. The monad $T_1$ for categories sends a
graph~$G$ to the graph $T_1G$ having the same vertices but possibly empty directed paths
$p$ as edges between $s(p)$ and $t(p)$. The monad $T_1$ is in fact strongly
cartesian, meaning that $\Shat \colon \graph \to \graph/T_1\done$ has a right
adjoint and that all the naturality squares of $\mu \colon T_1^2 \to T_1$ and
$\eta \colon 1 \to T_1$ are cartesian. Similar results hold true for higher
categories if one is willing to substitute $n$-globular sets
for graphs. In the case of $2$-categories this was first observed in
\cite{street;limits-indexed-by-category-valued-2-functors} and has since been
studied in various contexts, see e.\,g.\
\cite{batanin;computads-for-finitary-monads-on-globular-sets} for a
construction of computads for a general finitary monad on globular sets or
\cite{metayer;strict-omega-categories-are-monadic-over-polygraphs} for a 
proof of the monadicity of strict $\omega$-categories over computads.  In this
paragraph, following \cite{leinster;higher-operads-higher-categories}, we
sketch the relation between Street's $2$-computads, $2$-categories and globular
graphs.

\subsection*{\texorpdfstring{$2$}{2}-globular sets and \texorpdfstring{$2$}{2}-categories}
\label{sub:2_globular_sets_and_2_categories}
Consider the category $\dG_2$ generated by the graph
\begin{equation*}
  \begin{tikzpicture}[diagram]
    \matrix[objects] {
      |(0)| \bullet \& |(1)| \bullet \& |(2)| \bullet  \\
    };
    \path[maps,->] 
      (0) edge[eabove] node[above] {$\sigma_0$} (1)
      (0) edge[ebelow] node[below] {$\tau_0$} (1)
      (1) edge[eabove] node[above] {$\sigma_1$} (2)
      (1) edge[ebelow] node[below] {$\tau_1$} (2)
    ;
  \end{tikzpicture}
\end{equation*}
subject to relations $\sigma_1 \sigma_0 = \tau_1 \sigma_0$ and $\sigma_1\tau_0
= \tau_1\tau_0$.  
A \emph{$2$-globular set} is a presheaf $X \in \dGhat_2$ on $\dGhat_2$, i.\,e.\ a collection of three sets $X_0$,
$X_1$ and $X_2$ together with maps
\begin{equation*}
  \begin{tikzpicture}[diagram]
    \matrix[objects] {
      |(0)| X_0 \& |(1)|X_1 \& |(2)| X_2\\
    };
    \path[maps,->] 
      (1) edge[eabove] node[above] {$s_0$} (0)
      (1) edge[ebelow] node[below] {$t_0$} (0)
      (2) edge[eabove] node[above] {$s_1$} (1)
      (2) edge[ebelow] node[below] {$t_1$} (1)
    ;
  \end{tikzpicture}
\end{equation*}
such that $s_0s_1 = s_0t_1$ and $t_0s_1 = t_0t_1$.  

We will have the occasion to use a slightly different description of $\dGhat_2$
in terms of $2$-graphs. A $2$-graph $G$ consists of a set $V(G)$ of vertices
together with graphs $G_{x,y}$ for each pair $x,y \in V(G)$ of vertices.  A map
$G \to H$ of $2$-graphs is a map $f\colon V(G) \to V(H)$ together with
compatible maps $G_{x,y} \to H_{fx,fy}$ of graphs.  

Given such a $2$-graph $G$, we associate with it the globular set 
\begin{equation*}
  X_G = \left( 
  \begin{tikzpicture}[diagram]
    \matrix[objects] {
      |(v)| V(G) \& |(e)| \displaystyle\coprod_{x,y} (G_{x,y})_0 \& |(ee)| \displaystyle\coprod_{x,y} (G_{x,y})_1 \\
    };
    \path[maps,->] 
      (e) edge[eabove] node[above] {$s_0$} (v)
      (e) edge[ebelow] node[below] {$t_0$} (v)
      (ee) edge[eabove] node[above] {$s_1$} (e)
      (ee) edge[ebelow] node[below] {$t_1$} (e)
    ;
  \end{tikzpicture}
  \right),
\end{equation*}
where $s_0(u) = x$ and $t_0(u) = y$ for all $u \in (G_{x,y})_0$ and where the
maps $s_1$ and~$t_1$ are induced by the source and target maps of the graphs
$G_{x,y}$. The assignment $G \mapsto X_G$ can easily be seen to be functorial
in $G$.

Conversely, if $X$ is a $2$-globular set, we can associate with it a $2$-graph
as follows: For each pair $x,y \in X_0$ take the pullbacks
\begin{equation*}
  \begin{tikzpicture}[diagram]
    \matrix[objects] {
         |(v)|    (G_{x,y})_0        \& |(xy)| \done \\
       |(x1)| X_1 \& |(x0x0)| X_0 \times X_0 \\
    };
    \path[maps,->] 
      (x1) edge node[below] {$(s_0,t_0)$} (x0x0)
      (xy) edge node[right] {$(x,y)$} (x0x0)
      (v) edge (xy)
      (v) edge (x1)
    ;
  \end{tikzpicture}
  \quad\text{and}\quad
  \begin{tikzpicture}[diagram]
    \matrix[objects,wider] {
      |(e)| (G_{x,y})_1   \& |(xy)| \done \\
      |(x2)| X_2  \& |(x0x0)| X_0 \times X_0 \smash{.} \\
    };
    \path[maps,->] 
      (xy) edge node[right] {$(x,y)$} (x0x0)
  (x2) edge node[below] {$(s_0,t_0) \comp s_1$} (x0x0)
      (e) edge (x2)
      (e) edge (xy)
    ;
  \end{tikzpicture}
\end{equation*}
The bottom map in the right hand square satisfies $(s_0,t_0) \comp s_1 =
(s_0,t_0) \comp t_1$ as $X$ is a $2$-globular set by assumption. We therefore
get two factorisations
\begin{equation*}
  \begin{tikzpicture}[diagram]
    \matrix[objects] {
      |(e)| (G_{x,y})_1 \&   |(v)|    (G_{x,y})_0        \& |(xy)| \{ (x,y) \} \\
      |(x2)| X_2 \& |(x1)| X_1 \& |(x0x0)| X_0 \times X_0 \\
    };
    \path[maps,->] 
      (x1) edge node[below] {$(s_0,t_0)$} (x0x0)
      (x2) edge[eabove] node[above] {$s_1$} (x1)
      (x2) edge[ebelow] node[below] {$t_1$} (x1)
      (xy) edge node[right] {$(x,y)$} (x0x0)
      (v) edge (xy)
      (v) edge (x1)
      (e) edge[eabove] node[above] {$s_1$} (v)
      (e) edge[ebelow] node[below] {$t_1$} (v)
      (e) edge (x2)
    ;
  \end{tikzpicture}
\end{equation*}
of the right hand square through the left hand square and hence graphs
$(G_{x,y})_1 \toto (G_{x,y})_0$ for all pairs $x,y \in X_0$. The set $X_0$
together with the graphs $G_{x,y}$ obviously form a $2$-graph $G_X$ and it is
easy to check that the construction $G \mapsto G_X$ is functorial and inverse
to the functor $X \mapsto X_G$ constructed above.  We have thus proven the
following well-known proposition: 
\begin{proposition}\label{prop:two-graphs-are-2-globular-sets}
  The category\/ $\dGhat_2$ of\/ $2$-globular sets and the category\/
  $\graph_2$ of\/ $2$-graphs are equivalent.
\end{proposition}

Our interest in $2$-globular sets stems from the fact that the category of
small $2$-categories is monadic over $\dGhat_2$, see e.\,g.\ the end of
section~2 in \cite{street;limits-indexed-by-category-valued-2-functors} or
\cite[Appendix~F]{leinster;higher-operads-higher-categories}.  Let us describe
the monadic adjunction $F_2 \colon \graph_2 \toot \cat_2 \cocolon U_2$ in terms
of the category of $2$-graphs and the monadic adjunction $F_1 \colon \graph
\toot \cat \cocolon U_1$ between graphs and small categories.  The forgetful
functor $U_2 \colon \cat_2 \to \graph_2$ maps a $2$-category $\cA$ to the
$2$-graph that has the objects of $\cA$ as vertices and associated with any two
vertices $a,b \in \ob\cA$ the graph $U_1 \cA(a,b)$. The left
adjoint $F_2$ takes a $2$-graph $G$ to the $2$-category $F_2G$ with objects
$V(G)$ and 
\begin{equation*}
  F_2G(x,y) = \coprod_{x = x_0,\dots,x_n = y} F_1G_{x_0,x_1} \times \dots \times F_1G_{x_{n-1},x_n},
\end{equation*}
where the coproduct ranges over all positive integers $n$ and all sequences
$(x_0,\dots,x_n)$ of elements of $V(G)$ with $x_0 = x$ and $x_n = y$.

Using the fact that $U_1$ preserves coproducts and arbitrary limits, we obtain
the following description of the monad $T_2 = U_2F_2$ on $\graph_2$:
\begin{proposition}\label{prop:2-cat-monad}
  The monad $T_2$ for\/ $2$-categories takes a $2$-graph\/ $G$ to the\/ $2$-graph\/ $T_2G$ on the same objects but with 
  \begin{equation*}
    T_2G_{x,y} = \coprod_{x = x_0,\dots,x_n = y} T_1G_{x_0,x_1} \times \dots \times T_1G_{x_{n-1},x_n},
  \end{equation*}
  where $T_1$ denotes the monad on $\graph$ for categories. One has to be
  careful to add another summand $\done $ in the description of\/ $T_2$ in case
  that $x = y$ to cater for the identities of $F_2G$.
\end{proposition}

In fact, Leinster shows in
\cite[Appendix~F]{leinster;higher-operads-higher-categories} that $T_2$ (and
even $T_n$ for any $n \in \dN$) is a cartesian monad.

\subsection*{\texorpdfstring{$2$}{2}-computads}
\label{sub:2_computads}
The category of $2$-globular sets and the monad $T_2$ suffer from one serious
drawback: The shapes allowed to generate a $2$-category are rather restrictive.
Each $2$-cell in a $2$-globular set has only one source and one target
$1$-cell. However, it occurs frequently that one has more complicated diagrams
that should generate a free $2$-category. This defect can be remedied by
passing from $2$-globular sets to Street's category of computads at the expense
of more involved combinatorics occuring in the base category. 

The category of computads is a category of presheaves, see
\cite{carboni-johnstone;connected-limits-familial-representability-and-artin-glueing},
and one can either guess or compute the index category $\cC$ for $2$-computads
from the proof of Carboni and Johnstone. The index category is given
as a collage and hence explicitly computable. We circumvent this
technicality and define the category of $2$-computads as the category of
presheaves on the category $\cC$ that can be described as follows: The objects
of $\cC$ are $0$, $1$ and $\gamma_{n,m}$ for all pairs $n,m \in \dN$.  The
morphisms are generated by
\begin{enumerate}
  \item $\cC(0,1) = \{\sigma,\tau\}$,
  \item $\cC(1,\gamma_{n,m}) = \{ \sigma_1,\dots,\sigma_n \} \cup \{ \tau_1,\dots,\tau_m\}$
\end{enumerate}
subject to the relations
\begin{enumerate}
  \item $\sigma_1 \sigma = \tau_1 \sigma$ and $\sigma_n \tau = \tau_m \tau$,
  \item $\sigma_{i+1} \sigma = \sigma_{i} \tau$ for all $0 \leq i < n$,
  \item $\tau_{i+1} \sigma = \tau_i \tau$ for all $ 0\leq i < m$.
\end{enumerate}

A $2$-computad thus consists of a graph $G$ together with a family $B_{n,m}$ of
globs with $\dom B_{n,m}$ a path of length $n$ and $\cod B_{n,m}$ a path of
length $m$.  It follows from this description that the terminal computad
$\done$ has one vertex $\pt$ with one loop $t\colon \pt \to \pt$ and one glob
$B_{n,m}$ of each size. We can moreover consider any globular graph as a
$2$-computad, the set $B_{n,m}$ of globs given by the interior faces with appropriate
domains and codomains.

It was proven in \cite{street;limits-indexed-by-category-valued-2-functors}
that $2$-categories are monadic over $\cChat$ and it follows from the explicit
description of the monad $T_2$ on $\cChat$ for $2$-categories given there and
our description of the terminal computad $\done$ that the cells in $T_2 \done$
certainly include our globular graphs. One might thus hope for a
relation between globular graphs and the \emph{canonical arities} of the
monad $T_2$ in the sense of
\cite{weber;generic-morphisms-parametric-representations-and-weakly-cartesian-monads,weber;familial-2-functors-and-parametric-right-adjoints,berger-mellies-weber;monads-with-arities-and-their-associated-theories}.

\subsection*{Nerves of globular graphs revisited}
\label{sub:nerves_of_globular_graphs_revisited}

So far, we did not speak of morphisms between globular graphs. It turns out
that there are several choices and all of them have their merits. However, if
one wants to give a presentation of the nerve of globular graphs in terms of a
cosimplicial object, $2$-functors between the free $2$-categories on them are
the correct choice.

By virtue of the preceding paragraph, we consider globular graphs as certain
types of $2$-computads.  Let $B_n$ denote the globular graph from
\autoref{ex:nerve-bn}. The free $2$-category $F_2B_n$ on $B_n$ has $2$ objects
$s$, $t$ and $n+1$ $1$-cells $e_0,\dots,e_n \colon s\to t$ and a $2$-cell
$\lambda_{i,j} \colon e_i \to e_j$ whenever $i \leq j$.\footnote{Observe that
one obtains essentially the nerve of the globular graph $B_n$ upon applying the
nerve functor $\nerve\colon \cat \to \sset$ to each $\hom$-set in this free
$2$-category. This is, of course, no coincidence but true more generally.} It is immediate from
this description that $B_\ph$ is a cosimplicial object in the category whose
objects are globular graphs and whose morphisms are $2$-functors between the
free categories.

Consider a $2$-functor $f\colon F_2B_n \to F_2G$ into the free $2$-category on
some globular graph $G$ that preserves source and target. It maps each edge
$e_i$ of $B_n$ to a morphism in $F_2G$, i.\,e.\ to an $st$-path in $G$.
Moreover, each $2$-cell $\lambda_{i,j}$ is mapped to a $2$-cell in $F_2G$,
i.\,e.\ a pasting of globs of $G$. Taking the boundary of this pasting, we
obtain a glob $\gamma$ in $G$ witnessing $f(e_i) \leq f(e_j)$.
Summing up, we thus have the following proposition:
\begin{proposition}\label{prop:nerve-as-hom}
  The nerve\/ $\nerve(G)$ of a globular graph $G$ is isomorphic to the
  simplicial set $[F_2B_\ph,F_2G]_{s,t}$ of source and target preserving
  $2$-functors from the cosimplicial object $F_2B_\ph$ into $F_2G$. 
\end{proposition}
Note that \autoref{prop:nerve-as-hom} is equivalent to the assertion that
$\nerve(G)$ is naturally isomorphic to the nerve of the category $F_2G(s,t)$.

\section{Pasting diagrams}
\label{sub:pasting_diagrams}

\begin{definition}\label{def:pasting-diagram}
  A \emph{pasting diagram} $(G,\sS)$ consists of a globular graph\/ $G$ and a
  set\/ $\sS$ of globular subgraphs of\/ $G$ that contains all paths and is
  closed under taking subgraphs.
\end{definition}
\begin{remark}\label{rem:closed-under-taking-subgraphs}
  The condition that $\sS$ be closed under taking subgraphs in
  \autoref{def:pasting-diagram} is not essential to the notion of pasting
  diagram but merely convenient for some considerations below, see e.\,g.\ the
  formulation of \autoref{def:restricted-pd}.  Moreover, if $\sS$ is an
  arbitrary set of globular subgraphs of $G$, then there is a unique smallest
  set $\langle \sS \rangle \supseteq \sS$ of globular subgraphs of $G$ that is
  closed under taking subgraphs and contains all paths in $G$.  In this
  situation, we call $(G,\langle \sS \rangle)$ the pasting diagram
  \emph{generated} by $(G,\sS)$.  In fact, we often give only $(G,\sS)$ and
  understand that we actually mean the pasting diagram $(G,\langle \sS
  \rangle)$ generated by it.
\end{remark}

\begin{definition}\label{def:inclusion}
  Let $\Sigma = (G,\sS)$ and\/ $\Pi = (G,\sT)$ be two pasting diagrams on the
  same underlying graph\/ $G$. We then say that\/ $\Sigma \to \Pi$ is an
  \emph{inclusion of pasting diagrams} if\/ $\sS \subseteq \sT$.
\end{definition}
\begin{remark}
  We warn the reader that any inclusion $\Sigma \to \Pi$ in the sense of
  \autoref{def:inclusion} is necessarily the identity on the underlying graphs.
\end{remark}

\begin{example}
  For any globular graph $G$ there exist minimal and maximal pasting diagrams
  $\Sigma_{\min} = (G,\sS_{\min})$ and $\Pi_{\max} = (G,\sS_{\max})$ with
  underlying graph $G$. The collections $\sS_{\min}$ and $\sS_{\max}$ of
  globular subgraphs of $G$ are given by 
  \begin{equation*}
    \sS_{\min} = \bigl\{ A \subseteq G \mid A \text{ is a face of $G$ or a path in $G$} \bigr\}
  \end{equation*}
  and
  \begin{equation*}
    \sS_{\max} = \bigl\{ A \subseteq G \mid A \text{ an arbitrary globular subgraph of $G$} \bigr\},
  \end{equation*}
  respectively. The maximal pasting diagram on $G$ is generated by $\{G\}$, of
  course. Moreover, we have an obvious inclusion $\Sigma_{\min} \to \Pi_{\max}$. 
\end{example}

\begin{definition}\label{def:closed-under}
  Let\/ $\Sigma = (G,\sS)$ be a pasting diagram.
  \begin{enumerate}[label=(\alph*)]
    \item We call\/ $\Sigma$ \emph{generated by wide subgraphs} if there is a set\/ $\sR$ of wide subgraphs of\/ $G$ such that\/ $(G,\sS) = (G,\langle \sR
      \rangle)$. 
    \item 
      We call\/ $\Sigma$ \emph{complete} if\/ $\sS$ is closed under taking joins, that is,
      if\/ $H_1,H_2 \in \sS$ implies $H_1 \join H_2 \in \sS$ whenever the join $H_1 \join H_2$ is defined.
    \item We call\/ $\Sigma$ \emph{closed under taking
      subdivisions} if for any subdivision $K$ of some $H \in \sS$ we have $K
      \in \sS$. 
  \end{enumerate}
\end{definition}

\begin{remark}\label{def:minimal-complete}
  Given a pasting diagram $\Sigma = (G,\sS)$ on some globular graph $G$, there
  exists a minimal complete pasting diagram $\Sigma^c = (G,\sS^c)$ on $G$ such
  that $\sS \subseteq \sS^c$. 
\end{remark}

\begin{example}
  Both the minimal pasting diagram $\Sigma_{\min}$ and the maximal pasting
  diagram $\Pi_{\max}$ on some globular graph $G$ are closed under taking
  subdivisions. Moreover, $\Pi_{\max}$ is always complete and generated by wide
  subgraphs. However, $\Sigma_{\min}$ need neither be complete nor generated by
  wide subgraphs as can be seen for e.\,g.\ $G = B_1 \join B_1$, where
  $\Sigma_{\min}^c = \Pi_{\max}$.
\end{example}

\begin{lemma}\label{rem:generated-by-wide-subgraphs}
  Any complete pasting diagram is generated by wide subgraphs.
\end{lemma}
\begin{proof}
  Let $\Sigma = (G,\sS)$ be a complete pasting diagram and consider any $A \in
  \sS$.  There are paths $p$ and $q$ from $s$ to $s(A)$ and from $t(A)$ to $t$,
  respectively. We then have $B(A) = p \join A \join q \in \sS$ as $\sS$ is
  closed under taking joins and the set $\{ B(A) \mid A \in \sS \}$ is a
  generating set of wide subgraphs.
\end{proof}

\begin{definition}\label{def:restricted-pd}
  Let\/ $\Sigma = (G,\sS)$ be a pasting diagram and let\/ $H \subseteq G$ be a globular
  subgraph. We call the set 
  \begin{equation*}
    \sS_H = \{ A \in \sS \mid A \subseteq H \}.
  \end{equation*}
  the \emph{restriction} of\/ $\sS$ to\/ $H$ and refer to the pasting diagram\/
  $\Sigma_H = (H,\sS_H)$ as the \emph{restriction} of\/ $\Sigma$ to\/ $H$. In
  the case that $H = G_{x,y}$ for two vertices $x,y \in G$, we also write
  $\Sigma_{x,y}$ for the restriction of\/ $\Sigma$ to $H$.
\end{definition}

\begin{remark}
  Observe that if $\Sigma$ is complete, closed under taking subdivisions, or
  generated by wide subgraphs, then so are all of its restrictions. This
  implies in particular that $(\Sigma_H)^c = (\Sigma^c)_H$.
\end{remark}

\begin{definition}\label{def:join-pd}
  Let\/ $(G_1,\sS)$ and\/ $(G_2,\sT)$ be pasting diagrams generated by collections\/
  $\sS_0 \subseteq \sS$ and\/ $\sT_0 \subseteq \sT$ of wide subgraphs.  Their\/
  \emph{join}\/ $(G_1,\sS) \join (G_2,\sT)$ is the pasting diagram with
  underlying graph\/ $G_1 \join G_2$ that is generated by the set
  \begin{equation*}
    \sS_0 \join \sT_0 = \{ A \join B \mid A \in \sS_0 \text{ and } B \in \sS_1 \}.
  \end{equation*}
\end{definition}
\begin{remark}\label{rem:join-pd-well-defined}
  We have to justify that the join of two pasting diagrams as in
  \autoref{def:join-pd} is well-defined, i.\,e.\ independent of the choice of
  the generating sets $\sS_0$ and $\sT_0$. To this end, consider another pair
  $\sS_0' $ and $\sT_0' $ of collections of wide subgraphs that generate $\sS$
  and $\sT$, respectively. For each $A \in \sS_0$ and $B \in \sT_0$ we then
  find $A' \in \sS_0'$ and $B' \in \sT_0'$ such that $A \subseteq A'$ and $B
  \subseteq B'$. But this implies $A \join B \subseteq A' \join B'$. As the
  argument is symmetric in $\sS_0, \sT_0$ and $\sS_0', \sT_0'$, we may conclude
  that the pasting diagrams generated by $\sS_0 \join \sT_0$ and $\sS_0' \join
  \sT_0'$, respectively, coincide.
\end{remark}

\begin{remark}\label{rem:complete-and-join-decomposition}
  Suppose that $\Sigma$ is some pasting diagram on $G_1 \join G_2$ and let
  $\Sigma_i$ be the restriction of $\Sigma$ to $G_i$. We warn the reader that
  the obvious inclusion $\Sigma \subseteq \Sigma_1 \join \Sigma_2$ is generally
  strict as demonstrated by the minimal pasting diagram $\Sigma_{\min}$ on $B_1
  \join B_1$.  In the case that $\Sigma$ is complete, however, we always have an
  equality $\Sigma = \Sigma_1 \join \Sigma_2$. 
\end{remark}

\section{Nerves of Pasting Diagrams}
\label{sub:nerves_of_pasting_diagrams}

In this section we define the nerve $\nerve(\Sigma)$ of a pasting diagram
$\Sigma = (G,\sS)$ as a certain simplicial subset of $\nerve(G)$ and record
some elementary properties. We also observe that the pictorial calculus for
$\nerve(G)$ from \autoref{prop:nmarked-is-simplex} and
\ref{prop:marked-simplicial} carries over to nerves of complete pasting
diagrams. The larger part of this section, however, consists of technical
lemmata concerning the interplay of the nerve with operations such as
restriction of pasting diagrams or the join of two pasting diagrams that we
introduced above.  These basic results will be used throughout the remaining
chapters.

\subsection*{Definition and Elementary Properties}
We now give the promised definition of the nerve of a pasting diagram as a
subset of the nerve of the underlying globular graph.
\begin{definition}\label{def:nerve-of-pds}
  Let $\Sigma = (G,\sS)$ be a pasting diagram. Its nerve\/ $\nerve(\Sigma)$ is
  the simplicial subset of\/ $\nerve(G)$ consisting of those simplices $\sigma
  = (p_0 \leq \dots \leq p_n)$ with the property that there exists some $A \in
  \sS$ such that all the relations $p_{i-1} \leq p_{i}$ for $1 \leq i \leq n$
  have witnesses $\gamma_i \subseteq A$.
\end{definition}

\begin{example}\label{ex:spine-as-nerve}
  Let us consider a pasting diagram on the graph $B_n$ from
  \autoref{ex:nerve-bn}.  We then have $\nerve(B_n,\sS) = \Delta^1 \join \dots
  \join \Delta^1$ for $\sS$ the set of all interior faces of $B_n$.
\end{example}

\begin{lemma}\label{lem:nerves-of-subdivision-closed}
  Let $(G,\sS) \to (G,\sT)$ be an inclusion and suppose that $\sS$ is closed
  under taking subdivisions. If an inner face of a simplex of\/
  $\nerve(G,\sT)$ is contained in $\nerve(G,\sS)$, then so is the simplex
  itself. 
\end{lemma}
\begin{proof}
  Consider a simplex $\sigma \subseteq \nerve(G,\sT)$ with $d_i\sigma$ in
  $\nerve(G,\sS)$. Write 
  \begin{equation*}
    \sigma = (p_0 \leq \dots \leq p_n)
  \end{equation*}
  and consider the face $d_i \sigma$ along with a witness $\gamma \subseteq A$
  for $p_{i-1} \leq p_{i+1}$ and some $A \in \sS$. Cutting $\gamma$ along $p_i$
  provides us with witnesses $\gamma_i$, $\gamma_{i+1}$ for $p_{i-1} \leq p_i \leq
  p_{i+1}$.  Possibly passing to a subdivision of $A$ we may assume
  $\gamma_i,\gamma_{i+1} \subseteq A$ and thus find $\sigma \subseteq
  \nerve(G,\sS)$. 
\end{proof}

In the case that $\Sigma$ is complete, we have the following lemma that allows
us to describe $\nerve(\Sigma)$ in terms of admissible $n$-marked subgraphs of
$G$. 
\begin{lemma}\label{lem:pictorial-desc-nerve-pds}
  Consider a complete pasting diagram $\Sigma = (G,\sS)$. An $n$-simplex
  $\sigma$ with associated $n$-marked subgraph $(P_\sigma,\lambda_\sigma)$ of
  $\nerve(G)$ is contained in $\nerve(\Sigma)$ if and only if $P_\sigma \in
  \sS$. 
\end{lemma}
\begin{proof}
  Write $\sigma = (p_0 \leq \dots \leq p_n)$.  Recall that $P_\sigma =
  \bigcup_i p_i$.  It is immediate that $\sigma \in \nerve(\Sigma)$ if
  $P_\sigma \in \sS$.  Now suppose that $\sigma \in \nerve(\Sigma)$.  We then
  find some $A \in \sS$ that contains all the witnesses $\gamma_i$ for $p_{i-1}
  \leq p_i$. We also find possibly trivial paths $q_1$ from $s(G)$ to $s(A)$
  and $q_2$ from $t(A)$ to $t(G)$ such that all the paths $p_i$ can be written
  as $p_i = q_1 \cdot r_i \cdot q_2$ for some path $r_i$ from $s(A)$ to $t(A)$
  in $A$. It now follows that $R = \bigcup_i r_i \in \sS$ as $\sS$ is closed
  under taking subgraphs and this in turn implies $P = q_1 \join  R \join q_2
  \in \sS$ as $\sS$ is closed under taking joins.
\end{proof}
\begin{corollary}\label{cor:pictorial-pd-nerve}
  Let $\Sigma = (G,\sS)$ be a complete pasting diagram.  The $n$-simplices of\/
  $\nerve(\Sigma)$ are in bijection with admissible $n$-marked subgraphs
  $(P,\lambda)$ of\/ $G$ such that\/ $P \in \sS$. Moreover, the action of
  simplicial operators on these $n$-marked subgraphs may be computed as
  described in \autoref{prop:marked-simplicial}.
\end{corollary}

\subsection*{Nerves and the join operation}
In this paragraph we compute the nerve $\nerve(\Sigma_1 \join \Sigma_2)$ of the
join of two pasting diagrams that are generated by wide subgraphs. More
precisely, we show that $\nerve$ is a strong monoidal functor from the monoidal
category of pasting diagrams that are generated by wide subgraphs with the join
operation to the cartesian closed category of simplicial sets. This observation
will be used in \autoref{sec:scat} to construct composition laws for the
simplicial category that we associate with a pasting diagram.
\needspace{4\baselineskip}
\begin{proposition}\label{lem:nerve-of-join}
  Let\/ $\Sigma$ and\/ $\Pi$ be pasting diagrams that are generated by wide
  subgraphs. There then exists an isomorphism
  \begin{equation*}
    \phi_{\Sigma,\Pi} \colon \nerve(\Sigma) \times \nerve(\Pi) \to
    \nerve(\Sigma \join \Pi)
  \end{equation*}
  and these isomorphisms $\phi_{\Sigma,\Pi}$ equip the nerve with the structure
  of a strong monoidal functor from the monoidal category of pasting diagrams
  that are generated by wide subgraphs with the join operation to the 
  cartesian closed category of simplicial sets
\end{proposition}
\begin{proof}
  Given $n$-simplices $\sigma^k= (p_0^k \leq \dots \leq
  p_n^k)$ of $\nerve(\Sigma_k)$ we have an
  \hbox{$n$-simplex} $\sigma^1 \cdot \sigma^2 = (p_0^1 \cdot p_0^2 \leq \dots
  \leq p_n^1 \cdot p_n^2)$ in $\nerve( \Sigma_1 \join \Sigma_2)$. One
  easily verifies that the assignment $(\sigma^1,\sigma^2) \mapsto
  \sigma^1\cdot\sigma^2$ defines a map
  \begin{equation*}
    \nerve(\Sigma_1 ) \times
    \nerve(\Sigma_2 ) \to \nerve(\Sigma_1 \join \Sigma_2).
  \end{equation*}
  Let us show that this map is an isomorphism. To this end write $\Sigma_k =
  (G_k,\sS_k)$ and denote by $s_k$ and $t_k$ the source and target of $G_k$.
  Consider an $n$-simplex $\sigma = (p_0 \leq \dots \leq p_n)$ of~$\nerve(
  \Sigma_1 \join \Sigma_2)$.  Each path $p_i$ necessarily contains the cut
  vertex~$t_1 = s_2$ of $G_1 \join G_2$ and thus decomposes uniquely as $p_i =
  q_i^1 \cdot q_i^2$ for $s_k t_k$-paths $q_i^k$ in $G_k$.  We have thus found
  a unique pair $(\sigma^1,\sigma^2)$ of $n$-simplices $\sigma^k = (q_0^k \leq
  \dots \leq q_n^k)$ of $\nerve(\Sigma_k)$ such that $\sigma = \sigma^1 \cdot
  \sigma^2$. 

  Note that the pasting diagram $\Sigma_0$ on the trivial globular graph serves
  as unit in our monoidal category of pasting diagrams and we obviously have
  $\nerve(\Sigma_0) = \Delta^0$. Moreover, one of the axioms of a strong monoidal
  functor requires the diagram
 \begin{equation*}
    \begin{tikzpicture}[diagram]
      \matrix[objects] {
        |(prod)| \nerve(\Sigma_1) \times \nerve(\Sigma_2) \times \nerve(\Sigma_3) \&[+3em] |(prodjoin)| \nerve(\Sigma_1) \times \nerve(\Sigma_2 \join \Sigma_3)\\[-2ex]
        |(joinprod)| \nerve(\Sigma_1 \join \Sigma_2) \times \nerve(\Sigma_3) \& |(join)| \nerve(\Sigma_1 \join \Sigma_2 \join \Sigma_3) \\
      };
      \path[maps,->] 
        (prod) edge node[above] {$\nerve(\Sigma_{1}) \times \phi$} (prodjoin)
        (prod) edge node[left] {$\phi \times \nerve(\Sigma_2)$} (joinprod)
        (prodjoin) edge node[right] {$\phi$} (join)
        (joinprod) edge node[below] {$\phi$} (join)
      ;
    \end{tikzpicture}
  \end{equation*}
  to be commutative. Given our definition of $\phi$ this is immediate, though.
  The remaining axioms of a strong monoidal functor are left to the reader. 
\end{proof}

\begin{corollary}\label{cor:nerve-of-join}
  If\/ $\Sigma$ is a complete pasting diagram on\/ $G = G_1 \join G_2$, then\/
  $\Sigma = \Sigma_1 \join \Sigma_2$ and hence
  \begin{equation*}
    \nerve(\Sigma) \iso \nerve(\Sigma_1) \times \nerve(\Sigma_2),
  \end{equation*}
  where $\Sigma_i$ denotes the restriction of\/ $\Sigma$ to $G_i$.
\end{corollary}

\begin{remark}
  We did not specify the maps in the category that is to be the domain of the
  strong monoidal functor $\nerve$. One could choose e.\,g.\ maps $(H,\sS) \to
  (G,\sT)$ with $H \subseteq G$ a globular subgraph and $\sS \subseteq \sT$.
  The nerve is certainly functorial with respect to these maps and one easily
  verifies that the isomorphisms $\phi$ in \autoref{lem:nerve-of-join} are
  natural.
\end{remark}

\begin{remark}
  The hypothesis that the pasting diagrams in \autoref{lem:nerve-of-join} are
  generated by wide subgraphs is essential to the construction of the
  isomorphism $\phi$, for otherwise $\phi(\sigma^1,\sigma^2) = \sigma^1 \cdot
  \sigma^2$ need not even be a simplex of $\nerve(\Sigma_1 \join \Sigma_2)$.
\end{remark}

\begin{remark}\label{rem:pictorial-join-isos}
  The isomorphisms $\phi$ constructed in \autoref{lem:nerve-of-join} can also
  be described in terms of $n$-marked subgraphs. The image $\phi(\sigma,\tau)$
  of an $n$-simplex $(\sigma,\tau) \in \nerve(\Sigma_1) \times
  \nerve(\Sigma_2)$ with admissible $n$-marked subgraphs
  $(P_\sigma,\lambda_\sigma)$ and $(P_\tau,\lambda_\tau)$ has $(P_\sigma \join
  P_\tau, \lambda_\sigma \join \lambda_\tau)$ as its associated $n$-marked
  subgraph. Put differently, the isomorphisms $\phi$ are given by the join of
  the respective $n$-marked subgraphs.
\end{remark}

\subsection*{Unions and intersections}
In this paragraph we collect some technical yet easy lemmata concerning
intersections and unions of nerves of pasting diagrams. 
\begin{lemma}\label{lem:nerve-of-union}
  Consider a family\/ $(G,\sS_j)$,\/ $j \in J$, of pasting diagrams on the same
  underlying graph\/ $G$. We then have an equality
  \begin{equation*}
    \bigcup_{j \in J} \nerve(G,\sS_j) = \nerve\Bigl( G, \bigcup_{j \in J} \sS_j \Bigr)
  \end{equation*}
  of simplicial subsets of\/ $\nerve(G)$.
\end{lemma}
\begin{proof}
  Consider some simplex $\sigma = (p_0 \leq \dots \leq p_n)$ in $\nerve(G)$.
  We have $\sigma \in \bigcup_{j } \nerve(G,\sS_j)$ if and only if there exists
  some $j \in J$ and some $A \in \sS_j$ such that all the witnesses of $p_{i-1}
  \leq p_i$ are contained in $A$. This is obviously equivalent to the condition
  that there exists some $A\in \bigcup_j \sS_j$ such that all the witnesses of
  $p_{i-1} \leq p_i$ are contained in $A$, i.\,e.\ to the statement $\sigma \in
  \nerve(G,\bigcup_j \sS_j)$.
\end{proof}
\begin{lemma}\label{lem:nerve-of-restriction}
  Let\/ $(G,\sS)$ be a pasting diagram, $H \subseteq G$ a wide globular subgraph
  and\/ $\sS_H$ the restriction of\/ $\sS$ to\/ $H$. We then have a cartesian square 
  \begin{equation*}
    \begin{tikzpicture}[diagram]
      \matrix[objects,row sep={4.25em,between origins}] {
        |(hs)| \nerve(H,\sS_H) \& |(gs)| \nerve(G,\sS_G) \\
        |(h)| \nerve(H) \& |(g)| \nerve(G) \mathrlap{.} \\
      };
      \path[maps,->] 
        (hs) edge (gs)
        (hs) edge (h)
        (gs) edge (g)
        (h) edge (g)
      ;
    \end{tikzpicture}
  \end{equation*}
\end{lemma}
\begin{proof}
  Immediate from the definitions.
\end{proof}
\begin{corollary}\label{cor:nerve-pullback}
  Let\/ $(G,\sS) \subseteq (G,\sT)$ be an inclusion of pasting diagrams, $H \subseteq G$ a wide
  globular subgraph and\/ $\sS_H$ and $\sT_H$ the restrictions of\/ $\sS$ and
  $\sT$ to\/ $H$. We then have a cartesian square
  \begin{equation*}
    \begin{tikzpicture}[diagram]
      \matrix[objects] {
        |(hs)| \nerve(H,\sS_H) \& |(gs)| \nerve(G,\sS) \\[-0.75em]
        |(h)| \nerve(H,\sT_H) \& |(g)| \nerve(G,\sT) \mathrlap{.} \\
      };
      \path[maps,->] 
        (hs) edge (gs)
        (hs) edge (h)
        (gs) edge (g)
        (h) edge (g)
      ;
    \end{tikzpicture}
  \end{equation*}
\end{corollary}
\begin{lemma}\label{lemma:nerve-pushout-as-nerve}
  Let\/ $(G,\sS)$ and\/ $(H,\sT)$ be pasting diagrams such that\/ $H$ is
  a wide subgraph of\/ $G$. Further let $\sS_H$ be the restriction of\/
  $\sS$ to\/ $H$ and suppose $\sS_H \subseteq \sT$, i.\,e.\ $\sS_H = \sS
  \cap \sT$. Then
  \begin{equation*}
    \begin{tikzpicture}[diagram]
      \matrix[objects] {
        |(hs)| \nerve(H,\sS_H) \& |(gs)| \nerve(G,\sS) \\
        |(ht)| \nerve(H,\sT) \& |(gst)| \nerve(G,\sS \cup \sT) \\
      };
      \path[maps,->] 
        (hs) edge (gs)
        (hs) edge (ht)
        (gs) edge (gst)
        (ht) edge (gst)
      ;
    \end{tikzpicture}
  \end{equation*}
  is bicartesian.
\end{lemma}
\begin{proof}
  Let us first show that the square is cocartesian. We have jointly
  surjective inclusions $\nerve(G,\sS) \to \nerve(G,\sS \cup \sT)$ and
  $\nerve(H,\sT) \to \nerve(G,\sS \cup \sT)$ since $H \subseteq G$ is
  wide by assumption. It is furthermore easy to see that the intersection
  of $\nerve(H,\sT)$ and $\nerve(G,\sS)$ in $\nerve(G,\sS\cup\sT)$ is
  precisely $\nerve(H,\sS_H )$. In fact, if $\sigma = (p_0 \leq \dots
  \leq p_n)$ is a simplex in this intersection, then $p_i \subseteq H$
  for all $i$ and we consequently find some $A \in \sS$, $A \subseteq H$,
  that contains witnesses for all the relations $p_{i-1} \leq p_i$. This
  implies $\sigma \in \nerve(H,\sS_H)$.

  It now either follows from general properties of coherent
  categories or from \autoref{cor:nerve-pullback} that the square is cartesian,
  too.
\end{proof}

\begin{lemma}\label{lem:intersection-1}
  Consider a complete pasting diagram\/ $(G,\sS)$ on some globular graph~$G$
  with source\/ $s$ and target\/ $t$.  Let\/ $H_1$ and\/ $H_2$ be wide
  subgraphs of\/ $G$ and let\/ $H_0 = H_1 \cap H_2$. Denote by\/ $\sS_i$ the
  restriction of\/ $\sS$ to\/ $H_i$.
  \begin{enumerate}[label=(\alph*)]
    \item The intersection of\/ $\nerve(H_1,\sS_1)$ and\/ $\nerve(H_2,\sS_2)$
      in\/ $\nerve(G,\sS)$ is nonempty if and only if\/ $H_0$ contains a path
      from\/ $s$ to\/ $t$.
    \item If\/ $H_0$ is a wide globular subgraph, then 
      \begin{equation*}
        \begin{tikzpicture}[diagram]
          \matrix[objects] {
            |(hk)| \nerve(H_0,\sS_0) \& |(h)| \nerve(H_1,\sS_1) \\
            |(k)| \nerve(H_2,\sS_2) \& |(g)| \nerve(G,\sS) \\
          };
          \path[maps,->] 
            (hk) edge (h)
            (hk) edge (k)
            (k) edge (g)
            (h) edge (g)
          ;
        \end{tikzpicture}
      \end{equation*}
      is cartesian. 
  \end{enumerate}
\end{lemma}
\begin{proof}
  Part~(a) is obvious, for $st$-paths $p$ from $s$ to $t$ in $H_0$ are in
  bijective correspondence with $0$-simplices $(p)$ in $\nerve(\Sigma_1) \cap
  \nerve(\Sigma_2)$. 

  Let us now suppose that $H_0$ is a wide globular subgraph of $G$.  It is
  clear that $\nerve(G_0,\sS_0)$ is a simplicial subset of the intersection of
  $\nerve(H_1,\sS_1)$ and $\nerve(H_2,\sS_2)$. It thus suffices to show $
  \nerve(H_1,\sS_1)  \cap \nerve(H_2,\sS_2) \subseteq \nerve(H_0,\sS_0)$. To
  this end, consider any simplex $\sigma = (p_0 \leq \dots \leq p_n)$ in the
  intersection of $\nerve(H_1,\sS_1)$ and $\nerve(H_2,\sS_2)$. Observe that all
  the paths $p_i$ are paths in $H_0 = H_1 \cap H_2$ and that $\bigcup_i p_i \in
  \sS_k$ for $k \in \{1,2\}$ by virtue of \autoref{lem:pictorial-desc-nerve-pds}.  We
  thus have $\bigcup_i p_i \in \sS_1 \cap \sS_2 = \sS_0$ and conclude $\sigma
  \in \nerve(H_0,\sS_0)$ by \autoref{lem:pictorial-desc-nerve-pds}.
\end{proof}
\needspace{4\baselineskip}
\begin{corollary}\label{cor:nerve-join-intersection}
  Let\/ $\Sigma$ be a complete pasting diagram on a globular graph\/ $G$ with source $s$ and targe\/t $t$. Further let $x,y \in G$. If there
  exists a directed path from $x$ to $y$ in $G$, then 
  \begin{equation*}
    \begin{tikzpicture}[diagram]
      \matrix[objects] {
        |(xy)| \nerve(\Sigma_{s,x} \join \Sigma_{x,y} \join \Sigma_{y,t}) \& |(x)| \nerve(\Sigma_{s,x} \join \Sigma_{x,t} ) \\
        |(y)| \nerve( \Sigma_{s,y} \join \Sigma_{y,t} ) \& |(n)| \nerve(\Sigma) \\
      };
      \path[maps,->] 
        (xy) edge (x)
        (xy) edge (y)
        (y) edge (n)
        (x) edge (n)
      ;
    \end{tikzpicture}
  \end{equation*}
  is a cartesian square of simplicial sets, i.\,e.\
  \begin{equation*}
    \nerve(\Sigma_{s,x} \join \Sigma_{x,t}) \cap 
    \nerve(\Sigma_{s,y} \join \Sigma_{y,t}) 
    = \nerve( \Sigma_{s,x} \join \Sigma_{x,y} \join \Sigma_{y,t}).
  \end{equation*}
  Moreover, the intersection of\/ $\nerve(\Sigma_{s,x} \join \Sigma_{x,t})$ and
  $\nerve( \Sigma_{s,y} \join \Sigma_{y,t} )$ in $\nerve(\Sigma)$ is empty whenever $G$
  contains neither a directed path from $x$ to $y$ nor a directed path from $y$
  to $x$.
\end{corollary}
\begin{proof}
  Immediate from \autoref{lem:intersection-1} and \autoref{lem:intersection-of-xy-globs}.
\end{proof}

\subsection*{Formal partial composites}

In this short final paragraph, we define an operation \enquote{$\hc$} that
associates with an inclusion $\Sigma \to \Pi$ of pasting digarams a new pasting
diagram $\Sigma \hc \Pi$. In case that $\Pi$ is complete, $\Sigma \hc \Pi$ has
the property that the inclusion $\Sigma \to \Pi$ factors as 
\begin{equation*}
  \begin{tikzpicture}[diagram]
    \matrix[objects] {
      |(s)| \Sigma \&[-2em] \&[-2em] |(p)| \Pi \\[-2ex]
      \& |(sp)| \Sigma \hc \Pi \smash{.}\\
    };
    \path[maps,->] 
      (s) edge (p)
      (s) edge (sp)
      (sp) edge (p)
    ;
  \end{tikzpicture}
\end{equation*}
The pasting diagram $\Sigma \hc \Pi$ should hence be thought of as a formal
partial composite of the inclusion $\Sigma \to \Pi$. In fact, $\Sigma \hc \Pi$
features in exactly this role in \autoref{chap:global} and we advise the reader
to skip this section on a first reading and only come back to it after having
had a first look at the material in \autoref{chap:global}.

\begin{definition}\label{def:hc}
  Let\/ $\Sigma \to \Pi$ be an inclusion of pasting diagrams\/ $\Sigma = (G,\sS)$
  and\/ $\Pi = (G,\sT)$. We define a pasting diagram
  $\Sigma \hc \Pi = (G,\sS \hc \sT)$ on $G$ by
  \begin{equation*}
    \sS \hc \sT = \sS \cup\ \bigcup_{\mathclap{x \in G \smallsetminus \{s,t\}}}\ \sT_{s,x} \join \sT_{x,t}.
  \end{equation*}
\end{definition}

\begin{remark}
  We warn the reader that generally $\sS \hc \sT \nsubseteq \sT$ as can already
  be seen for $\Sigma = \Pi$ some non-complete pasting diagram. However, $\sS
  \hc \sT \subseteq \sT$ whenever $\sT$ is complete, see
  \autoref{lem:hc-is-full}.
\end{remark}

\begin{example}\label{ex:hc-example}
  Let us compute $\Sigma \hc \Pi$ for $\Sigma = (G,\sS_{\min}^c)$ the minimal complete and
  $\Pi = (G,\sS_{\max})$ the maximal pasting diagram on the globular graph $G$
  shown in \autoref{subfig:hc-graph}.  The definition of $\sS_{\min}^c \hc
  \sS_{\max}$ gives us
  \begin{equation*}
    \sS_{\min}^c \hc \sS_{\max} = \sS_{\min}^c \cup \bigl( \sS_{{\max},0,1} \join
    \sS_{{\max},1,2} \bigr).
  \end{equation*}
  The set $\sS_{{\max},0,1} \join \sS_{{\max},1,2}$ is the set of all globular
  subgraphs of $G_{0,1} \join G_{1,2}$, i.\,e.\ of the graph shown in
  \autoref{subfig:hc-subgraph}. One now easily sees that the only subgraphs in
  $(\sS \hc \sT) \smallsetminus \sS$ are $G_{0,1} \join G_{1,2}$ itself and the
  two of its subgraphs shown in \autoref{subfig:hc-new} and
  \ref{subfig:hc-glob}. 
  \begin{figure}
  \centering
  \begin{subfigure}[b]{0.45\textwidth}
    \centering
    \begin{tikzpicture}
      \node[fill,circle,inner sep=2pt] (s) at (0,0) {};
      \node[fill,circle,inner sep=2pt] (m) at (2,0) {}; 
      \node[fill,circle,inner sep=2pt] (t) at (4,0) {}; 
      \node[below] at (s) {$0$};
      \node[below] at (m) {$1$};
      \node[below] at (t) {$2$};
      \draw (m) edge[->,out=45,in=135] (t);
      \draw (m) edge[->,out=0,in=180] (t);
      \draw (m) edge[->,out=-45,in=-135] (t);
      \draw (s) edge[->] (m);
      \draw (s) edge[->,out=60,in=120] (t);
  \end{tikzpicture}
    \caption{A globular graph $G$}\label{subfig:hc-graph}
  \end{subfigure}
\hfil
  \begin{subfigure}[b]{0.45\textwidth}
    \centering
    \begin{tikzpicture}
      \node[fill,circle,inner sep=2pt] (s) at (0,0) {};
      \node[fill,circle,inner sep=2pt] (m) at (2,0) {}; 
      \node[fill,circle,inner sep=2pt] (t) at (4,0) {}; 
      \node[below] at (s) {$0$};
      \node[below] at (m) {$1$};
      \node[below] at (t) {$2$};
      \draw (m) edge[->,out=45,in=135] (t);
      \draw (m) edge[->,out=0,in=180] (t);
      \draw (m) edge[->,out=-45,in=-135] (t);
      \draw (s) edge[->] (m);
      \draw[white] (s) edge[->,out=60,in=120] (t);
  \end{tikzpicture}
    \caption{The subgraph $G_{0,1} \join G_{1,2}$ of $G$.}\label{subfig:hc-subgraph}
  \end{subfigure}
    \newline
  \begin{subfigure}[b]{0.45\textwidth}
    \centering
    \begin{tikzpicture}
      \node[fill,circle,inner sep=2pt] (s) at (0,0) {};
      \node[fill,circle,inner sep=2pt] (m) at (2,0) {}; 
      \node[fill,circle,inner sep=2pt] (t) at (4,0) {}; 
      \node[below] at (s) {$0$};
      \node[below] at (m) {$1$};
      \node[below] at (t) {$2$};
      \draw (m) edge[->,out=45,in=135] (t);
      \draw (m) edge[->,out=-45,in=-135] (t);
      \draw (s) edge[->] (m);
      \draw[white] (s) edge[->,out=60,in=120] (t);
  \end{tikzpicture}
    \caption{A wide subgraph of $G_{0,1} \join G_{1,2}$.} \label{subfig:hc-new}
  \end{subfigure}
  \hfil
  \begin{subfigure}[b]{0.45\textwidth}
    \centering
    \begin{tikzpicture}
      \node[fill,circle,inner sep=2pt] (s) at (0,0) {};
      \node[fill,circle,inner sep=2pt] (m) at (2,0) {}; 
      \node[below] at (s) {$1$};
      \node[below] at (m) {$2$};
      \draw (s) edge[->,out=45,in=135] (m);
      \draw (s) edge[->,out=-45,in=-135] (m);
  \end{tikzpicture}
    \caption{A glob $\gamma$ of $G_{0,1} \join G_{1,2}$.}\label{subfig:hc-glob}
  \end{subfigure}
  \caption{The graphs occuring in \autoref{ex:hc-example}.}\label{fig:inclusion-counterexample}
  \end{figure}
\end{example}

\begin{remark}\label{rem:hc-and-subgraphs}
  Consider an inclusion $\Sigma \to \Pi$ of pasting diagrams on some graph $G$.
  For any globular subgraph $H \subseteq G$, there exists a canonical inclusion
  \begin{equation*}
   (\Sigma_H \hc \Pi_H) \to (\Sigma \hc \Pi)_H .
  \end{equation*}
  This inclusion is strict in general, even for subgraphs $H$ of the form $H =
  G_{x,y}$. This can already be seen in \autoref{ex:hc-example}, where the glob
  shown in \autoref{subfig:hc-glob} is an element of $(\sS_{\min} \hc
  \sS_{\max})_{1,2}$ but $\sS_{\min,1,2} \hc \sS_{\max,1,2}$ coincides with
  $\sS_{\min,1,2}$ simply because there are no vertices $x \in G_{1,2}
  \smallsetminus \{1,2\}$.
\end{remark}

We end this paragraph with some technical but easy observations on
\enquote{$\hc$} that will be of good use in \autoref{chap:global}.
\needspace{4\baselineskip}
\begin{lemma}\label{lem:hc-is-complete}
  Let\/ $\Sigma \to \Pi$ be an inclusion of pasting digarams. If\/ $\Pi$ is
  complete, then so is\/ $\Sigma \hc \Pi$.
\end{lemma}
\begin{proof}
  Let us write $\Pi = (G,\sT)$ and $\Sigma = (G,\sS)$ and consider $A,B \in \sS
  \hc \sT$ with $t(A) = s(B) = x$. Observe that $A \subseteq G_{s,x}$ and $B
  \subseteq G_{x,t}$, i.\,e.\ $A \in \sT_{s,y} \join \sT_{y,x}$ and $B \in
  \sT_{x,,z} \join \sT_{z,t}$. But this implies $A \join B \in \sT_{s,x} \join
  \sT_{x,t} \subseteq \sS \hc \sT$ by completeness of $\Pi$.
\end{proof}
\begin{lemma}\label{lem:hc-is-full}
  Consider an inclusion $\Sigma \to \Pi$ of pasting diagrams on some globular
  graph\/ $G$ with source $s$ and target $t$. If\/ $\Pi$ is complete and $x$ and
  $y$ are two vertices of\/ $G$ such that $(x,y) \neq (s,t)$, then $(\Sigma \hc
  \Pi)_{x,y} = \Pi_{x,y}$.
\end{lemma}
\begin{proof}
  Let us write $\Pi = (G,\sT)$ and $\Sigma = (G,\sS)$. Let us first show that
  $\sT_{x,y} \subseteq (\sS \hc \sT)_{x,y}$. To this end, consider any $A \in
  \sT$ with $A \subseteq G_{x,y}$. If $x \neq s$ then $A \subseteq G_{x,y}
  \subseteq G_{x,t}$ and hence $A \in \sT_{x,t} \subseteq \sT_{s,x} \join
  \sT_{x,t} \subseteq \sS \hc \sT$. The case $y \neq t$ is handled analogously.

  Conversely, if $A \in (\sS \hc \sT)_{x,y}$, then $A \in \sS_{x,y} \subseteq
  \sT_{x,y}$ or $A \subseteq G_{x,y}$ with $A \in \sT_{s,z} \join \sT_{z,t}$
  for some vertex $z \in G \smallsetminus \{x,y\}$.  In the latter case, we
  conclude $A \in \sT_{x,y}$ by completeness of $\Pi$. 
\end{proof}
\begin{lemma}\label{lem:nerve-of-hc}
  Let\/ $\Sigma \to \Pi$ be an inclusion of pasting diagrams on the underlying
  graph $G$ with source $s$ and target $t$. We then have an equality
  \begin{equation*}
    \nerve(\Sigma \hc \Pi) = \nerve(\Sigma) \cup\ \bigcup_{\mathclap{x \in G \smallsetminus\{s,t\}}}\ \nerve(\Pi_{s,x} \join \Pi_{x,t})
  \end{equation*}
  of simplicial subsets of\/ $\nerve(G)$.
\end{lemma}
\begin{proof}
  This follows from the definition of $\Sigma \hc \Pi$ and
  \autoref{lem:nerve-of-union}.
\end{proof}

\chapter{Pasting Diagrams and Simplicial Categories}\label{chap:labeling}
In this chapter we associate with any complete pasting diagram $\Sigma$ a
simplicial category $\dC[\Sigma]$. It turns out that $\dC[\Pi_{\max}]$, where
$\Pi_{\max}$ denotes the maximal pasting diagram on some globular graph $G$, is
nothing but the simplicial category obtained from the free $2$-category on $G$
by local application of the nerve functor, see \autoref{rem:dc-maximal-pd}. We
are therefore led to think of functors $u \colon \dC[\Sigma] \to \dA$ with a
pasting diagram $\Sigma_{\min} \subseteq \Sigma \subseteq \Pi_{\max}$ as
mediating between the mere specification of maps and cells in $\dA$ at the
level of $\Sigma_{\min}$ and a fully coherent composition $\dC[\Pi_{\max}] \to
\dA$. There is one problem in this picture, though. In general, $\Sigma_{\min}$
is not complete and our definition of $\dC[\Sigma]$ does not work. In order to
overcome this difficulty, we introduce labelings of globular graphs in a
simplicial category $\dA$ so as to capture the idea of a compatible
specification of maps and cells in $\dA$. Our main theorem in this chapter then
essentially closes the gap between labelings of a globular graph $G$ in some
simplicial category $\dA$ and simplicial functors $\dC[\Sigma_{\min}^c] \to
\dA$.
\newpage
\labelingthm*
The proof of this theorem relies on a description of $\dC[\Sigma_{\min}^c]$ in
terms of products of low-dimensional simplices in $\dC[\Sigma_{\min}^c](x,y)$.
After having given the basic definitions sketched above in \autoref{sec:scat},
the whole of \autoref{sec:minimal-complete} is devoted to this description of
$\dC[\Sigma_{\min}^c]$.  The chapter then ends with a proof of
\autoref{thm:labeling} in \autoref{sec:labelings}.

\section{The simplicial category associated with a complete pasting diagram}\label{sec:scat}
This short section introduces the simplicial category $\dC[\Sigma]$
associated with a complete pasting diagram $\Sigma$. Moreover, we give the definition of a labeling of a globular graph in some simplicial category $\dA$.

\begin{definition}\label{def:c-of-sigma}
  Let $\Sigma$ be a complete pasting diagram. We
  associate with $\Sigma$ a simplicial category $\dC[\Sigma]$ with objects the
  vertices of $\Sigma$ and mapping spaces $ \dC[\Sigma](x,y) =
  \nerve(\Sigma_{x,y})$. The identities are given by the isomorphisms 
  \begin{equation*}
    \Delta^0 \iso \nerve(\Sigma_{x,x}) = \dC[\Sigma](x,x)
  \end{equation*}
  and the composition laws are given by
  \begin{equation*}
    \begin{tikzpicture}[diagram]
      \matrix[objects] {
        |(cxyz)| \dC[\Sigma](x,y) \times \dC[\Sigma](y,z) \& \& |(cxz)| \dC[\Sigma](x,z) \\[-4ex]
      |(nxyz)| \nerve(\Sigma_{x,y}) \times \nerve(\Sigma_{y,z}) \& 
      |(join)| \nerve( \Sigma_{x,y} \join \Sigma_{y,z} ) \& 
      |(nxz)|  \nerve(\Sigma_{x,z}) \smash{.} \\
      };
      \path[maps,->] 
        (cxyz) edge (cxz)
        (nxyz) edge node[above] {$\sim$} (join)
        (join) edge node[below] {inclusion} (nxz)
      ;
      \node[anchor=center,rotate=90] at ($ (cxyz) ! 0.5 ! (nxyz) $) {$=$};
      \node[anchor=center,rotate=90] at ($ (cxz) ! 0.5 ! (nxz) $) {$=$};
    \end{tikzpicture}
  \end{equation*}
\end{definition}

\begin{remark}\label{rem:csigma-well-defined}
  The fact that $\dC[\Sigma]$ is well-defined, i.\,e.\ that the composition
  laws are associative and unital, follows immediately from the fact that
  $\nerve$ is a strong monoidal functor, see \autoref{lem:nerve-of-join}.
\end{remark}
\begin{remark}
  The necessity of the condition that $\Sigma$ be complete in
  \autoref{def:c-of-sigma} is somewhat subtle. In fact, the isomorphism
  \begin{equation*}
    \nerve(\Sigma_{x,y}) \times \nerve(\Sigma_{y,z}) \to \nerve(\Sigma_{x,y}
  \join \Sigma_{y,z})
  \end{equation*}
  exists for all $\Sigma$ that are generated by wide subgraphs. However, as
  already pointed out in \autoref{rem:complete-and-join-decomposition}, for non-complete $\Sigma$ it might
  very well happen that $\Sigma_{x,y} \join \Sigma_{y,z} \nsubseteq
  \Sigma_{x,z}$, i.\,e.\ that the composition laws of $\dC[\Sigma]$ are not
  well-defined.
\end{remark}

\begin{remark}
  The assignment $\Sigma \mapsto \dC[\Sigma]$ is obviously functorial in
  arbitrary inclusions of complete pasting diagrams, that is, we have a functor
  from the category of complete pasting diagrams and inclusions into the
  category of simplicial categories.
\end{remark}

\begin{example}\label{rem:dc-maximal-pd}
  Let us compute the category $\dC[\Pi_{\max}]$ for the maximal pasting diagram
  on some globular graph $G$.  The objects of $\dC = \dC[\Pi_{\max}]$ are the
  vertices of $G$ and the mapping spaces are $\dC(x,y) = \nerve(\Pi_{\max,x,y})
  = \nerve(G_{x,y})$.  According to \autoref{prop:nerve-as-hom}, we can
  identify $\nerve(G_{x,y})$ with the nerve of the category $F_2G(x,y)$, where
  $F_2G$ denotes the free $2$-category on the globular graph $G$ considered as
  a $2$-computad. Altogether, we thus find that $\dC[\Pi_{\max}]$ has objects
  the vertices of $G$ and mapping spaces 
  \begin{equation*}
   \dC[\Pi_{\max}](x,y) = \nerve\bigl(
  F_2G(x,y) \bigr).
  \end{equation*}
  We leave it to the reader to verify that the compositions in
  $\dC[\Pi_{\max}]$ are induced by those of $F_2G$. Altogether, we thus find
  that $\dC[\Pi_{\max}]$ is obtained from $F_2G$ by applying the nerve functor
  $\cat \to \sset$ locally to each of the categories $F_2G(x,y)$.
\end{example}

\begin{remark}\label{rem:dc-graphical-description}
    The composition laws in the category $\dC[\Sigma]$ for some complete
    pasting diagram $\Sigma= (G,\sS)$ admit a neat description in terms of
    $n$-marked subgraphs. We know from \autoref{cor:pictorial-pd-nerve} that
    $n$-simplices in $\dC[\Sigma](x,y)$ correspond to admissible $n$-marked
    subgraphs $(P_\sigma,\lambda_\sigma)$ with $P_\sigma \in \sS_{x,y}$.
    Moreover, \autoref{rem:pictorial-join-isos} tells us that the composition
    $\sigma \comp \tau$ of two such simplices is nothing but $(P_\tau \join
    P_\sigma, \lambda_\tau\join \lambda_\sigma)$, where $\lambda_\tau \join
    \lambda_\sigma$ is given by $\lambda_\tau$ on $P_\tau$ and by
    $\lambda_\sigma$ on $P_\sigma$.
\end{remark}

\begin{remark}\label{rem:simplicial-computad}
  The categories $\dC[\Sigma]$ are simplicial computads in the sense of
  \cite{riehl-verity;homotopy-coherent-adjunctions-and-the-formal-theory-of-monads},
  that is, any $n$-simplex $\sigma \in \dC[\Sigma](x,y)$ has a unique
  decomposition $\sigma = (\sigma_1 \alpha_1) \comp \dots \comp (\sigma_a
  \alpha_a)$, where the $\sigma_i$ are nondegenerate atomic $n_i$-simplices and
  the $\alpha_i$ are degeneracy operators in $\dDelta$. Here, an atomic
  $n$-simplex is a simplex that cannot be written as a nontrivial composition
  in $\dC[\Sigma]$. We will use the decomposition of $\sigma$ into atomic
  $n$-simplices more or less explicitly in our proof of
  \autoref{thm:labeling}.
  
  The above decomposition of some simplex $\sigma$ is easy to get hold of using
  the pictorial description of $\dC[\Sigma]$ given in the preceding remark. It
  is clear that an $n$-simplex $\sigma$ is atomic if and only if $P_\sigma$ is
  $2$-connected or a single edge.  Given any admissible $n$-marked subgraph
  $(P_\sigma,\lambda_\sigma)$ in $\dC[\Sigma](x,y)$, we therefore write
  $P_\sigma = P_1 \join \dots \join P_a$ with each $P_i$ a $2$-connected
  globular subgraph and thus obtain $\sigma = \tau_a \comp \dots \comp \tau_1$
  with $(P_{\tau_i},\lambda_{\tau_i}) = (P_i, \lambda_\sigma|_{P_i})$. The
  Eilenberg-Zilber lemma now yields nondegenerate simplices $\sigma_i$ and
  degeneracy operators $\alpha_i$ with $\tau_i = \sigma_i \alpha_i$ and it is
  clear from the description of the simplicial operators in
  \autoref{prop:marked-simplicial} that $\sigma_i$ is atomic, too.
\end{remark}

The following definition formalises the notion of a $2$-dimensional diagram in
a simplicial category. 
\begin{definition}\label{def:labeling}
  Let\/ $G$ be a globular graph. A\/ \emph{labeling} $\Lambda$ of\/ $G$ in some
  simplicial category $\dA$ consists of the following data:
  \begin{enumerate}
    \item An object $\Lambda x \in \dA$ for each vertex $x \in G$.
    \item A\/ $0$-simplex $\Lambda e \in \dA(\Lambda x,\Lambda y)$ for each edge $e$ from $x$ to $y$ in $G$.
    \item A\/ $1$-simplex $\Lambda \phi \in \dA(\Lambda x,\Lambda y)$ for each interior face\/
      $\phi$ of\/ $G$ with\/ $s(\phi) =x$ and $t(\phi) = y$. 
  \end{enumerate}
  This data is subject to the condition that for each interior face $\phi$ of\/ $G$ we have equalities
  \begin{equation*}
    d_0 \Lambda \phi = \Lambda e_r \comp \dots \comp \Lambda e_1
    \quad\text{and}\quad
    d_1 \Lambda \phi = \Lambda f_s \comp \dots \comp \Lambda f_1,
  \end{equation*}
  where $\dom(\phi) = f_1 \cdot \dots \cdot f_s$ and $\cod(\phi) = e_1\cdot \dots \cdot e_r$.
\end{definition}

\begin{example}
  \begin{enumerate}[label=(\alph*)]
    \item 
  A labeling $\Lambda$ of the graph $B_n$ from \autoref{ex:nerve-bn} in some
  simplicial category $\dA$ consists of the choice of two objects $\Lambda s,
  \Lambda t \in \dA$ and a map $\Delta^1 \join \dots \join \Delta^1 \to
  \dA(\Lambda s, \Lambda t)$ of simplicial sets. In the case that $\dA$ is
  enriched over quasi-categories, a labeling of $B_n$ is therefore nothing but
  a string of $n$ composable cells in $\dA(\Lambda s,\Lambda t)$.
\item A labeling $\Lambda$ of the graph $B_1 \join B_1$ in
  \autoref{fig:b1-times-b1} in some simplicial category $\dA$ consists of the
      choice of three vertices $\Lambda x$, $\Lambda y$ and $\Lambda z$ in
      $\dA$ together with two maps $\Delta^1 \to \dA(\Lambda x, \Lambda y)$ and
      $\Delta^1 \to \dA(\Lambda y, \Lambda z)$. In anticipation of \autoref{sec:minimal-complete}, the reader might want to
      compare the notion of a labeling of $B_1 \join B_1$ with that of a
      functor $\dC[\Sigma_{\min}^c] \to \dA$ from the simplicial category
      associated with the minimal complete pasting diagram on $B_1 \join B_1$.
  \end{enumerate}
\end{example}

\begin{remark}\label{rem:labeling}
  Each functor $u\colon \dC[\Sigma] \to \dA$ determines a labeling $\Lambda_u$
  of the graph $G$ underlying $\Sigma$ in $\dA$. The labeling is given by
  $\Lambda_u x = u(x)$, $\Lambda_u e = u(e)$ and $\Lambda_u \phi = u(\phi)$,
  where $e$ and $\phi = (\dom \phi \leq \cod \phi)$ are considered as $0$- and
  $1$-simplices in $\dC[\Sigma](s(e),t(e))$ and $\dC[\Sigma](s(\phi),t(\phi))$,
  respectively.
\end{remark}

\section{The simplicial category associated with a minimal complete pasting diagram}\label{sec:minimal-complete}
In this section we prepare the ground for the proof of \autoref{thm:labeling}
in the following section and give an explicit description of the category
$\dC[\Sigma]$ in the case that $\Sigma$ is the minimal complete pasting diagram
on some globular graph $G$. In a certain sense, we show that $\dC[\Sigma]$ is
freely generated by $G$ and additional witnesses of the Godement interchange
law, that is, any simplex $\sigma \in \dC[\Sigma](x,y)$ of dimension $n \geq 2$
sits inside some cube $(\Delta^1)^a \subseteq \dC[\Sigma](x,y)$ that
corresponds to the different orders of composition of a diagram such as
\begin{equation*}
  \begin{tikzpicture}[diagram]
    \matrix[objects] {
      |(1)| \bullet \& |(2)| \bullet \& |(3)| \bullet\&[-2em] \cdots \&[-2em] |(4)| \bullet \& |(5)| \bullet\\
    };
    \path[maps,->] 
      (1) edge[bend left] (2)
      (1) edge[bend right] (2)
      (2) edge[bend left] (3)
      (2) edge[bend right] (3)
      (4) edge[bend left] (5)
      (4) edge[bend right] (5)
    ;
    \node[left]  at (1) {$x$};
    \node[right]  at (5) {$y$};
  \end{tikzpicture}
\end{equation*}
whose $a$ globs are faces of $G$, see
\autoref{prop:godement-representation-simplex}. The argument leading to this
presentation relies on the description of simplices in $\nerve(\Sigma)$ in
terms of admissible $n$-marked subgraphs $(P_\sigma,\lambda_\sigma)$ from
\autoref{cor:pictorial-pd-nerve}. 

This section has two parts. In the first part we prove the statement sketched
above, i.\,e.\ that any $n$-simplex $\sigma$ in $\dC(x,y)$ is contained in some
cube $(\Delta^1)^a$. In the second part we then determine the action of
simplicial operators $\alpha$ on the simplices of $\dC(x,y)$ in terms of the
maps $\Delta^n \to (\Delta^1)^a$. Especially this latter part is rather
technical and the reader might want to skip some of the proofs until she or he
has had a look at the proof of \autoref{thm:labeling} in the forthcoming
section, where these technicalities are crucial.

\newpage
Fix a globular graph $G$ and let $\Sigma = \Sigma_{\min}^c = (G,\sS)$ be the
minimal complete pasting diagram on $G$. Throughout this section we abbreviate
$\dC = \dC[\Sigma]$.  Recall from \autoref{lem:pictorial-desc-nerve-pds} that
there is a one-to-one-correspondence between $n$-simplices $\sigma \in
\nerve(\Sigma)$ and $n$-marked subgraphs $(P_\sigma,\lambda_\sigma)$, where
$P_\sigma \in \sS$ and $\lambda_\sigma \colon \Phi(P_\sigma) \to \{1,\dots,n\}$
is admissible in the sense of \autoref{def:admissible}. However, any $P_\sigma
\in \sS = \sS_{\min}$ is a join of edges and faces of $G$, so that
admissibility of $\lambda_\sigma$ turns out to be an empty condition in the
case at hand, because there are no interior faces $\phi$ and $\psi$ in
$P_\sigma$ with some edge $e$ in both $\cod(\phi)$ and $\dom(\psi)$.  We thus
identify $n$-simplices in $\nerve(\Sigma)$ or $\nerve(\Sigma_{x,y}) = \dC(x,y)$
with $n$-marked subgraphs $(P_\sigma,\lambda_\sigma)$, where $P_\sigma \in
\sS_{x,y}$ and $\lambda_\sigma \colon \Phi(P_\sigma) \to \{1,\dots,n\}$ is an
arbitrary map. 

Let us now consider an $n$-simplex $\sigma = (p_0 \leq \dots \leq p_n) \in \dC(x,y)$ and
let $(P_\sigma,\lambda_\sigma)$ be the corresponding $n$-marked subgraph.
Observe that $P_\sigma$ can be written uniquely as a join $P_\sigma = P_1 \join
\dots \join P_a$ with each $P_i$ an edge or a face of $P_\sigma$ since $\Sigma
= (G,\sS)$ is the minimal complete pasting diagram on $G$. 
With
\begin{equation*}
  \epsilon_i(\sigma) = \begin{cases}
    0 & \text{if } P_i \text{ is an edge of $P$},\\
    1 & \text{if } P_i \text{ is a face of $P$},
  \end{cases}
\end{equation*}
each $P_i$ corresponds to an $\epsilon_i(\sigma)$-dimensional simplex $\sigma_i
\in \dC(sP_i,tP_i)$, namely $\sigma_i = (P_i)$ if $P_i$ is an edge and
$\sigma_i = (\dom P_i \leq \cod P_i)$ if $P_i$ is a face. 
We thus obtain a map
\begin{equation*}
  \Delta^{\epsilon_1(\sigma)} \times \Delta^{\epsilon_2(\sigma)} \times \dots \times
  \Delta^{\epsilon_a(\sigma)} \xto{(\sigma_1,\dots,\sigma_a)} 
  \dC(sP_1,tP_1) \times \dots \times \dC(sP_a,tP_a).
\end{equation*}
As this map is crucial to the rest of this paragraph, we introduce the
abbreviations 
\begin{equation}\label{eq:delta-eps-def}
  \Delta^{\epsilon(\sigma)} = \Delta^{\epsilon_1(\sigma)} \times
  \Delta^{\epsilon_2(\sigma)} \times \dots \times \Delta^{\epsilon_a(\sigma)}
\end{equation}
and
\begin{equation}\label{eq:dc-sigma-def}
  \dC(\sigma) = \dC(sP_1,tP_1) \times \dots \times \dC(sP_a,tP_a).
\end{equation}
\begin{remark}
  The reader should note that the above discussion is merely a special case of
  the decomposition $\sigma = (\sigma_a \alpha_a) \comp \dots \comp
  (\sigma_1\alpha_1)$ of an $n$-simplex that we already sketched in
  \autoref{rem:simplicial-computad}. In fact, the decomposition $P_\sigma = P_1
  \join \dots \join P_a$ of $P_\sigma$ into faces and edges corresponds to a
  decomposition $\sigma = \tau_a \comp \dots \comp \tau_1$ of $\sigma$ into
  atomic $n$-simplices. Moreover, writing $\tau_i = \sigma_i \alpha_i$ with
  $\sigma_i$ nondegenerate recovers the construction of the simplices
  $\sigma_i$ given above.
\end{remark}
\begin{proposition}\label{prop:godement-representation-simplex}
  Let $\sigma \in \dC(x,y)$ be an arbitrary $n$-simplex. There then
  exists a unique map\/ $\what{\sigma}\colon \Delta^n \to \Delta^{\epsilon(\sigma)}$
  such that 
  \begin{equation*}
    \begin{tikzpicture}[diagram]
      \matrix[objects] {
        |(dp)| \Delta^{\epsilon(\sigma)} \&[+3em] |(sp)| \dC(\sigma)  \\
        |(dn)| \Delta^n \& |(sxy)| \dC(x,y) \\
      };
      \path[maps,->] 
        (dp) edge node[above] {$ (\sigma_1, \dots, \sigma_a) $} (sp)
        (dn) edge node[left] {$ \what{\sigma} $} (dp)
        (dn) edge node[below] {$ \sigma $} (sxy)
        (sp) edge node[right] {composition} (sxy)
      ;
    \end{tikzpicture}
  \end{equation*}
  commutes. 
\end{proposition}
\begin{proof}
  We keep the notation fixed that we use throughout this section. That
  means in particular that $\sigma = (p_0 \leq \dots \leq p_n)$ has
  $(P_\sigma,\lambda_\sigma)$ as its associated $n$-marked subgraph and that $P_\sigma = P_1
  \join \dots \join P_a$ is the decomposition of $P_\sigma$ into edges and
  faces of $P$. Moreover, $\sigma_i = (P_i)$ if $P_i$ is an edge and $\sigma_i
  = (\dom P_i \leq \cod P_i)$ if $P_i$ is a face.

  Let us consider an arbitrary map $\what{\sigma} \colon \Delta^n \to
  \Delta^{\epsilon(\sigma)}$.  By the universal property of products, the map
  $\what{\sigma}$ is uniquely determined by the compositions
  \begin{equation*}
    \beta_i \colon \Delta^n \to \Delta^{\epsilon(\sigma)} \to \Delta^{\epsilon_i(\sigma)},
  \end{equation*}
  that is, by maps $\beta_i \colon [n] \to [\epsilon_i(\sigma)]$.  The composition
  \begin{equation*} 
    (\sigma_1, \dots, \sigma_a) \comp \what{\sigma}
  \end{equation*}
  thus classifies the $n$-simplex $( \sigma_1 \beta_1, \dots, \sigma_a
  \beta_a)$ of $\dC(\sigma)$. Given the explicit description of the simplices
  $\sigma_i$ at the beginning of this proof, one now easily computes the
  $n$-simplices $\tau_i = \sigma_i \beta_i$ as $\tau_i = (q_0^i \leq \dots \leq
  q_n^i)$ with
  \begin{equation*}
    q_j^i = \begin{cases}
      P_i & \text{if } \epsilon_i(\sigma) = 0, \\
      \dom P_i & \text{if } \epsilon_i(\sigma) = 1 \text{ and } \beta_i(j) = 0,\\
      \cod P_i & \text{if } \epsilon_i(\sigma) = 1 \text{ and } \beta_i(j) = 1.
    \end{cases}
  \end{equation*}
  The composition $\tau = \tau_r \comp \ldots \comp \tau_1 = (q_0 \leq \dots \leq q_n) \in \dC(x,y)$ of
  these simplices $\tau_i$ therefore has
  \begin{equation*}
    q_j = q_j^1 \cdot \ldots \cdot q_j^r
  \end{equation*}
  by the definition of composition in $\dC$.
  Note that $\what{\sigma}$ renders the diagram 
  \begin{equation}\label{dia:what-sigma}
    \begin{tikzpicture}[diagram]
      \matrix[objects] {
        |(dp)| \Delta^{\epsilon(\sigma)} \&[+3em] |(sp)| \dC(\sigma) \\ 
        |(dn)| \Delta^n \& |(sxy)| \dC(x,y) \\
      };
      \path[maps,->] 
        (dp) edge node[above] {$ (\sigma_1, \dots, \sigma_r) $} (sp)
        (dn) edge node[left] {$ \what{\sigma} $} (dp)
        (dn) edge node[below] {$ \sigma $} (sxy)
        (sp) edge node[right] {composition} (sxy)
      ;
    \end{tikzpicture}
  \end{equation}
  commutative if and only if $\tau = \sigma$, that is, if and only if $q_j =
  p_j$ for all $0 \leq j \leq n$.  According to \autoref{rem:pi-characterisation}, the path
  $p_j$ occuring in $\sigma = (p_0 \leq \dots \leq p_n)$ can be characterised
  as the unique wide path in $P$ with exactly the faces $P_i$ with
  $\lambda_\sigma(P_i) \leq j$ to its left. However, given the above
  description of $q_j$ it is easy to see that a face $P_i$ lies to the left of
  $q_j$ if and only if $\beta_i(j) = 1$.  Altogether, we may thus conclude that
  the map $\what{\sigma}$ induced by the maps
  \begin{equation*}
    \beta_i(j) = \begin{cases}
      0 & \text{if }j < \lambda_\sigma(P_i) \\
      1 & \text{if }j \geq \lambda_\sigma(P_i)
    \end{cases}
  \end{equation*}
  is the unique map that renders \eqref{dia:what-sigma} commutative.
\end{proof}
\begin{corollary}\label{cor:godement-composition}
  Let $\sigma \in \dC(x,y)$ and $\tau \in \dC(y,z)$ be two composable $n$-simplices.
  The map $\widehat{\tau \comp \sigma}$ is the map 
  \begin{equation*}
    (\what{\sigma}, \what{\tau}) \colon \Delta^n \to 
    \Delta^{\epsilon(\sigma)} \times \Delta^{\epsilon(\tau)} = \Delta^{\epsilon(\tau \comp \sigma)}.
  \end{equation*}
\end{corollary}
\begin{proof}
  The composite $\tau \comp \sigma$ in $\dC$ corresponds to the join $(P_\sigma \join P_\tau,\lambda_\sigma \join \lambda_\tau)$ of $n$-marked subgraphs, see
  \autoref{rem:dc-graphical-description}. This implies 
  $\Delta^{\epsilon(\tau \comp \sigma)} = \Delta^{\epsilon(\sigma)} \times
  \Delta^{\epsilon(\tau)}$, $\dC(\tau \comp \sigma) = \dC(\sigma) \times
  \dC(\tau)$ and $( (\tau\sigma)_1, \dots, (\tau\sigma)_c) = (\sigma_1, \dots,
  \sigma_a, \tau_1, \dots, \tau_b)$.  The diagram
  \begin{equation*}
    \begin{tikzpicture}[diagram]
      \matrix[objects] {
        |(dprod)| \Delta^{\epsilon(\sigma)} \times \Delta^{\epsilon(\tau)} \&[+4em]
        |(cprod)| \dC(\sigma) \times \dC(\tau) \&[+2em] \\
        |(dn)| \Delta^n \& |(c2)| \dC(x,y) \times \dC(y,z) \& |(c)| \dC(x,z) \\
      };
      \path[maps,->] 
        (dn) edge node[left] {$(\what{\sigma}, \what{\tau})$} (dprod)
        (dprod) edge node[above] {$(\sigma_1,\dots,\sigma_a,\tau_1,\dots,\tau_b)$} (cprod)
        (cprod) edge node[above right] {composition} (c)
        (cprod) edge node[left] {composition $\times$ composition} (c2)
        (c2) edge node[below] {composition} (c)
        (dn) edge node[below] {$(\sigma,\tau)$} (c2)
      ;
    \end{tikzpicture}
  \end{equation*}
  commutes by associativity of composition in $\dC$ and the definition of
  $\what{\sigma}$ and $\what{\tau}$.  It now follows from the uniqueness
  asserted in \autoref{prop:godement-representation} that $\widehat{\tau \comp
  \sigma} = (\what{\tau}, \what{\sigma})$.
\end{proof}

Let us now consider a simplicial operator $\alpha \colon [m] \to [n]$ and the
$n$-marked subgraph $(P_{\sigma\alpha},\lambda_{\sigma\alpha})$ associated with
$\sigma \alpha$. The steps given in \autoref{prop:marked-simplicial} to
compute $(P_{\sigma\alpha},\lambda_{\sigma\alpha})$ from
$(P_\sigma,\lambda_\sigma)$ then admit some simplifications because any edge of
$P_\sigma, P_{\sigma\alpha} \in \sS = \sS_{\min}^c$ is incident with the
exterior face. More precisely, one obtains the following procedure to determine
$(P_{\sigma\alpha},\lambda_{\sigma\alpha})$ in terms of $(P_\sigma,\lambda_\sigma)$:
  \begin{enumerate}
    \item Remove the edges and interior vertices of all paths $\dom(\phi)$, where $\phi$ is some interior face of $P_\sigma$ with $\lambda_\sigma(\phi) \leq \alpha(0)$.
    \item Remove the edges and interior vertices of all paths $\cod(\phi)$, where $\phi$ is some interior face of $P_\sigma$ with $\lambda_\sigma(\phi) > \alpha(m)$.
    \item Define $\lambda_{\sigma\alpha}$ on the remaining graph by
      \begin{equation*}
        \lambda_{\sigma\alpha}(\phi) = \min\bigl\{\, k \in \{1,\dots,m\} \mid \alpha(k) \geq \lambda_\sigma(\phi) \bigr\}.
      \end{equation*}
  \end{enumerate}

  As both $P_\sigma$ and $P_{\sigma\alpha}$ are wide subgraphs of $G_{x,y}$, the following 
description of the decomposition $P_{\sigma\alpha} = Q_1 \join \dots \join Q_b$ into edges and faces of $P_{\sigma\alpha}$ in terms of $P_\sigma = P_1 \join \dots \join P_a$ is immediate.
\begin{lemma}\label{lem:decomposition-description}
  Let $\sigma$ be an $n$-simplex in $\dC(x,y)$ and let $P_\sigma = P_1 \join
  \dots \join P_a$ be the decomposition of its associated $n$-marked subgraph
  into edges and faces. Further let $\alpha \colon [m] \to [n]$ be a simplicial
  operator in $\dDelta$ and let $P_{\sigma\alpha} = Q_1 \join \dots \join Q_b$
  be the decomposition of $P_{\sigma\alpha}$ into edges and faces.  For any $i
  \in \{1,\dots,a\}$ we then have the following relations between these
  decompositions:
  \begin{enumerate}
    \item If\/ $P_i$ is an edge, then $P_i = Q_j$ for some suitable $j$.
    \item If\/ $P_i$ is an interior face $\phi$ with\/ $\alpha(0) <
      \lambda_\sigma(\phi) \leq \alpha(m)$, then $P_i = Q_j$ for some suitable
      $j$. Moreover, $k < \lambda_{\sigma\alpha}(\phi)$ if and only if\/
      $\alpha(k) < \lambda_\sigma(\phi)$.
    \item If\/ $P_i$ is an interior face with\/ $\lambda_\sigma(P_i) \leq \alpha(0)$, then $\cod P_i = Q_j \join \dots \join Q_{j+k}$ for suitable $j$ and $k$.
    \item If\/ $P_i$ is an interior face with\/ $\lambda_\sigma(P_i) > \alpha(m)$, then $\dom P_i = Q_j \join \dots \join Q_{j+k}$ for suitable $j$ and $k$.
  \end{enumerate}
\end{lemma}

Let us fix the notation of \autoref{lem:decomposition-description} for the
remainder of this section. In order to be able to describe the action of
simplicial operators on the simplices of $\dC(x,y)$ in terms of the maps
$\what{\sigma}$ and $\widehat{\sigma\alpha}$, we need to relate
$\Delta^{\epsilon(\sigma)}$ and $\Delta^{\epsilon(\sigma\alpha)}$ in some
sensible way. 
Recall from \eqref{eq:delta-eps-def} that 
\begin{equation*}
  \Delta^{\epsilon(\alpha)} = \Delta^{\epsilon_1(\alpha)} \times \dots \times
  \Delta^{\epsilon_a(\alpha)}.
\end{equation*}
From the description of $(P_{\sigma\alpha},\lambda_{\sigma\alpha})$ given in
\autoref{lem:decomposition-description} we fabricate a map $\epsilon(\alpha)
\colon \Delta^{\epsilon(\sigma\alpha)} \to \Delta^{\epsilon(\sigma)}$ whose
component $\epsilon(\alpha)_i$ at $\Delta^{\epsilon_i(\sigma)}$ is the identity
$\Delta^{\epsilon_j(\sigma\alpha)} = \Delta^{\epsilon_i(\sigma)}$ if $P_i =
Q_j$ for some $j$ and whose components $\epsilon(\alpha)_i$ at
$\Delta^{\epsilon_i(\sigma)}$ with $P_i \nsubseteq Q$ are given by
\begin{equation*}
  \begin{cases}
    \Delta^{\epsilon_j(\sigma\alpha)} \times \dots \times \Delta^{\epsilon_{j+k}(\sigma\alpha)} \iso \Delta^0 \xto{d_0} \Delta^{\epsilon_i(\sigma)} & \text{if } Q_j \join \dots \join Q_{j+k} = \cod P_i, \\
    \Delta^{\epsilon_j(\sigma\alpha)} \times \dots \times \Delta^{\epsilon_{j+k}(\sigma\alpha)} \iso \Delta^0\xto{d_1} \Delta^{\epsilon_i(\sigma)} & \text{if } Q_j \join \dots \join Q_{j+k} = \dom P_i. \\
  \end{cases}
\end{equation*}

\begin{remark}\label{rem:epsilon-nat}
  The maps $\epsilon(\alpha)$ are functorial in $\alpha$ in the sense that we
  have $\epsilon(\alpha \comp \beta) = \epsilon(\alpha) \comp \epsilon(\beta)$
  for all composable simplicial operators $\alpha$ and $\beta$.
\end{remark}

We record the following technical remark in order to isolate all technical
properties of the  maps $\epsilon(\alpha)$ from the actual proof of
\autoref{thm:labeling} in the forthcoming section. However, this remark might
also provide the reader with some intuition for the whys and wherefores of the
maps $\epsilon(\alpha)$.
\begin{remark}\label{rem:epsilon-lambda}
  Suppose that $\Lambda$ is a labeling of $\Sigma = \Sigma_{\min}^c$ in some
  simplicial category $\dA$. Further suppose that $\sigma$ is some $n$-simplex
  in $\dC(x,y)$ and that $\alpha\colon [m] \to [n]$ is a simplicial operator in
  $\dDelta$. If we let 
  \begin{equation*}
    \dA(\sigma) = \dA(\Lambda sP_1,\Lambda tP_1) \times \dots \times \dA(\Lambda sP_a, \Lambda tP_a),
  \end{equation*}
  then there are obvious partial composition maps
  \begin{equation*}
    \varpi \colon \dA( \sigma\alpha ) \to \dA(\sigma)
  \end{equation*}
  given by the identities on those $\dA(\Lambda sP_i,\Lambda tP_i)$ with $P_i = Q_j$ for some
  $j$ and by the compositions
  \begin{equation*}
    \dA(\Lambda sQ_j, \Lambda tQ_j) \times \dots \times \dA( \Lambda sQ_{j+k}, \Lambda tQ_{j+k} ) \to \dA(\Lambda sQ_j,\Lambda tQ_{j+k})
  \end{equation*}
  whenever $P_i$ is a face and $Q_j \join \dots \join  Q_{j+k}$ is either $\dom
  P_i$ or $\cod P_i$.

  Given the relation between $P_\sigma = P_1 \join \dots \join P_a$ and $P_{\sigma\alpha} = Q_1 \join \dots \join Q_{b}$ in \autoref{lem:decomposition-description}, it is immediate that
  \begin{equation*}
    \varpi\bigl( \Lambda (\sigma\alpha)_j \bigr)_i
    = \begin{cases}
      \Lambda \sigma_i & \text{if } P_i = Q_j \text{ for some } j, \\
      \Lambda Q_{j+k} \comp \dots \comp \Lambda Q_{j} & \text{if } \dom P_i = Q_j \join \dots \join Q_{j+k} \\
      & \text{or } \cod P_i = Q_j \join \dots \join Q_{j+k}.
    \end{cases}
  \end{equation*}
  With our definition $\sigma_i = (\dom P_i \leq \cod P_i)$ and the definition of a labeling $\Lambda$, one then sees that 
  \begin{align*}
    \varpi\bigl( \Lambda (\sigma\alpha)_j \bigr)_i
    &= 
      \begin{cases}
      \Lambda \sigma_i & \text{if } P_i = Q_j \text{ for some } j, \\
      d_1( \Lambda \sigma_i ) & \text{if } \dom P_i = Q_j \join \dots \join Q_{j+k}, \\
      d_0( \Lambda \sigma_i ) & \text{if } \cod P_i = Q_j \join \dots \join Q_{j+k}.
      \end{cases}
  \end{align*}

  Moreover, it follows from our construction of $\epsilon(\alpha)$ from above
  that $\sigma_i \comp \epsilon(\alpha)$ classifies the simplex given by
  \begin{equation*}
    \Lambda \sigma_i \cdot \epsilon(\alpha)_i  
    = \begin{cases}
      \Lambda \sigma_i & \text{if } P_i = Q_j \text{ for some } j,\\
      d_1( \Lambda \sigma_i )& \text{if } \dom P_i = Q_j \join \dots \join Q_{j+k},\\
      d_0( \Lambda \sigma_i )& \text{if } \cod P_i = Q_j \join \dots \join Q_{j+k}.
    \end{cases}
  \end{equation*}
  Altogether, this proves that the diagram
  \begin{equation*}
  \begin{tikzpicture}[diagram]
      \matrix[objects] {
        |(dsa)| \Delta^{\epsilon(\sigma\alpha)} \&[+6em] |(sp)| \dA(\sigma\alpha) \\[-2ex]
        |(ds)| \Delta^{\epsilon(\sigma)}  \& |(sq)| \dA(\sigma)\\
      };
      \path[maps,->] 
      (dsa) edge node[left] {$\epsilon(\alpha)$} (ds)
      (dsa) edge node[above] {$\bigl( \Lambda(\sigma\alpha)_1, \dots, \Lambda (\sigma\alpha)_b \bigr)$} (sp)
      (ds) edge node[below] {$( \Lambda\sigma_1, \dots, \Lambda\sigma_a)$} (sq)
      (sp) edge node[right] {$\varpi$} (sq)
      ;
    \end{tikzpicture}
  \end{equation*}
  commutes.
\end{remark}

Let us now give a description of the action of simplicial operators on the
simplices of $\dC(x,y)$ in terms of the maps $\what{\sigma}$.
\begin{proposition}\label{prop:godement-representation}
  The diagram
  \begin{equation*}
    \begin{tikzpicture}[diagram]
      \matrix[objects] {
        |(dm)| \Delta^m \& |(dsa)| \Delta^{\epsilon(\sigma\alpha)} \\
        |(dn)| \Delta^n \& |(ds)| \Delta^{\epsilon(\sigma)}  \\
      };
      \path[maps,->] 
      (dm) edge node[above] {$\widehat{\sigma\alpha}$} (dsa)
      (dm) edge node[left] {$\alpha$} (dn)
      (dn) edge node[below] {$\what{\sigma}$} (ds)
      (dsa) edge node[right] {$\epsilon(\alpha)$} (ds)
      ;
    \end{tikzpicture}
  \end{equation*}
  commutes for any simplicial operator $\alpha\colon [m] \to [n]$ in $\dDelta$ and any
  $n$-simplex $\sigma \in \dC(x,y)$.
\end{proposition}
\begin{proof}
  We still keep the notation from \autoref{lem:decomposition-description}.
  According to the proof of \autoref{prop:godement-representation-simplex}, the
  composition
  \begin{equation*}
    \Delta^m \xto{\alpha} \Delta^n \xto{\what{\sigma}} \Delta^{\epsilon(\sigma)} \to \Delta^{\epsilon_i(\sigma)}
  \end{equation*}
  corresponds to the map $[m] \to [0]$ if $\epsilon_i(\sigma) = 0$ and to $\beta_i \colon [m] \to [1]$ given by
  \begin{equation*}
    \beta_i(j) = \begin{cases}
      0 & \text{if } \alpha(j) < \lambda_\sigma(P_i), \\
      1 & \text{if } \alpha(j) \geq \lambda_\sigma(P_i), 
    \end{cases}
  \end{equation*}
  if $\epsilon_i(\sigma) = 1$.
  In order to show that the diagram in the lemma commutes, it
  suffices to show that these maps coincide with the maps $\gamma_i\colon [m] \to
  [\epsilon_i(\sigma)]$ corresponding to the compositions 
  \begin{equation}\label{eq:long-composition}
    \Delta^m \xto{\widehat{\alpha\sigma}} \Delta^{\epsilon(\sigma\alpha)} \xto{\epsilon(\alpha)} \Delta^{\epsilon(\alpha)} \to \Delta^{\epsilon_i(\alpha)}.
  \end{equation}
  To this end, we distinguish the following cases that correspond to the
  different cases in \autoref{lem:decomposition-description}:
  \begin{enumerate}
    \item If $\epsilon_i(\sigma) = 0$, i.\,e.\ if $P_i$ is an edge, then
      $[\epsilon_i(\alpha)] = [0]$ is terminal and we are done.
    \item If $\epsilon_i(\sigma) = 1$, i.\,e.\ if $P_i$ is an interior face and
      $\alpha(0) < \lambda_\sigma(P_i) \leq \alpha(m)$, then $P_i = Q_k$ for
      some $k$. The component $\epsilon(\alpha)_i$ is then the identity on
      $\Delta^{\epsilon_k(\sigma\alpha)} = \Delta^{\epsilon_i(\alpha)}$. The
      map $[m] \to [1]$ corresponding to the composition
      \eqref{eq:long-composition} is therefore given by
      \begin{equation*}
        \gamma_i(j) = \begin{cases}
          0 & \text{if } j < \lambda_{\sigma\alpha}(P_i) \\
          1 & \text{if } j \geq \lambda_{\sigma\alpha}(P_i). 
        \end{cases}
      \end{equation*}
      We conclude by the equivalence 
      \begin{equation*}
        j < \lambda_{\sigma\alpha(P_i)} \quad\text{if and only if} \quad \alpha(j) < \lambda_\sigma(P_i)
      \end{equation*}
      stated in \autoref{lem:decomposition-description}.
    \item If $\epsilon_i(\sigma) = 1$ and $\lambda_\sigma(P_i) \leq \alpha(0)$,
      then $\cod P_i = Q_k \join \dots \join Q_{k+l}$ and $\epsilon(\alpha)_i$
      is the map
      \begin{equation*}
        \Delta^{\epsilon_k(\sigma\alpha)} \times \dots \times
        \Delta^{\epsilon_{k+l}(\sigma\alpha)} \xto{d_0}
        \Delta^{\epsilon_i(\sigma)}
      \end{equation*}
      The map $[m] \to [1]$ corresponding to the composition
      \eqref{eq:long-composition} is hence the constant map $1$. This coincides
      with $\beta_i$ since $\lambda_\sigma(P_i) \leq \alpha(0)$, i.\,e.\
      $\alpha(j) \geq \lambda_{\sigma}(P_i)$ for all $j \in [m]$.
    \item The remaining case that $\epsilon_i(\sigma) = 1$ and
      $\lambda_\sigma(P_i) > \alpha(m)$ can be handled as the preceding case.
  \end{enumerate}
  This finishes the proof of \autoref{prop:godement-representation}.

\end{proof}

\section{Labelings of Globular Graphs}\label{sec:labelings}

This section solely consists of a proof of 
\labelingthm
\begin{proof}
  Let us begin by introducing some notation and recalling some facts from
  \autoref{sec:minimal-complete}. Given a simplex $\sigma$ in
  $\dC[\Sigma](x,y)$ with associated $n$-marked subgraph
  $(P_\sigma,\lambda_\sigma)$, we decompose 
  \begin{equation*}
    P_\sigma = P_1 \join \dots \join
  P_a
  \end{equation*}
  with each $P_i$ an edge or an interior face of $P_\sigma$, see
  \autoref{sec:minimal-complete}. We then have simplices $\sigma_i = (P_i)$ if
  $i$ is an edge and $\sigma_i = (\dom P_i \leq P_i)$ if $P_i$ is an interior
  face. Similar to the abbreviation 
  \begin{equation*}
    \dC[\Sigma](\sigma) = \dC[\Sigma](sP_1,tP_1) \times \dots \times
    \dC[\Sigma](sP_a,tP_a)
  \end{equation*}
  introduced in \eqref{eq:dc-sigma-def} above, we write
  \begin{equation*}
    \dA(\sigma) = \dA(\Lambda sP_1, \Lambda tP_1) \times \dots \times \dA(
    \Lambda sP_a, \Lambda tP_a).
  \end{equation*}

  Let us now show that a simplicial functor $u\colon \dC[\Sigma] \to \dA$ is
  uniquely determined by its associated labeling $\Lambda_u$. It is obvious
  that the labeling $\Lambda_u$ completely determines the action of $u$ on
  objects.  Moreover, $\Lambda_u$ also determines the action of $u$ an all the
  simplices $\sigma \in \dC[\Sigma](x,y)$ of the form $\sigma = (e)$ for an
  edge $e$ or $\sigma = (\dom \phi \leq \cod \phi)$ for an interior face
  $\phi$.  Now consider some arbitrary $n$-simplex $\sigma \in
  \dC[\Sigma](x,y)$.  We then have a commutative diagram
  \begin{equation*}
    \begin{tikzpicture}[diagram]
      \matrix[objects] {
        |(dp)| \Delta^{\epsilon(\sigma)} \&[+3em] |(sp)| \dC[\Sigma](\sigma) \&[+2em] |(ap)| \dA(\sigma) \\[-1em]
        |(dn)| \Delta^n \& |(sxy)| \dC[\Sigma](x,y) \& |(axy)| \dA(ux,uy) \\
      };
      \path[maps,->] 
        (dp) edge node[above] {$ (\sigma_{1}, \dots, \sigma_{r}) $} (sp)
        (dn) edge node[left] {$ \what{\sigma} $} (dp)
        (dn) edge node[below] {$ \sigma $} (sxy)
        (sp) edge node[above] {$(u, \dots, u)$} (ap)
        (ap) edge node[right] {composition} (axy)
        (sxy) edge node[below] {$u$} (axy)
        (sp) edge node[left] {composition} (sxy)
      ;
    \end{tikzpicture}
  \end{equation*}
  by \autoref{prop:godement-representation-simplex} and functoriality of $u$.
  The top row in this diagram is 
  \begin{equation*}
    (u,\dots,u) \comp (\sigma_{1}, \dots, \sigma_{r}) 
    = (u \sigma_{1}, \dots, u \sigma_{r} ) 
    = \bigl(\Lambda_u P_1 , \dots, \Lambda_u P_r\bigr).
  \end{equation*}
  The image $u(\sigma)$ of $\sigma$ is hence the $n$-simplex in $\dA(ux,uy)$
  classified by the composition
  \begin{equation*}
    \begin{tikzpicture}[diagram]
      \matrix[objects] {
        |(dn)| \Delta^n \&[-1em] |(dp)| \Delta^{\epsilon(\sigma)}  \&[+5em] |(ap)| \dA(\sigma) \&[+3em] |(axy)| \dA(ux,uy) \\
      };
      \path[maps,->] 
        (dp) edge node[above] {$ \bigl( \Lambda_uP_1, \dots, \Lambda_uP_r \bigr) $} (ap)
        (dn) edge node[above] {$\what{\sigma} $} (dp)
        (ap) edge node[above] {composition} (axy)
      ;
    \end{tikzpicture}
  \end{equation*}
  and this composition is uniquely determined by $\Lambda_u$. Altogether, this
  proves that $u \mapsto \Lambda_u$ is injective.

  In order to prove surjectivity of the assignment $u \mapsto \Lambda_u$, we
  have to construct for each labelling $\Lambda$ of $G$ in $\dA$  a simplicial
  functor $u\colon \dC[\Sigma] \to \dA$ such that
  $\Lambda = \Lambda_u$. The argument so far 
  already hints at the correct definition of $u$. In fact, we have no choice
  but to define $u(x) = \Lambda x$ on objects and to define the image
  $u(\sigma)$ of some $n$-simplex $\sigma \in \dC[\Sigma](x,y)$ as the simplex
  classified by the map 
  \begin{equation*}
    \begin{tikzpicture}[diagram]
      \matrix[objects] {
        |(dn)| \Delta^n \&[-1em] |(dp)| \Delta^{\epsilon(\sigma)}  \&[+5em] |(ap)| \dA(\sigma) \&[+3em] |(axy)| \dA(ux,uy)\smash{.} \\
      };
      \path[maps,->] 
        (dp) edge node[above] {$ \bigl( \Lambda P_1, \dots, \Lambda P_r \bigr) $} (ap)
        (dn) edge node[above] {$\what{\sigma} $} (dp)
        (ap) edge node[above] {composition} (axy)
      ;
    \end{tikzpicture}
  \end{equation*}
  Observe that if these assignments indeed define a functor $u$, then $\Lambda
  = \Lambda_u$.  Let us verify that these assignments define a simplicial
  functor. For two composable $n$-simplices $\sigma
  \in \dC[\Sigma](x,y)$ and $\tau \in \dC[\Sigma](y,z)$ we have a commutative diagram
  \begin{equation*}
    \begin{tikzpicture}[diagram]
      \matrix[objects] {
        |(dts)| \Delta^{\epsilon(\tau \comp \sigma)} \&[+8em] 
        |(ats)| \dA(\tau \comp \sigma) \\[-6ex]
        |(dprod)| \Delta^{\epsilon(\sigma)} \times \Delta^{\epsilon(\tau)}  \&[+8em]
        |(asat)| \dA(\sigma) \times \dA(\tau) \\[-1.5em]
         \& |(a2)| \dA(x,y) \times \dA(y,z) \\[-1.5em]
        |(dn)| \Delta^n \&   |(a)| \dA(ux,uz)\\
      };
      \path[maps,->] 
        (dprod) edge node[above] {$(\Lambda \sigma_1, \dots, \Lambda \sigma_a, \Lambda \tau_1, \dots, \Lambda \tau_b)$} (asat)
        (asat) edge node[right] {composition $\times$ composition} (a2)
        (dn) edge node[left] {$(\what{\sigma}, \what{\tau}) = \widehat{\tau \comp \sigma}$} (dprod)
        (a2) edge node[right] {composition} (a)
        (dn) edge node[below] {$u( \tau \comp \sigma)$} (a)
        (dn) edge node[above left] {$(u(\sigma), u(\tau))$} (a2)
        (dts) edge node[above] {$(\Lambda (\tau\comp \sigma)_1, \dots, \Lambda (\tau\comp\sigma)_{a+b})$} (ats)
      ;
      \node[anchor=center,rotate=90] at ($ (dprod) ! 0.5 ! (dts) $) {$=$};
      \node[anchor=center,rotate=90] at ($ (asat) ! 0.5 ! (ats) $) {$=$};
    \end{tikzpicture}
  \end{equation*}
  by \autoref{cor:godement-composition}. We thus conclude $u(\tau \comp \sigma) = u(\tau) \comp
  u(\sigma)$ by associativity of composition in $\dA$. Moreover, for any simplicial operator $\alpha$ we have a diagram
  \begin{equation*}
     \begin{tikzpicture}[diagram]
      \matrix[objects] {
        |(dm)| \Delta^m \&[-1em] |(dsa)| \Delta^{\epsilon(\sigma\alpha)}  \&[+5em] |(asa)| \dA(\sigma\alpha)  \& \\
        |(dn)| \Delta^n \& |(ds)| \Delta^{\epsilon(\sigma)}  \& |(as)| \dA(\sigma) \& |(axy)| \dA(x,y) \\
      };
      \path[maps,->] 
        (ds) edge node[below] {$ \bigl( \Lambda \sigma_1, \dots, \Lambda \sigma_a \bigr) $} (as)
        (dn) edge node[below] {$\what{\sigma} $} (ds)
        (dsa) edge node[above] {$ \bigl( \Lambda (\sigma\alpha)_1, \dots, \Lambda ( \sigma\alpha)_b \bigr) $} (asa)
        (dm) edge node[above] {$\widehat{\sigma\alpha} $} (dsa)
        (dm) edge node[left] {$\alpha$} (dn)
        (dsa) edge node[right] {$\epsilon(\alpha)$} (ds)
        (asa) edge node[left] {composition} (as)
        (as) edge node[below] {composition} (axy)
        (asa) edge node[above right] {composition} (axy)
      ;
    \end{tikzpicture}
  \end{equation*}
  in which the triangle commutes by associativity of composition in $\dA$ and
  in which the two squares on the left hand side commute by
  \autoref{rem:epsilon-lambda} and \autoref{prop:godement-representation}.  
\end{proof}

\chapter{Global Lifting Properties of Pasting Diagrams}\label{chap:global}

Throughout this chapter we fix a class $\sR$ of maps in $\sset$ that contains all isomorphisms and let $\sL$ denote the class of maps that have the left lifting property against all maps in $\sR$.
We call
the elements of $\sL$ and $\sR$ $\sL$-maps and $\sR$-maps, respectively.  A
simplicial functor $u \colon \dA \to \dB$ is a local $\sR$-functor if all the
maps $u\colon \dA(a,a') \to \dB(ua,ua')$, $a \in \dA$, are $\sR$-maps.  The main
result of this chapter is the following characterisation of those inclusions
between pasting diagrams that induce simplicial functors that have the left
lifting property against all local $\sR$-functors:
\globalthm*

The proof of \autoref{thm:global} relies on three independent steps. As a first
step, in \autoref{sec:sufficiency}, we prove that the hypothesis on the
inclusion $\Sigma \to \Pi$ in \autoref{thm:global} is sufficient: 
\begin{restatable*}{proposition}{liftingsufficient}\label{prop:lifting-sufficient}
  Let $\Sigma \to \Pi$ be an inclusion of complete pasting diagrams such that\/
  $\nerve(\Sigma_{x,y} \hc \Pi_{x,y}) \to \nerve(\Pi_{x,y})$ is an $\sL$-map
  for all vertices $x,y$ of\/ $\Sigma$. The functor\/ $\dC[\Sigma] \to \dC[\Pi]$
  then has the left lifting property against all local $\sR$-functors.
\end{restatable*}

The other two steps are necessary to show that the hypothesis on $\Sigma \to
\Pi$ in \autoref{thm:global} is necessary. To this end, we verify in
\autoref{sec:necessity1} that the hypothesis is necessary at least for the map
$\nerve(\Sigma \hc \Pi) \to \nerve(\Pi)$:
\begin{restatable*}{proposition}{liftingnecessary}\label{prop:lifting-necessary}
  Let $\Sigma \to \Pi$ be an inclusion of complete pasting diagrams. If\/
  $\dC[\Sigma] \to \dC[\Pi]$ has the left lifting property against all local
  $\sR$-functors, then $\nerve(\Sigma \hc \Pi) \to \nerve(\Pi)$ is an
  $\sL$-map.
\end{restatable*}

The final step in the proof of \autoref{thm:global} is then to descend from the
maps $\nerve(\Sigma \hc \Pi) \to \nerve(\Pi)$ handled by
\autoref{prop:lifting-necessary} to all the maps $\nerve(\Sigma_{x,y} \hc
\Pi_{x,y}) \to \nerve(\Pi_{x,y})$ with $x,y \in G$. This is achieved by the
following proposition whose proof we give in \autoref{sec:necessity2}.
\begin{restatable*}{proposition}{liftingrestriction}\label{prop:lifting-restriction}
  If\/ $\Sigma \to \Pi$ is an inclusion of complete pasting diagrams such that\/
  $\dC[\Sigma] \to \dC[\Pi]$ has the left lifting property against all local
  $\sR$-functors, then so do the functors $\dC[\Sigma_{x,y}] \to
  \dC[\Pi_{x,y}]$ for all vertices $x,y$ of\/ $\Sigma$.
\end{restatable*}

In the final section \autoref{sec:proof-b} we then give the proof of
\autoref{thm:global} that we sketched in this introduction.

\section{Sufficiency of the hypothesis}\label{sec:sufficiency}
We keep the classes $\sL$ and $\sR$ of maps in $\sset$ with $\sL =  {}^\lifts
\sR$ fixed.  This section is devoted to a proof of the following proposition
that we already stated above. 
\liftingsufficient
In order to prove \autoref{prop:lifting-sufficient}, we thus have to show that for each square
\begin{equation}\label{dia:lifting-problem}
  \begin{tikzpicture}[diagram]
    \matrix[objects] {
      |(s)| \dC[\Sigma] \& |(b)| \dB \\
      |(t)| \dC[\Pi] \& |(a)| \dA \\
    };
    \path[maps,->] 
      (s) edge node[above] {$u$} (b)
      (s) edge (t)
      (b) edge node[right] {$p$} (a)
      (t) edge node[below] {$v$} (a)
    ;
  \end{tikzpicture}
\end{equation}
in which $p$ is a local $\sR$-functor, there exists a simplicial functor $\ell \colon \dC[\Pi] \to \dB$ such that
\begin{equation}\label{dia:lifting-solution}
  \begin{tikzpicture}[diagram]
    \matrix[objects] {
      |(s)| \dC[\Sigma] \& |(b)| \dB \\
      |(t)| \dC[\Pi] \& |(a)| \dA \\
    };
    \path[maps,->] 
      (s) edge node[above] {$u$} (b)
      (s) edge (t)
      (b) edge node[right] {$p$} (a)
      (t) edge node[below] {$v$} (a)
      (t) edge node {$\ell$} (b)
    ;
  \end{tikzpicture}
\end{equation}
commutes. Throughout this section, we keep an inclusion $\Sigma \to \Pi$ of
pasting diagrams satisfying the hypothesis of \autoref{prop:lifting-sufficient}
and a square such as \eqref{dia:lifting-problem} fixed and assume $p$ to be a
local $\sR$-functor. We let $G$ denote the graph underlying both $\Sigma$ and
$\Pi$. Observe that $\dC[\Sigma] \to \dC[\Pi]$ is the identity on objects, so
that $pu(x) = v(x)$ for any vertex $x \in G$. We thus define $\ell(x) = u(x)$
on objects and this definition lets \eqref{dia:lifting-solution} commute on the
level of objects.  It remains to construct for each two vertices $x,z \in G$
suitably functorial lifts $\ell_{x,z}$ as in the diagram 
  \begin{equation}\label{dia:lifting-local}
    \begin{tikzpicture}[diagram]
      \matrix[objects] {
        |(nsxy)| \nerve(\Sigma_{x,z}) \&[-2em] |(s)| \dC[\Sigma](x,z) \& |(b)| \dB(ux,uz) \\
        |(npxy)| \nerve(\Pi_{x,z}) \& |(t)| \dC[\Pi](x,z) \& |(a)| \dA(vx,vz)\smash{.} \\
      };
      \path[maps,->] 
        (s) edge node[above] {$u$} (b)
        (s) edge (t)
        (b) edge node[right] {$p$} (a)
        (t) edge node[below] {$v$} (a)
        (t) edge node[above left] {$\ell_{x,z}$} (b)
        (nsxy) edge node[left] {inclusion} (npxy)
      ;
      \node at ($ (nsxy.east) ! 0.5 ! (s.west) $) {$=$};
      \node at ($ (npxy.east) ! 0.5 ! (t.west) $) {$=$};
    \end{tikzpicture}
  \end{equation}
  Here,
  functoriality means nothing but that for all vertices $x$ of $G$ the map $\ell_{x,x}
  \colon \nerve(\Pi_{x,x}) = \Delta^0 \to \dB(ux,ux)$ classifies the identity
  $\id_{ux} \in \dB(ux,ux)_0$, and that the square
  \begin{equation}\label{dia:functoriality}
    \begin{tikzpicture}[diagram]
      \matrix[objects] {
        |(xyz)| \nerve(\Pi_{x,y}) \times \nerve(\Pi_{y,z}) \&[+2em]
        |(bxyz)| \dB(ux,uy) \times \dB(uy,uz) \\
        |(xz)| \nerve(\Pi_{x,z}) \& 
        |(bxz)| \dB(ux,uz) \\
      };
      \path[maps,->] 
        (xyz) edge node[above] {$\ell_{x,y} \times \ell_{y,z}$} (bxyz)
        (xyz) edge node[left] {composition} (xz)
        (xz) edge node[below] {$\ell_{x,z}$} (bxz)
        (bxyz) edge node[right] {composition} (bxz)
      ;
    \end{tikzpicture}
  \end{equation}
  commutes for all vertices $x,y,z \in G$. The constraint on $\ell_{x,x}$ is
  easily satisfied by simply defining $\ell_{x,x}$ to be the map classifying
  the identity in $\dB(ux,ux)_0$. It is then clear that
  \eqref{dia:functoriality} commutes for all $x,y,z \in G$ with $x = y$ or $y =
  z$. Moreover, the following lemma identifies a further case, where
  \eqref{dia:functoriality} commutes for trivial reasons.
  \begin{lemma}\label{lem:trivial-functoriality}
    The square \eqref{dia:functoriality} commutes whenever $y \notin G_{x,z}$.
  \end{lemma}
  \begin{proof}
    Suppose $y \notin G_{x,z}$. There is then no directed path from $x$ to $y$
    or no directed path from $y$ to $z$ in $G$. This implies in particular that
    $\nerve(\Pi_{x,y}) = \emptyset$ or $\nerve(\Pi_{y,z}) = \emptyset$ and thus
    $\nerve(\Pi_{x,y}) \times \nerve(\Pi_{y,z}) = \emptyset$. Hence,
    \eqref{dia:functoriality} commutes simply because its upper left corner is
    an initial object. 
  \end{proof}
  
  Let us sum up our discussion so far: In order to solve the lifting problem
  \eqref{dia:lifting-problem} it suffices to devise maps $\ell_{x,z} \colon
  \nerve(\Pi_{x,z}) \to \dB(ux,uz)$ as in \eqref{dia:lifting-local} for all
  pairs $x,z$ of distinct vertices of $G$ such that diagram
  \eqref{dia:functoriality} commutes for all vertices $x,z \in G$ and $y \in
  G_{x,z} \smallsetminus \{x,z\}$. We construct these maps $\ell_{x,z}$ by
  recursion over the partial order \enquote{$\preceq$} on the set of pairs
  $(x,z)$ of vertices of $G$ given by
  \begin{equation*}
    (x,z) \preceq (x',z') \quad\text{if and only if}\quad G_{x,z} \subseteq G_{x',z'}.
  \end{equation*}

  To this end, let us fix two vertices $x,z \in G$ and assume that all
  $\ell_{a,b}$ with $(a,b) \prec (x,z)$ have already been constructed. Further
  assume that these $\ell_{a,b}$ are functorial in the sense that
  \eqref{dia:functoriality} with $a,b,c$ instead of $x,y,z$ commutes for all
  $a,c \in G$ with $(a,c) \prec (x,z)$ and $b \in G_{a,c}$.  Given this data,
  we have to furnish a map $\ell_{x,z} \colon \nerve(\Pi_{x,z}) \to \dB(ux,uz)$
  such that \eqref{dia:lifting-local} commutes and such that
  \eqref{dia:functoriality} commutes for all $y \in G_{x,z}$.

  Observe that for any vertex $y \in G_{x,z}\smallsetminus\{x,z\}$ we have
  $(x,y), (y,z) \prec (x,z)$. This means in particular that we have for any
  such vertex $y$ the maps $\ell_{x,y}$ and $\ell_{y,z}$ at our disposal. We
  thus define for any such vertex $y$ the map 
  \begin{equation*}
   h_y\colon \nerve(\Pi_{x,y} \join \Pi_{y,z}) \to \dB(ux,uz) 
  \end{equation*}
  as the composition
  \begin{equation*}
    \begin{tikzpicture}[diagram]
      \matrix[objects] {
        |(join)| \nerve(\Pi_{x,y} \join \Pi_{y,z}) \&[-2.3em] |(prod)| \nerve(\Pi_{x,y}) \times
      \nerve(\Pi_{y,z}) \&[+1.2em] |(bprod)| \dB(ux,uy) \times
    \dB(uy,uz) \\
      \& \& |(b)| \dB(ux,uz)\smash{.} \\
      };
      \path[maps,->] 
        (prod) edge node[above] {$\ell_{x,y} \times \ell_{y,z}$} (bprod)
        (bprod) edge node[right] {composition} (b)
        (join) edge[] node[below left] {$h_y$} (b)
      ;
      \node at ($ (join.east) ! 0.5 ! (prod.west) $) {$\iso$};
    \end{tikzpicture}
  \end{equation*}
  We know from \autoref{lem:nerve-of-hc} that
  \begin{equation*}
    \nerve(\Sigma_{x,y} \hc \Pi_{x,y}) 
    = \nerve(\Sigma_{x,z})\ \cup\ \bigcup_{\mathclap{y \in G_{x,z} \smallsetminus \{x,z\}}}\ \nerve(\Pi_{x,y} \join \Pi_{y,z}),
  \end{equation*}
  where the unions are taken in $\nerve(\Pi_{x,z})$.  Observe that the domain
  of $h_y$ is a simplicial subset of $\nerve(\Sigma_{x,z} \hc \Pi_{x,z})$.  The
  following two lemmata therefore guarantee that these maps $h_y$ together with
  $u \colon \dC[\Sigma](x,z) \to \dB(ux,uz)$ induce a map $h\colon
  \nerve(\Sigma_{x,z} \hc \Pi_{x,z} ) \to \dB(ux,uz)$ and this permits us to
  eventually utilise the hypothesis on the inclusion $\Sigma \to \Pi$ in the
  statement of \autoref{prop:lifting-sufficient}.
  \begin{restatable}{lemma}{uhylemma}\label{lem:uhy}
    Let $y \in G_{x,z} \smallsetminus \{x,z\}$. The map $h_y$ and 
    \begin{equation*}
      u \colon \nerve(\Sigma_{x,z}) = \dC[\Sigma](x,z) \to \dB(ux,uz)
    \end{equation*}
    coincide on the intersection of their domains.
  \end{restatable}
  \begin{restatable}{lemma}{hyhylemma}\label{lem:hyhyp}
    Let $y,y' \in G_{x,z} \smallsetminus \{x,z\}$. The maps $h_y$ and $h_{y'}$
    coincide on the intersection of their domains.
  \end{restatable}
  We defer the proof of these lemmata until the end of this section so as not
  to distract from the main argument.  In order to understand the composition
  $p \comp h$, where $h \colon \nerve(\Sigma_{x,z} \hc \Pi_{x,z}) \to \dB(ux,uz)
  $ is the map whose existence is guaranteed by the above two
  lemmata, it suffices to understand the composition $p \comp u$ and all the
  compositions $p \comp h_y$, where $y \in G_{x,z} \smallsetminus \{x,z\}$. The
  former composition is known from \eqref{dia:lifting-problem} while the latter
  composition can easily be computed as in the following remark:
  \begin{remark}\label{rem:vhy}
    The diagram
    \begin{equation*}
      \begin{tikzpicture}[diagram]
        \matrix[objects] {
          |(prod)| \nerve(\Pi_{x,y} ) \times \nerve(\Pi_{y,z}) \&[+1em]
          |(bprod)| \dB(ux,uy) \times \dB(uy,uz) \&
          |(b)| \dB(ux,uz)  \\
          \& |(aprod)| \dA(ux,uy) \times \dA(uy,uz) 
          \& |(a)| \dA(ux,uz) \\
        };
        \path[maps,->] 
          (bprod) edge node[above] {comp.} (b)
          (aprod) edge node[below] {comp.} (a)
          (bprod) edge node[left] {$p \times p$} (aprod)
          (b) edge node[right] {$p$} (a)
          (prod) edge node[above] {$\ell_{x,y} \times \ell_{y,z}$} (bprod)
          (prod) edge node[below left] {$v \times v$} (aprod)
        ;
      \end{tikzpicture}
    \end{equation*}
    commutes by functoriality of $p$ and \eqref{dia:lifting-local} for
    $\ell_{x,y}$ and $\ell_{y,z}$. The composition $p \comp h_y$ therefore appears in the commutative diagram 
    \begin{equation*}
      \begin{tikzpicture}[diagram]
        \matrix[objects] {
          |(prod)| \nerve(\Pi_{x,y} ) \times \nerve(\Pi_{y,z}) \&[-3em] |(join)| \nerve(\Pi_{x,y} \join \Pi_{y,z}) \& |(n)| \nerve(\Pi_{x,z}) \\ 
          |(aprod)| \dA(ux,uy) \times \dA(uy,uz) \& \& |(a)| \dA(ux,uz) \\
        };
        \path[maps,->] 
          (prod) edge[bend left] node[above] {composition in $\dC[\Pi]$} (n)
          (n) edge node[right] {$v$} (a)
          (join) edge node[above] {inclusion} (n)
          (prod) edge node[left] {$v \times v$} (aprod)
          (aprod) edge node[below] {composition} (a)
          (join) edge node[below left] {$p \comp h_y$} (a)
        ;
        \node at ($ (prod.east) !0.5 !(join.west) $) {$\iso$};
      \end{tikzpicture}
    \end{equation*}
    and we conclude that $p \comp h_y$ is nothing but the composition of $v$
    and the canonical inclusion $\nerve(\Pi_{x,y} \join \Pi_{y,z}) \to
    \nerve(\Pi_{x,z})$.
  \end{remark}
  
  Altogether, we see that the diagram 
  \begin{equation}\label{dia:recursive-lifting}
    \begin{tikzpicture}[diagram]
      \matrix[objects] {
        |(s)| \nerve(\Sigma_{x,z}) \&
        |(st)| \nerve(\Sigma_{x,z} \hc \Pi_{x,z}) \& |(b)| \dB(ux,uz) \\
        \& |(t)| \nerve(\Pi_{x,z}) \& |(a)| \dA(vx,vy) \\
      };
      \path[maps,->] 
        (s) edge (st)
        (s) edge [bend left] node[above] {$u$} (b)
        (s) edge (t)
        (st) edge node[above] {$h$} (b)
        (st) edge (t)
        (t) edge node[below] {$v$} (a)
        (b) edge node[right] {$p$} (a)
      ;
    \end{tikzpicture}
  \end{equation}
  commutes. Moreover, $p$ is an $\sR$-map by assumption
  and the left vertical map $\nerve(\Sigma_{x,z} \hc \Pi_{x,z}) \to \nerve(\Pi_{x,z})$ is an
  $\sL$-map by our hypothesis on the inclusion $\Sigma \to \Pi$ in the
  statement of \autoref{prop:lifting-sufficient}. We therefore find a lift
  $\ell_{x,z}\colon
  \nerve(\Pi_{x,z}) \to \dB(ux,uz)$ in \eqref{dia:recursive-lifting} and this lift renders \eqref{dia:lifting-local} commutative, too. 
  It only remains to check that the lift $\ell_{x,z}$ in \eqref{dia:recursive-lifting} is functorial
  in the sense that all the diagrams \eqref{dia:functoriality} with $y \in
  G_{x,z} \smallsetminus \{x,z\}$ commute. 
  To this end, let us consider the diagram 
  \begin{equation*}
    \begin{tikzpicture}[diagram]
      \matrix[objects] {
        |(pxyz)| \nerve(\Pi_{x,y}) \times \nerve(\Pi_{y,z})  \&[-3em]
        |(join)| \nerve(\Pi_{x,y} \join \Pi_{y,z}) 
      \&[-2.1em] |(hc)| \nerve(\Sigma_{x,z} \hc \Pi_{x,z}) \&[-2.6em] |(pxz)| \nerve(\Pi_{x,z}) \\
        |(bxyz)| \dB(ux,uy) \times \dB(uy,uz) \& \& \&|(bxz)| \dB(ux,uz) \smash{.}\\
      };
      \path[maps,->] 
        (pxyz) edge node[left] {$\ell_{x,y} \times \ell_{y,z}$} (bxyz)
        (bxyz) edge node[below] {composition} (bxz)
        (hc) edge node[below left] {$h$} (bxz)
        (pxz) edge node[right] {$\ell_{x,z}$} (bxz)
        (hc) edge (pxz)
        (join) edge (hc)
      ;
      \node at ($ (pxyz.east) ! 0.5 ! (join.west) $) {$\iso$};
    \end{tikzpicture}
  \end{equation*}
  The triangle commutes by our construction of $\ell_{x,z}$ as a lift in
  \eqref{dia:recursive-lifting} and the quadriliteral commutes because both
  compositions from $\nerve(\Pi_{x,y} \join \Pi_{y,z})$ to $\dB(ux,uz)$ are
  nothing but the map $h_y$. Functoriality of $\ell_{x,z}$ therefore follows
  from the simple observation that the top row in this diagram is 
  composition in $\dC[\Pi]$. This finishes the proof of
  \autoref{prop:lifting-sufficient} except for the fact that we owe the reader
  proofs for \autoref{lem:uhy} and \ref{lem:hyhyp}. 
    \uhylemma*
    \begin{proof}[Proof of \autoref{lem:uhy}]
      The domain of $h_y$ is $\nerve(\Pi_{x,y} \join \Pi_{y,z})$ and the domain
      of $u$ is $\nerve(\Sigma_{x,z})$. Their intersection in
      $\nerve(\Pi_{x,z})$ is nothing but $\nerve(\Sigma_{x,y} \join
      \Sigma_{y,z})$ by \autoref{cor:nerve-pullback}.
      The restrictions of $u$ and $h_y$ to this intersection feature as the bottom-left and top-right compositions in the diagram 
      \begin{equation*}
        \begin{tikzpicture}[diagram]
          \matrix[objects] {
            |(joins)| \nerve( \Sigma_{x,y} \join \Sigma_{y,z})
            \& |(joinp)| \nerve( \Pi_{x,y} \join \Pi_{y,z} )
            \&[-2em] |(prodp)| \nerve( \Pi_{x,y} ) \times \nerve(\Pi_{y,z}) \\[-1.75ex]
            |(prods)| \nerve(\Sigma_{x,y}) \times \nerve(\Sigma_{y,z})\& \& |(prodb)| \dB(ux,uy) \times \dB(uy,uz) \\[-1.75ex]
            |(s)| \nerve(\Sigma_{x,z}) \& \& |(b)| \dB(ux,uz)\smash{.} \\
          };
          \path[maps,->] 
            (joins) edge (joinp)
            (joins) edge node[left=.6ex,anchor=center,rotate=90] {$\sim$} (prods)
            (prodb) edge node[right] {composition} (b)
            (prodp) edge node[right] {$\ell_{x,y} \times \ell_{y,z}$} (prodb)
            (prods) edge node[below] {$u \times u$} (prodb)
            (prods) edge node[left] {composition} (s)
            (s) edge node[below] {$u$} (b)
          ;
          \node at ($ (joinp.east) ! 0.5 ! (prodp.west) $) {$\iso$};
        \end{tikzpicture}
      \end{equation*}
      This is obvious for $h_y$ and follows for $u$ from the definition of
      composition in $\dC[\Sigma]$. The bottom square in this diagram commutes
      by functoriality of $u$ and top square commutes by naturality of the
      isomorphisms $\nerve(\Sigma_{x,y} \join \Sigma_{y,z}) \iso
      \nerve(\Sigma_{x,y}) \times \nerve(\Sigma_{y,z})$ together with the
      standing assumption that $\ell_{x,y}$ and $\ell_{y,z}$ render the diagram
      \eqref{dia:lifting-local} commutative. In fact, the upper left triangle
      in \eqref{dia:lifting-local} suffices for the diagram at hand.
  \end{proof}

  \hyhylemma*
  \begin{proof}[Proof of \autoref{lem:hyhyp}]
    \autoref{cor:nerve-join-intersection} tells us that $\dom(h_y)$ and $\dom(h_{y'})$ intersect in $\nerve(\Pi_{x,z})$ if and only if 
    $G$ contains a directed path between $y$ and $y'$. As the
    situation is symmetric, we may suppose that $G$ contains a directed path
    from $y$ to $y'$. According to \autoref{cor:nerve-join-intersection}, the
    intersection of $\dom(h_y)$ and $\dom(h_{y'})$ is then given by
      $\nerve(\Pi_{x,y} \join
      \Pi_{y,y'} \join \Pi_{y',z}) $.
    By naturality of the isomorphisms $\nerve(\Pi_{x,y} \join \Pi_{y,z})
    \iso \nerve(\Pi_{x,y}) \times \nerve(\Pi_{y,z})$ the statement of the lemma
    thus reduces to the commutativity of the diagram
  \begin{equation*}
    \begin{tikzpicture}[diagram]
      \matrix[objects] {
        \&[-5.5em]|(prod)| \nerve( \Pi_{x,y} ) \times \nerve( \Pi_{y,y'} ) \times \nerve( \Pi_{y',z} ) \&[-5.5em] \\
         |(zp)| \nerve(\Pi_{x,y'}) \times \nerve(\Pi_{y',z}) \& \& |(z)| \nerve(\Pi_{x,y}) \times \nerve(\Pi_{y,z}) \\[-4ex]
         \& |(bzzp)| \dB(ux,uy) \times \dB(uy,uy') \times \dB(uy',uz) \\
         |(bzp)| \dB(ux,uy') \times \dB(uy',uz) \& \& |(bz)| \dB(ux,uy) \times \dB(uy,uz) \\
         \& |(b)| \dB(ux,uz)\smash{,}\\
      };
      \path[maps,->] 
        (bzp) edge node[below left] {composition} (b)
        (bz) edge node[below right] {composition} (b)
        (prod) edge node[above left] {composition $\times \nerve(\Pi_{y',z})$} (zp)
        (prod) edge node[above right] {$\nerve(\Pi_{x,y}) \times$ composition} (z)
        (prod) edge node {$\ell_{x,y} \times \ell_{y,y'} \times \ell_{y',z}$} (bzzp)
        (zp) edge node[left] {$\ell_{x,y'} \times \ell_{y',z}$} (bzp)
        (z) edge node[right] {$\ell_{x,y} \times \ell_{y,z}$} (bz)
        (bzzp) edge node[] {composition $\times \dB(uy',uz)$} (bzp)
        (bzzp) edge node[] {$\dB(ux,uy) \times $ composition} (bz)
      ;
    \end{tikzpicture}
  \end{equation*}
    as the compositions on the left and right hand side of this diagram are --
    up to coherent isomorphism -- nothing but the restrictions of $h_{y'}$ and
    $h_{y}$ to the intersection of their domains. However, the upper two
    squares in this diagram commute by functoriality of our lifts
    $\ell_{\ph,\ph}$ and the lower square commutes by associativity of
    composition in $\dB$. This concludes the proof.
  \end{proof}

  \section{Necessity of the hypothesis, part 1}\label{sec:necessity1}
  In this section, we give a proof of the following proposition:
  \liftingnecessary*

  The proof of \autoref{prop:lifting-necessary} works by exhibiting for each
  commutative square
  \begin{equation}\label{dia:hc-lifting-problem}
    \begin{tikzpicture}[diagram]
      \matrix[objects] {
        |(sp)| \nerve(\Sigma \hc \Pi) \& |(x)| X \\
        |(p)| \nerve(\Pi) \& |(y)| Y\smash{,} \\
      };
      \path[maps,->] 
        (sp) edge node[above] {$u$} (x)
        (sp) edge (p)
        (p) edge node[below] {$v$} (y)
        (x) edge node[right] {$p$} (y)
      ;
    \end{tikzpicture}
  \end{equation}
  in which $p$ is an $\sR$-map, 
  another commutative square 
  \begin{equation}\label{dia:dc-lifting-problem}
      \begin{tikzpicture}[diagram]
        \matrix[objects] {
          |(s)| \dC[\Sigma] \& |(u)| \dC[\Sigma \hc \Pi]_{/u} \\
          |(p)| \dC[\Pi] \& |(pu)| \dC[\Sigma \hc \Pi]_{/pu}\smash{,} \\
        };
        \path[maps,->] 
          (s) edge (u)
          (s) edge (p)
          (p) edge node[below] {$v/pu$} (pu)
          (u) edge node[right] {$\dC_{/p}$} (pu)
        ;
      \end{tikzpicture}
  \end{equation}
    in which $\dC_{/p}$ is a local $\sR$-functor, that has the property that
    any solution to the lifting problem \eqref{dia:dc-lifting-problem} gives a
    solution to the original lifting problem \eqref{dia:hc-lifting-problem}. 

    Even though the definitions and remarks of this section are only relevant
    to the proof of \autoref{prop:lifting-necessary}, we still prefer to give
    them in their general form.
  \begin{definition}\label{def:over-simplicial-set}
    Let $\Sigma$ be a complete pasting diagram with source $s$ and target~$t$.
    Further let $u\colon \nerve(\Sigma) \to X$ be a map of simplicial sets. We
    define a simplicial category $\dC[\Sigma]_{/u}$ as follows:
      The objects of\/ $\dC[\Sigma]_{/u}$ are the vertices of $\Sigma$ and the 
      mapping spaces of\/ $\dC[\Sigma]_{/u}$ are given by
      \begin{equation*}
        \dC[\Sigma]_{u}(x,y) = \begin{cases}
          \dC[\Sigma](x,y) &\text{if } (x,y) \neq (s,t), \\
          X & \text{if } (x,y) = (s,t).
        \end{cases}
      \end{equation*}
      The composition laws in $\dC[\Sigma]_{/u}$ are those of\/ $\dC[\Sigma]$
      except for the compositions
      \begin{equation*}
        \dC[\Sigma]_{/u}(s,x) \times \dC[\Sigma]_{/u}(x,t) \to \dC[\Sigma]_{/u}(s,t),
      \end{equation*}
      which are given by 
      \begin{equation*}
        \begin{tikzpicture}[diagram]
          \matrix[objects] {
            |(csxt)| \dC[\Sigma]_{/u}(s,x) \times \dC[\Sigma]_{/u}(x,t) \&[+3em] \& |(cst)| \dC[\Sigma]_{/u}(s,t) \\[-4ex]
            |(nsxt)| \nerve(\Sigma_{s,x}) \times \nerve(\Sigma_{x,t}) \& |(nst)| \nerve(\Sigma_{s,t}) \& |(x)| X\smash{.} \\
          };
          \path[maps,->] 
            (csxt) edge (cst)
            (nsxt) edge node[below] {composition in $\dC[\Sigma]$} (nst)
            (nst) edge node[below] {$u$} (x)
          ;
          \node[anchor=center,rotate=90]  at ($ (csxt) ! 0.5 ! (nsxt) $) {$=$};
          \node[anchor=center,rotate=90]  at ($ (cst) ! 0.5 ! (x) $) {$=$};
        \end{tikzpicture}
      \end{equation*}
      if\/ $x \notin \{s,t\}$ and by the isomorphisms $\Delta^0 \times X \iso X$ and $X \times \Delta^0 \iso X$ if\/ $x = s$ or $x = t$, respectively.
  \end{definition}
  \begin{remark}
    We should justify that $\dC[\Sigma]_{/u}$ is well-defined, i.\,e.\ that the
    composition laws for $\dC[\Sigma]_{/u}$ are associative and unital. This,
    however, is immediate from the definition and left to the reader.
  \end{remark}

  \begin{remark}
  Observe that we have a canonical functor $u\colon \dC[\Sigma] \to
    \dC[\Sigma]_{/u}$ that is the identity on objects and on all mapping spaces
    $\dC[\Sigma](x,y)$ with $(x,y) \neq (s,t)$, and acts on $\dC[\Sigma](s,t) =
    \nerve(\Sigma)$ by $u$.
  \end{remark}

  \begin{remark}\label{rem:over-u-functorial}
    If $p \colon X \to Y$ is another map of simplicial sets, then there is a
    canonical functor $\dC_{/p} \colon \dC[\Sigma]_{/u} \to \dC[\Sigma]_{/pu}$
    that is the identity on objects and on all mapping spaces
    $\dC[\Sigma]_{/u}(x,y)$ with $(x,y) \neq (s,t)$, and acts on
    $\dC[\Sigma]_{/u}(s,t) = X$ by $p$. The reader should note that $\dC_{/p}$ is
    a local $\sR$-functor whenever $p \colon X \to Y$ is an $\sR$-map.
  \end{remark}

  \begin{definition}\label{def:voverw}
    Let $v \colon \nerve(\Pi) \to Y$ and $w \colon \nerve(\Sigma \hc \Pi) \to
    Y$ be two maps of simplicial sets.  We define a functor $v/w \colon
    \dC[\Pi] \to \dC[\Sigma \hc \Pi]_{/w}$ as follows: The functor $v/w$
    is the identity on objects. The action of\/ $v/w$ on mapping spaces
    $\dC[\Pi](x,y)$ with $(x,y) \neq (s,t)$ is given by the identity on 
    \begin{equation*}
      \dC[\Pi](x,y) = \nerve(\Pi_{x,y}) = \nerve\bigl( (\Sigma \hc \Pi)_{x,y} \bigr) = \dC[\Sigma \hc \Pi]_{/w}(x,y),
    \end{equation*}
    where the second equality follows from \autoref{lem:hc-is-full} and the action
    of\/ $v/w$ on the mapping space $\dC[\Pi](s,t)$ is given by the
    composition
    \begin{equation*}
      \dC[\Pi](s,t) = \nerve(\Pi) \xto{v} Y = \dC[\Sigma \hc \Pi]_{/w}(s,t).
    \end{equation*}
  \end{definition}

  We leave it to the reader to verify that the functor $v/w$ in
  \autoref{def:voverw} is well-defined. \autoref{def:voverw} was the last piece
  missing in order to be able to transform squares such as
  \eqref{dia:hc-lifting-problem} into squares \eqref{dia:dc-lifting-problem}.
  This is accomplished by the following lemma:
  \begin{lemma}\label{lemma:square-transform}
    Suppose that a square such as \eqref{dia:hc-lifting-problem} is given.
    The diagram
    \begin{equation*}
      \begin{tikzpicture}[diagram]
        \matrix[objects] {
          |(s)| \dC[\Sigma \hc \Pi] \& |(u)| \dC[\Sigma \hc \Pi]_{/u} \\
          |(p)| \dC[\Pi] \& |(pu)| \dC[\Sigma \hc \Pi]_{/pu} \\
        };
        \path[maps,->] 
          (s) edge node[above] {$u$} (u)
          (s) edge (p)
          (u) edge node[right] {$\dC_{/p}$} (pu)
          (p) edge node[below] {$v/pu$} (pu)
        ;
      \end{tikzpicture}
    \end{equation*}
    commutes. 
  \end{lemma}
  \begin{proof}
    Both compositions are the identity on objects and all mapping spaces
    $\dC[\Sigma \hc \Pi](x,y)$ with $(x,y) \neq (s,t)$. On the remaining
    mapping space $\dC[\Sigma \hc \Pi](s,t)$, the diagram given in the lemma is
    nothing but \eqref{dia:hc-lifting-problem}.
  \end{proof}

  \begin{lemma}\label{lemma:solution-transport}
    Suppose a square such as \eqref{dia:hc-lifting-problem} to be given. Suppose further that the diagram
    \begin{equation*}
      \begin{tikzpicture}[diagram]
        \matrix[objects] {
          |(s)| \dC[\Sigma] \& |(hc)| \dC[\Sigma \hc \Pi] \& |(u)| \dC[\Sigma \hc \Pi]_{/u} \\
          |(p)| \dC[\Pi] \& \& |(pu)| \dC[\Sigma \hc \Pi]_{/pu} \\
        };
        \path[maps,->] 
          (s) edge (hc)
          (hc) edge node[above] {$u$} (u)
          (s) edge (p)
          (u) edge node[right] {$\dC_{/p}$} (pu)
          (p) edge node[below] {$v/pu$} (pu)
          (p) edge node {$\ell$} (u)
        ;
      \end{tikzpicture}
    \end{equation*}
    commutes. The map $\ell_{s,t} \colon \nerve(\Pi) \to X$ given by the
    action of\/ $\ell$ on the mapping space $\dC[\Pi](s,t)$ then solves the
    lifting problem \eqref{dia:hc-lifting-problem}.
  \end{lemma}
  \begin{proof}
    It is immediate that $p \ell_{s,t} = v$ from our definitions of $\dC/p$ and $v/w$.  Let us verify that the
    composition
    \begin{equation*}
      \nerve(\Sigma \hc \Pi) \to \nerve(\Pi) \xto{\ell_{s,t}} X
    \end{equation*}
    is equal to $u$. By \autoref{lem:nerve-of-restriction}, we may write the
    domain of this composition as
    \begin{equation*}
      \nerve(\Sigma \hc \Pi) = \nerve(\Sigma) \cup\ \bigcup_{\mathclap{x \in G \smallsetminus\{s,t\}}}\ \nerve(\Pi_{s,x} \join \Pi_{x,t})
    \end{equation*}
    and it thus suffices to show that the map
    \begin{equation}\label{eq:easy-composition}
      \nerve(\Sigma) \to \nerve(\Pi) \xto{\ell_{s,t}} X 
    \end{equation}
    and all the maps 
    \begin{equation}\label{eq:complicated-composition}
      \nerve(\Pi_{s,x} \join \Pi_{x,t} ) \to \nerve(\Pi) \xto{\ell_{s,t}} X
    \end{equation}
    with $x \in G \smallsetminus \{s,t\}$ coincide with the restrictions of $u$
    to their domain. For the former map \eqref{eq:easy-composition} this
    follows immediately from the commutative square given in the lemma.  For
    the latter maps, we observe that $\ell_{s,x}$ and $\ell_{x,t}$ are
    necessarily identity maps and that functoriality of $\ell$ hence implies
    that the square 
    \begin{equation*}
      \begin{tikzpicture}[diagram]
        \matrix[objects] {
          |(prod)| \nerve(\Pi_{s,x}) \times \nerve(\Pi_{x,t}) \&[-2.5em]
          |(join)| \nerve( \Pi_{s,x} \join \Pi_{x,t} ) \&[+1em] \&[-1em] |(np)| \nerve(\Pi) \\
          |(prod2)| \nerve(\Pi_{s,x}) \times \nerve(\Pi_{x,t}) \& |(join2)| \nerve( \Pi_{s,x} \join \Pi_{x,t} )\& |(hc)| \nerve(\Sigma \hc \Pi) \& |(n)| X\\
        };
        \path[maps,->] 
          (prod) edge node[left] {$\ell_{s,x} \times \ell_{x,t}$} (prod2)
          (join) edge node[right] {$1$} (join2)
          (join) edge node[above] {inclusion} (np)
          (join2) edge node[below] {inclusion} (hc)
          (hc) edge node[below] {$u$} (n)
          (np) edge node[right] {$\ell_{s,t}$} (n)
        ;
        \node at ($ (prod.east) ! 0.5 ! (join.west) $) {$\iso$};
        \node at ($ (prod2.east) ! 0.5 ! (join2.west) $) {$\iso$};
      \end{tikzpicture}
    \end{equation*}
    commutes, i.\,e.\ that the maps \eqref{eq:complicated-composition} coincide
    with the restriction of $u$ to their domain. This finishes the proof.
  \end{proof}

  We close this section with the proof of \autoref{prop:lifting-necessary} that
  we sketched at its beginning.
  \liftingnecessary
  \begin{proof}
    Suppose that 
    \begin{equation*}
    \begin{tikzpicture}[diagram]
      \matrix[objects] {
        |(sp)| \nerve(\Sigma \hc \Pi) \& |(x)| X \\
        |(p)| \nerve(\Pi) \& |(y)| Y\smash{,} \\
      };
      \path[maps,->] 
        (sp) edge node[above] {$u$} (x)
        (sp) edge (p)
        (p) edge node[below] {$v$} (y)
        (x) edge node[right] {$p$} (y)
      ;
    \end{tikzpicture}
  \end{equation*}
    is a commutative diagram of simplicial sets in which $p$ is an $\sR$-map.
    We then have a commutative diagram 
    \begin{equation*}
      \begin{tikzpicture}[diagram]
        \matrix[objects] {
          |(s)| \dC[\Sigma] \& |(u)| \dC[\Sigma \hc \Pi]_{/u} \\
          |(p)| \dC[\Pi] \& |(pu)| \dC[\Sigma \hc \Pi]_{/pu}\smash{,} \\
        };
        \path[maps,->] 
          (s) edge (u)
          (s) edge (p)
          (p) edge node[below] {$v/pu$} (pu)
          (u) edge node[right] {$\dC_{/p}$} (pu)
        ;
      \end{tikzpicture}
    \end{equation*}
    by \autoref{lemma:square-transform} and the functor $\dC_{/p}$ is a local
    $\sR$-functor by \autoref{rem:over-u-functorial}. As $\dC[\Sigma] \to
    \dC[\Pi]$ has the left lifting property against all local $\sR$-functors,
    we find a lift $\ell \colon \dC[\Pi] \to \dC[\Sigma \hc \Pi]_{/u}$ in this
    diagram. However, any such lift $\ell$ induces a lift $\ell_{s,t}$ in the
    original diagram of simplicial sets by \autoref{lemma:solution-transport}.
  \end{proof}

  \section{Necessity of the hypothesis, part 2}\label{sec:necessity2}

  This section is devoted to a proof of the following third proposition
  announced at the beginning of this chapter:
  \liftingrestriction*

  The proof of \autoref{prop:lifting-restriction} is very similar in spirit to the proof of \autoref{prop:lifting-necessary} in the preceding section in the sense that we transform some given lifting problem of the form
  \begin{equation}\label{dia:restriction-lifting}
    \begin{tikzpicture}[diagram]
      \matrix[objects] {
        |(s)| \dC[\Sigma_{x,y}] \& |(b)| \dB \\
        |(p)| \dC[\Pi_{x,y}] \& |(a)| \dA \\
      };
      \path[maps,->] 
        (s) edge node[above] {$u$} (b)
        (s) edge (p)
        (p) edge node[below] {$v$} (a)
        (b) edge node[right] {$p$} (a)
      ;
    \end{tikzpicture}
  \end{equation}
  into a lifting problem 
  \begin{equation}\label{dia:lozenge-lifting}
    \begin{tikzpicture}[diagram]
      \matrix[objects] {
        |(s)| \dC[\Sigma] \& |(b)| \dB_\lozenge \\
        |(p)| \dC[\Pi] \& |(a)| \dA_\lozenge \\
      };
      \path[maps,->] 
        (s) edge (b)
        (s) edge (p)
        (p) edge (a)
        (b) edge node[right] {$p_\lozenge$} (a)
      ;
    \end{tikzpicture}
  \end{equation}
  in such a way that $p_\lozenge$ is a local $\sR$-functor whenever $p$ is.
  Moreover, we show that any solution to the lifting problem
  \eqref{dia:lozenge-lifting} can be translated back into a solution to the
  lifting problem \eqref{dia:restriction-lifting}. However, the technical details are completely different. 

  \begin{definition}\label{def:lozenge}
  Let $\dA$ be a simplicial category. 
  We define a simplicial category $\dA_\lozenge$ as follows:
    The set of objects of\/ $\dA_\lozenge $ is the disjoint
      union of the set of objects of $\dA$ and $\{s,t\}$.
    The mapping spaces of\/ $\dA_\lozenge$ are given by
      \begin{equation*}
        \dA_\lozenge(x,y) = 
        \begin{cases}
          \dA(x,y) &\text{if } x,y \in \dA \\
          \Delta^0 &\text{if } x = s \text{ or } y = t \\
          \emptyset &\text{otherwise}
        \end{cases}
      \end{equation*}
      and the compositions are those of\/ $\dA$ and the maps determined by the
      universal properties of the initial object $\emptyset$ or the terminal
      object $\Delta^0$ of\/ $\sset$.
  \end{definition}

  \begin{remark}\label{rem:lozenge-functoriality}
    The category $\dA_\lozenge$ should be pictured as the category $\dA$ with
    an initial object $s$ and a terminal object $t$ freely adjoined. 
    Consistent with this picture, the assignment $\dA \mapsto \dA_\lozenge$ is functorial in $\dA$ and there are obvious inclusions $\dA \to \dA_\lozenge$ such that the diagrams
    \begin{equation*}
      \begin{tikzpicture}[diagram]
        \matrix[objects] {
          |(a)| \dA \& |(al)| \dA_\lozenge \\
          |(b)| \dB \& |(bl)| \dB_\lozenge \\
        };
        \path[maps,->] 
          (a) edge (al)
          (a) edge node[left] {$f$} (b)
          (b) edge (bl)
          (al) edge node[right] {$f_\lozenge$} (bl)
        ;
      \end{tikzpicture}
    \end{equation*}
    commute. Moreover, as $f_\lozenge$ is the identity on all mapping spaces
    $\dA(x,y)$ with $\{x,y\} \nsubseteq \dA$, it is immediate that $f \mapsto
    f_\lozenge$ takes local $\sR$-functors to local $\sR$-functors.
  \end{remark}

  Let us now consider a span
    \begin{equation*}
      \begin{tikzpicture}[diagram]
        \matrix[objects] {
          |(a)| \dA \& |(b)| \dB \& |(c)| \dC \\
        };
        \path[maps,->] 
          (b) edge node[above] {$u$} (a)
          (b) edge node[above] {$i$} (c)
        ;
      \end{tikzpicture}
    \end{equation*}
    of simplicial categories with $i$ fully faithful and injective on objects.
    We want to construct a functor $\what{u} \colon \dC \to \dA_\lozenge$ such
    that
    \begin{equation*}
      \begin{tikzpicture}[diagram]
        \matrix[objects] {
          |(xy)| \dB \& |(a)| \dA \\
          |(s)| \dC \& |(al)| \dA_\lozenge\\
        };
        \path[maps,->] 
          (xy) edge node[above] {$u$} (a)
          (xy) edge node[left] {$i$} (s)
          (a) edge (al)
          (s) edge node[below] {$\what{u}$} (al)
        ;
      \end{tikzpicture}
    \end{equation*}
    commutes. For this purpose, we consider a partition $\ob \dC = C_0 \cup C_1
    \cup C_2$, $C_1 = \ob \dB$,  of the objects of $\dC$ such that $\dC(c,c') =
    \emptyset$ whenever $x \in C_i$, $y \in C_j$ with $j < i$.  We then define
    $\what{u}\colon \dC \to \dA_\lozenge$ by 
    \begin{equation*}
      \what{u}(c) =  \begin{cases}
        s & \text{if } c \in C_0 \\
        uc & \text{if } c \in \ob\dB \\
        t & \text{if } c \in C_2
      \end{cases}
    \end{equation*}
    on objects. We then have no choice but to define $\what{u} \colon \dC(c,c') \to \dA_\lozenge( \what{u}c, \what{u}c')$ on mapping spaces by
    \begin{equation*}
      \left\{\
        \begin{aligned}
          &\dC(c,c') = \dB(c,c') \xto{u} \dA(uc,uc') = \dA_\lozenge(\what{u}c,\what{u}c') && \quad \text{if } c,c ' \in \ob \dB, \\
          &\dC(c,c') \to \Delta^0 = \dA_\lozenge(\what{u}c,\what{u}c')&&\quad
          \text{if } c \in C_0 \text{ or } c' \in C_2, \\
          &\dC(c,c') = \emptyset \to \dA_\lozenge(\what{u}c,\what{u}c')&&
          \quad\text{otherwise}.
        \end{aligned} \right.
    \end{equation*}
    It is now easily verified that this assignment indeed defines a simplicial
    functor $\what{u} \colon \dC \to \dA_{\lozenge}$. Moreover, the
    construction of $\what{u}$ from a given partition $\dC = C_0 \cup C_1 \cup
    C_2$ is functorial in a certain sense:
    \begin{remark}
      Consider a morphism
    \begin{equation*}
      \begin{tikzpicture}[diagram]
        \matrix[objects] {
          |(a)| \dA \& |(b)| \dB \& |(c)| \dC \\[-2ex]
          |(ap)| \dA' \& |(bp)| \dB' \& |(cp)| \dC' \\
        };
        \path[maps,->] 
          (b) edge node[above] {$u$} (a)
          (b) edge node[above] {$i$} (c)
          (bp) edge node[below] {$v$} (ap)
          (bp) edge node[below] {$j$} (cp)
          (a) edge node[left] {$f$} (ap)
          (b) edge node[left] {$g$} (bp)
          (c) edge node[right] {$h$} (cp)
        ;
      \end{tikzpicture}
    \end{equation*}
    of spans. Suppose that $i$ and $j$ are injective on objects and fully
      faithful. Further assume that $g$ and $h$ are bijective on objects.
      Finally, let  
      \begin{equation*}
        \ob \dC' = C_0' \cup C_1' \cup C_2', \quad C_1' = \ob \dB'
      \end{equation*}
      be a partition of the objects of $\dC'$ such that $\dC'(c,c') =
      \emptyset$ whenever $c \in C_i$, $c' \in C_j$ and $j < i$. By the
      preceding discussion, this partition gives rise to a functor $\what{v}
      \colon \dC' \to \dA'_\lozenge$.

      But we also obtain a partition $\ob\dC = C_0 \cup C_1 \cup C_2$ with $C_i = h^{-1}C_i'$ as $g$ and $h$ are bijective on objects. It is now easy to see that this partition gives rise to a functor $\what{u} \colon \dC \to \dA_\lozenge$.
      Moreover, looking at the construction of $\what{u}$ and $\what{v}$, one easily verifies that the square
    \begin{equation*}
    \begin{tikzpicture}[diagram]
        \matrix[objects] {
          |(s)| \dC \& |(al)| \dA_\lozenge\\[-2ex]
          |(sp)| \dC' \& |(alp)| \dA'_\lozenge\\
        };
        \path[maps,->] 
          (s) edge node[above] {$\what{u}$} (al)
          (sp) edge node[below] {$\what{v}$} (alp)
          (al) edge node[right] {$f_\lozenge$} (alp)
          (s) edge node[left] {$h$} (sp)
        ;
      \end{tikzpicture}
    \end{equation*}
    commutes.
  \end{remark}

  In order to put this digression on the properties of the assignment $\dA
  \mapsto \dA_\lozenge$ to good use, we have to be able to produce certain
  partitions of the objects of categories like $\dC[\Sigma]$. This is achieved
  by the following lemma:
  \begin{lemma}\label{lem:vertex-partition}
    Let $G$ be a globular graph. For any two vertices $x,y \in G$ there exists
    a partition $V = V_0 \cup V_1 \cup V_2$ such that $V_1$ is the set of
    vertices of $G_{x,y}$ and such that $G$ contains no paths $p$ from $u$ to
    $v$ with $u \in V_i$, $v \in V_j$ and $j < i$.
  \end{lemma}
  \begin{proof}
    Let $V$ denote the set of vertices of $G$ and let $V_1$ be the set of vertices of $G_{x,y}$.
    Let $V_0$ be the set of vertices $v$ of $V \smallsetminus V_1$ such
    that there is a directed path $p$ from $v$ to some vertex $w \in V_1$.
    Observe that there is no path from some vertex $w' \in V_1$ to some
    vertex $v \in V_0$, for one could then fabricate a path from $x$ to $y$
    through $v$, see \autoref{fig:no-v1-back}.
    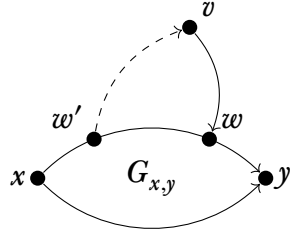
\begin{figure}
          \centering
        \begin{tikzpicture}
            \node[fill,circle,inner sep=2pt] (x) at (2,0) {};
            \node[fill,circle,inner sep=2pt] (y) at (5,0) {};
            \node[fill,circle,inner sep=2pt] (v) at (4,2) {};
            \node[left] at (x) {$x$};
            \node[right] at (y) {$y$};
            \node[above right] at (v) {$v$};
            \path (x) edge[->,out=45,in=135] coordinate[pos=0.25] (wpos) coordinate[pos=0.5] (oben) coordinate[pos=0.75] (wppos) (y);
            \path (x) edge[->,out=-45,in=-135] coordinate (unten) (y);
            \node at ($ (oben.south) ! 0.5 !(unten.north) $) {$G_{x,y}$};
            \node[fill,circle,inner sep=2pt] (w) at (wpos) {};
            \node[fill,circle,inner sep=2pt] (wp) at (wppos) {};
            \node[above left] at (w) {$w'$};
            \node[above right] at (wp) {$w$};
            \path (w) edge[dashed,bend left,->] (v) 
              (v) edge[bend left,->] (wp);
        \end{tikzpicture}
          \caption{A vertex $v \in V_0$ with its path to some vertex $w \in G_{x,y}$. If the dashed path exists, then $v \in G_{x,y}$.}\label{fig:no-v1-back}
    \end{figure}

    Now let $V_2$ be the complement of $V_0 \cup V_1$ in $V$. Note that there
    cannot be any directed path from a vertex $v \in V_2$ to some vertex $w \in
    V_0 \cup V_1$ since this would imply $w \in V_0$.

    Altogether, we have thus found a partition $V_0 \cup V_1 \cup V_2$ with the
    desired properties. 
  \end{proof}

  Let us now return to our actual concern, namely the promised proof of the
  following proposition:
  \liftingrestriction
  \begin{proof}
    Consider a lifting problem
    \begin{equation*}
    \begin{tikzpicture}[diagram]
      \matrix[objects] {
        |(s)| \dC[\Sigma_{x,y}] \& |(b)| \dB \\
        |(p)| \dC[\Pi_{x,y}] \& |(a)| \dA \\
      };
      \path[maps,->] 
        (s) edge node[above] {$u$} (b)
        (s) edge (p)
        (p) edge node[below] {$v$} (a)
        (b) edge node[right] {$p$} (a)
      ;
    \end{tikzpicture}
  \end{equation*}
    with $p$ a local $\sR$-functor. 
    We then get a morphism 
    \begin{equation*}
      \begin{tikzpicture}
       \matrix[objects] {
         |(b)| \dB \& |(s)| \dC[\Sigma_{x,y}] \& |(sfull)| \dC[\Sigma] \\
        |(a)| \dA \& |(p)| \dC[\Pi_{x,y}] \& |(pfull)| \dC[\Pi] \\
      };
      \path[maps,->] 
        (s) edge node[above] {$u$} (b)
        (s) edge (p)
        (p) edge node[below] {$v$} (a)
        (b) edge node[left] {$p$} (a)
        (s) edge node[above] {$i$} (sfull)
        (p) edge node[below] {$j$} (pfull)
        (sfull) edge (pfull)
      ;
    \end{tikzpicture}
    \end{equation*}
    of spans in which the two vertical maps on the right hand side are
    bijective on objects and in which $i$ and $j$ are injective on objects and
    fully faithful. The partition of $\ob \dC[\Pi] = \ob \dC[\Sigma]$ constructed in \autoref{lem:vertex-partition} now gives rise to the commutative square on the right hand side of the diagram
\begin{equation*}
    \begin{tikzpicture}[diagram]
      \matrix[objects] {
        |(sr)| \dC[\Sigma_{x,y}] \&  |(s)| \dC[\Sigma] \& |(b)| \dB_\lozenge \\
        |(pr)| \dC[\Pi_{x,y}] \& |(p)| \dC[\Pi] \& |(a)| \dA_\lozenge\smash{.} \\
      };
      \path[maps,->] 
        (sr) edge (s)
        (pr) edge (p)
        (sr) edge (pr)
        (s) edge node[above] {$\what{u}$} (b)
        (s) edge (p)
        (p) edge node[below] {$\what{v}$}  (a)
        (b) edge node[right] {$p_\lozenge$} (a)
      ;
    \end{tikzpicture}
\end{equation*}
    Note that $p_\lozenge$ is a local $\sR$-functor by
    \autoref{rem:lozenge-functoriality}.  Note further that the square on the
    left in this diagram commutes, too, and we thus find a solution to the
    original lifting problem that we started with.
  \end{proof}

  \section{The proof of Theorem B}\label{sec:proof-b}

  Finally, we have all pieces together to give the proof of 
  \globalthm

  \begin{proof}
    Sufficiency of the condition in the theorem is
    \autoref{prop:lifting-sufficient} and necessity can be deduced as follows:
    Let $\Sigma \to \Pi$ be an inclusion between complete pasting diagrams such
    that $\dC[\Sigma] \to \dC[\Pi]$ has the left lifting property against all
    local $\sR$-functors. According to \autoref{prop:lifting-restriction}, the
    functors $\dC[\Sigma_{x,y}] \to \dC[\Pi_{x,y}]$, where $x,y \in \Sigma$ are
    arbitrary vertices, then have the left lifting property against all local
    $\sR$-functors, too. It now follows from \autoref{prop:lifting-necessary}
    that the maps $\nerve(\Sigma_{x,y} \hc \Pi_{x,y} ) \to \nerve(\Pi_{x,y})$
    are $\sL$-maps for all vertices $x,y \in \Sigma$.
  \end{proof}

\chapter{Local Lifting Properties of Pasting Diagrams}\label{chap:local}

This whole chapter is devoted to a proof of 
\localthm

Our proof of \autoref{thm:local} closely parallels Power's proof of his
$2$-categorical pasting theorem in
\cite{power;a-2-categorical-pasting-theorem}. More precisely, both proofs work
by induction on the size of the underlying graph $G$ and both inductive steps
rest on \autoref{lemma:glob-graph:face}. In our case, there are more technical
difficulties to be overcome, though.

We explain the basic strategy and handle some trivial cases of the proof in \autoref{sec:setup}. Afterwards, in \autoref{sec:induction} we reduce the problem to a combinatorial problem that we solve in \autoref{sec:fillable}

\section{The basic setup}\label{sec:setup}

Let us give a quick outline of the proof of \autoref{thm:local}. We proceed by
induction on the number of edges of the graph $G$ underlying both $\Sigma$ and $\Pi$.
The base clause is basically trivial. For the
induction step, we reduce to the case that~$G$ is $2$-connected by decomposing
complete pasting diagrams on non $2$-connected graphs as
\begin{equation*}
  \nerve(\Sigma_1 \join \Sigma_2) \iso \nerve(\Sigma_1) \times \nerve(\Sigma_2).
\end{equation*}
The assumption that $G$ is $2$-connected then allows us to avail ourselves of
\autoref{lemma:glob-graph:face}, and construct certain globular subgraphs
$G_0,G_1,G_2$ of $G$ that are smaller than $G$. We may thus apply our induction
hypothesis to the restriction of the pasting diagrams $\Sigma$ and $\Pi$ in
\autoref{thm:local}. The mid anodyne maps thus obtained can be glued so as to
obtain a factorisation
\begin{equation*}
  \nerve(\Sigma) \to X_0 \to \nerve(\Pi)
\end{equation*}
in which $\nerve(\Sigma) \to X_0$ is mid anodyne. In order to finish the proof
of \autoref{thm:local}, we thus have to exhibit $X_0 \to \nerve(\Pi)$ as mid
anodyne, too. We solve this problem in \autoref{sec:fillable} by a careful
analysis of those simplices of $\nerve(\Pi)$ that are not contained in $X_0$.
This analysis eventually leads to a filtration of $X_0 \to \nerve(\Pi)$ by
pushouts of inclusions $\Lambda^n_i \to \Delta^n$ with $0 < i < n$. The map $X_0 \to \nerve(\Pi)$ is thus mid anodyne and this finishes the proof.

\begin{remark}\label{rem:vertical-composition-closed-under-stuff}
  Observe that the inclusion $\Sigma_{\min}^c \to \Pi_{\max}$ of the minimal complete pasting diagram into the maximal pasting diagram satisfies the hypothesis of \autoref{thm:local}.
  
  Moreover, if $\Sigma \to \Pi$ is an inclusion satisfying the hypothesis of
  \autoref{thm:local}, then so do the inclusions $\Sigma \to \Sigma \hc \Pi$,
  $\Sigma \hc \Pi \to \Pi$ and all restrictions of $\Sigma \to \Pi$ to globular
  subgraphs of the form $G_{x,y}$ of the underlying graph $G$.
\end{remark}

As announced above, we do induction on the size of $G$, the base clause being trivial:
\begin{lemma}\label{lem:induction-one-interior-face}
  The map $\nerve(\Sigma) \to \nerve(\Pi)$ is mid anodyne for all inclusions
  $\Sigma  \to \Pi$ satisfying the hypotheses of \autoref{thm:local} on some
  globular graph $G$ with at most one interior face.
\end{lemma}
\begin{proof}
  If $G$ has no interior face, there is nothing to show as $G$ is nothing but a
  possibly trivial directed path. Similarly, if $G$ has only one interior face,
  then $\sS = \sT$ and the map in question is an
  identity.
\end{proof}
\begin{lemma}\label{lem:induction-cut-vertices}
  Suppose that\/ $G$ decomposes as\/ $G = G_1 \join G_2$. For\/ $i \in \{1,2\}$
  denote the restrictions of\/ $\Sigma$ and\/ $\Pi$ to\/
  $G_i$ by\/ $\Sigma_i$ and\/ $\Pi_i$, respectively. If the maps
  \begin{equation*}
    \nerve(\Sigma_i) \to \nerve(\Pi_i) , \quad i \in \{1,2\},
  \end{equation*}
  are mid anodyne, then so is the map $\nerve(\Sigma) \to \nerve(\Pi)$.
\end{lemma}
\begin{proof}
  Follows from the assumption that $\Sigma$ and $\Pi$ are complete,
  \autoref{cor:nerve-of-join} and the fact that mid anodyne maps are closed
  under taking products. 
\end{proof}

\section{The inductive step}\label{sec:induction}

Given \autoref{lem:induction-cut-vertices}, we may suppose that $G$ is
$2$-connected and has at least two faces.  According to
\autoref{lemma:glob-graph:face}, there then exists a face $\phi$ of $G$ such
that $\dom\phi \subseteq \dom G$. Consider the following subgraphs $G_i$ of
$G$:
\begin{enumerate}
  \item $G_1$ is the subgraph of $G$ consisting of all directed paths
    containing either $\dom \phi$ or $\codom \phi$ as subpath.
  \item $G_2$ is the subgraph of $G$ consisting of all directed paths that
    contain no edge of $\dom\phi$, that is, $G_2$ is obtained from $G$ by removing the interior vertices of $\dom(\phi)$.
  \item $G_0$ is the intersection of $G_1$ and $G_2$, i.\,e.\ the subgraph of
    $G$ consisting of all directed paths containing $\cod \phi$ but no edge of
    $\dom \phi$. 
\end{enumerate}
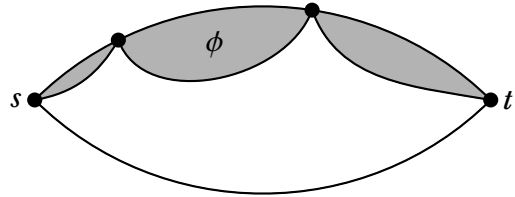
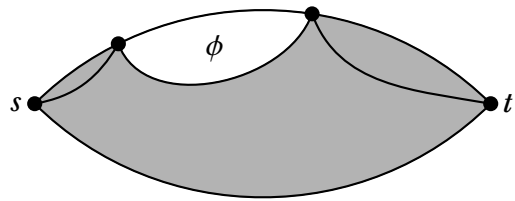
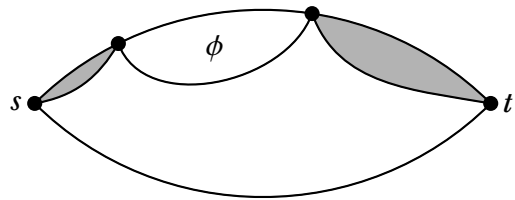
\begin{figure}
  \centering
  \begin{subfigure}[b]{0.8\textwidth}
    \centering
    \begin{tikzpicture}
      \coordinate (s) at (0,0);
      \coordinate (t) at (6,0);
      \node[left] at (s) {$s$};
      \node[right] at (t) {$t$};
      \draw (s.center) edge[out=45,in=135] coordinate[pos=0.2] (sphi) coordinate[pos=0.6] (tphi) coordinate[pos=0.4] (mphi) (t.center);
      \draw (s.center) edge[out=-45,in=-135]  (t.center);
      \draw (sphi.center) edge[out=-70,in=-110] coordinate (mphi2) (tphi.center);
      \draw (s.center) edge[out=10,in=-120] (sphi.center);
      \draw (tphi.center) edge[out=-70,in=170] (t.center);
      \fill[black!30] 
        (t.center) to[out=170,in=-70]  (tphi.center)
            to[out=-110,in=-70] (sphi.center)
            to[out=-120,in=10]  (s.center)
            to[out=45,in=135] (t.center);
      \draw[thick] (s.center) to[out=45,in=135] (t.center);
      \draw[thick] (s.center) to[out=-45,in=-135] (t.center);
      \draw[thick] (sphi.center) to[out=-70,in=-110] (tphi.center);
      \node[] at ($ (mphi) ! 0.5 ! (mphi2) $) {$\phi$};
      \draw[thick] (s.center) to[out=10,in=-120] (sphi.center);
      \draw[thick] (tphi.center) to[out=-70,in=170] (t.center);
      \node[fill,circle,inner sep=2pt] (s) at (0,0) {};
      \node[fill,circle,inner sep=2pt] (t) at (6,0) {}; 
      \node[fill,circle,inner sep=2pt] at (sphi) {};
      \node[fill,circle,inner sep=2pt] at (tphi) {};
  \end{tikzpicture}
    \caption{$G_1$}
  \end{subfigure}
  \vskip\baselineskip
  \begin{subfigure}[b]{0.8\textwidth}
    \centering
    \begin{tikzpicture}
      \coordinate (s) at (0,0);
      \coordinate (t) at (6,0);
      \node[left] at (s) {$s$};
      \node[right] at (t) {$t$};
      \draw (s.center) edge[out=45,in=135] coordinate[pos=0.2] (sphi) coordinate[pos=0.6] (tphi) coordinate[pos=0.4] (mphi) (t.center);
      \draw (s.center) edge[out=-45,in=-135]  (t.center);
      \draw (sphi.center) edge[out=-70,in=-110] coordinate (mphi2) (tphi.center);
      \draw (s.center) edge[out=10,in=-120] (sphi.center);
      \draw (tphi.center) edge[out=-70,in=170] (t.center);
      \fill[black!30] 
        (s.center) to[out=45,in=-155]  (sphi.center)
            to[out=-70,in=-110] (tphi.center)
            to[out=-12,in=135]  (t.center)
            to[out=-135,in=-45] (s.center);
      \draw[thick] (s.center) to[out=45,in=135] (t.center);
      \draw[thick] (s.center) to[out=-45,in=-135] (t.center);
      \draw[thick] (sphi.center) to[out=-70,in=-110] (tphi.center);
      \node[] at ($ (mphi) ! 0.5 ! (mphi2) $) {$\phi$};
      \draw[thick] (s.center) to[out=10,in=-120] (sphi.center);
      \draw[thick] (tphi.center) to[out=-70,in=170] (t.center);
      \node[fill,circle,inner sep=2pt] (s) at (0,0) {};
      \node[fill,circle,inner sep=2pt] (t) at (6,0) {}; 
      \node[fill,circle,inner sep=2pt] at (sphi) {};
      \node[fill,circle,inner sep=2pt] at (tphi) {};
  \end{tikzpicture}
    \caption{$G_2$}
  \end{subfigure}
  \vskip\baselineskip
  \begin{subfigure}[b]{0.8\textwidth}
    \centering
    \begin{tikzpicture}
      \coordinate (s) at (0,0);
      \coordinate (t) at (6,0);
      \node[left] at (s) {$s$};
      \node[right] at (t) {$t$};
      \draw (s.center) edge[out=45,in=135] coordinate[pos=0.2] (sphi) coordinate[pos=0.6] (tphi) coordinate[pos=0.4] (mphi) (t.center);
      \draw (s.center) edge[out=-45,in=-135]  (t.center);
      \draw (sphi.center) edge[out=-70,in=-110] coordinate (mphi2) (tphi.center);
      \draw (s.center) edge[out=10,in=-120] (sphi.center);
      \draw (tphi.center) edge[out=-70,in=170] (t.center);
      \fill[black!30] 
            (tphi.center)
            to[out=-12,in=135]  (t.center)
            to[out=170,in=-70] (tphi.center);
      \fill[black!30] 
            (sphi.center)
            to[out=-120,in=10]  (s.center)
            to[out=45,in=-155] (sphi.center);
      \draw[thick] (s.center) to[out=45,in=135] (t.center);
      \draw[thick] (s.center) to[out=-45,in=-135] (t.center);
      \draw[thick] (sphi.center) to[out=-70,in=-110] (tphi.center);
      \node[] at ($ (mphi) ! 0.5 ! (mphi2) $) {$\phi$};
      \draw[thick] (s.center) to[out=10,in=-120] (sphi.center);
      \draw[thick] (tphi.center) to[out=-70,in=170] (t.center);
      \node[fill,circle,inner sep=2pt] (s) at (0,0) {};
      \node[fill,circle,inner sep=2pt] (t) at (6,0) {}; 
      \node[fill,circle,inner sep=2pt] at (sphi) {};
      \node[fill,circle,inner sep=2pt] at (tphi) {};
  \end{tikzpicture}
    \caption{$G_0$}
  \end{subfigure}
  \caption{Schematic pictures of the graphs $G_1$, $G_2$ and $G_0$.}\label{fig:g0g1g2}
\end{figure}

In \autoref{fig:g0g1g2} one can see a sketch of how these graphs roughly
look like.  Note that the graphs $G_i$ are globular by
\autoref{prop:ps:power}. Moreover, they all have fewer edges than $G$ ---
a fact that we record in the following lemma.
\begin{lemma}\label{lem:gi-globular-and-smaller}
  The graphs $G_i$ have fewer edges than $G$.
\end{lemma}
\begin{proof}
  As $\dom \phi$ consists of at least one edge, both $G_0$ and $G_2$ have fewer
  edges than $G$. Finally, if $G$ and $G_1$ had the same number of edges, then
  $s(\phi)$ or $t(\phi)$ would be a cut vertex of $G$ as we assume $G$ to have
  at least two faces.  This contradicts our assumption that $G$ is
  $2$-connected, though.
\end{proof}

For $i \in \{0,1,2\}$ denote by $\Sigma_i$ and $\Pi_i$ the restrictions of
$\Sigma$ and $\Pi$ to $G_i$.  We gather the maps
$\nerve(\Sigma_i) \to \nerve(\Pi_i)$ induced on the respective
nerves in the commutative diagram
\begin{equation}\label{eq:smaller-nerve-maps}
  \begin{tikzpicture}[diagram]
    \matrix[objects] {
      |(g1-f)| \nerve(\Sigma_1) \& |(g0-f)| \nerve(\Sigma_0)  \& |(g2-f)| \nerve(\Sigma_2) \\
      |(g1-t)| \nerve(\Pi_1 ) \& |(g0-t)| \nerve(\Pi_0)  \& |(g2-t)| \nerve(\Pi_2) \smash{.} \\
    };
    \path[maps,->] 
      (g0-f) edge (g1-f)
      (g0-f) edge (g2-f)
      (g0-t) edge (g1-t)
      (g0-t) edge (g2-t)
      (g0-f) edge (g0-t)
      (g1-f) edge (g1-t)
      (g2-f) edge (g2-t)
    ;
  \end{tikzpicture}
\end{equation}
\begin{lemma}\label{lem:nerve-pushout-in-proof}
  If the maps $\nerve(\Sigma') \to \nerve(\Pi')$ are mid anodyne for all
  inclusions $\Sigma' \to \Pi'$ satisfying the hypotheses of
  \autoref{thm:local}  on globular graphs $G'$ with fewer edges than $G$, then
  the map 
  \begin{equation*}
    \nerve(\Sigma_1) \coprod_{\nerve(\Sigma_0)} \nerve(\Sigma_2) \to 
    \nerve(\Pi_1) \coprod_{\nerve(\Pi_0)} \nerve(\Pi_2) 
  \end{equation*} 
  induced by \eqref{eq:smaller-nerve-maps}
  is mid anodyne, too.
\end{lemma}
\begin{proof}
  The maps 
  \begin{equation*}
    \nerve(\Sigma_i) \coprod_{\nerve( \Sigma_0)} \nerve(\Pi_0) \to \nerve(\Pi_i)
  \end{equation*}
  with $i \in \{1,2\}$ are mid anodyne by the stability of mid anodyne maps
  under pushout and the fact that mid anodyne maps have the right cancellation
  property, see
  \cite[Theorem~E]{stevenson;stability-for-inner-fibrations-revisited}. Using
  this observation, one may check the desired lifting property by hand.
\end{proof}

Now consider the pushout 
\begin{equation}\label{eq:pushout-in-proof}
  \begin{tikzpicture}[diagram]
    \matrix[objects] {
      |(colims)| \colim \nerve(\Sigma_\ph) \& |(s)| \nerve(\Sigma) \\
      |(colimp)| \colim \nerve(\Pi_\ph) \& |(p)|  \displaystyle\nerve(\Sigma)\ \coprod_{ \mathclap{\colim\nerve(\Sigma_\ph)}}\ \colim \nerve(\Pi_\ph) \\
    };
    \path[maps,->] 
      (colims) edge (colimp)
      (colims) edge (s)
      (colimp) edge (p)
      (s) edge (p)
    ;
  \end{tikzpicture}
\end{equation}
of the map in \autoref{lem:nerve-pushout-in-proof} along the
inclusion $\colim \nerve(\Sigma_\ph) \to \nerve(\Sigma)$. 
The vertical map on the right hand side in \eqref{eq:pushout-in-proof} is mid anodyne by the stability of
mid anodyne maps under pushouts and \autoref{lem:nerve-pushout-in-proof}.  Moreover, the inclusion $\nerve(\Sigma) \to
\nerve(\Pi)$ factors through this map and to finish the proof of
\autoref{thm:local}, it therefore suffices to show that the canonical map
\begin{equation*}
  \nerve(\Sigma) \  \coprod_{ \mathclap{ \colim \nerve(\Sigma_\ph)  } }\ \colim
  \nerve(\Pi_\ph ) \to \nerve(\Pi) 
\end{equation*}
is mid anodyne. This will be accomplished in the next section.

\section{Fillable simplices}\label{sec:fillable}

In order to finish the proof of \autoref{thm:local}, we have to show that the map
\begin{equation}\label{eq:last-map}
  \nerve(\Sigma) \  \coprod_{ \mathclap{ \colim \nerve(\Sigma_\ph)  } }\ \colim
  \nerve(\Pi_\ph ) \to \nerve(\Pi) 
\end{equation}
induced by \eqref{eq:smaller-nerve-maps} is mid anodyne.
We will achieve this by a direct inspection of the
simplices of the domain and codomain of \eqref{eq:last-map}.  For this
purpose, let us first recollect the construction of $G_1$. We chose a
face $\phi$ of $G$ with $\dom \phi \subseteq \dom G$ and let $G_1$ be the
subgraph of $G$ that consists of all paths that contain either $\dom
\phi$ or $\cod \phi$ as subpath.  It is immediate from this description
that we have the following characterisation of the simplices of
$\nerve(\Pi_1)$ considered as a simplicial subset of $\nerve(\Pi)$:
\begin{lemma}\label{lem:g1-simplices}
  Consider an $n$-simplex\/ $\sigma = (p_0 \leq \dots \leq p_n)$
  of\/~$\nerve(\Pi)$.  If all the paths\/ $p_i$ contain either\/ $\dom
  \phi$ or\/ $\cod \phi$, then\/ $\sigma \in \nerve(\Pi_1)$.
\end{lemma}
The graph $G_2$ was defined as the subgraph of $G$ that consists of all paths
that have no edge in common with $\dom \phi$. The following lemma and its
corollary provide us with a description of the simplices of
$\nerve(\Pi_2)$.
\begin{lemma}\label{lem:g2-simplices}
  Consider paths\/ $p \leq q$ in\/ $G$ and suppose that\/ $\dom \phi$ and\/ $p$
  have no common edge.  Then the same holds true for\/ $\dom \phi$ and\/ $q$. 
\end{lemma}
\begin{proof}
  We may assume $p < q$. For any path $q$ with $p < q$, we
  may write $p = a \cdot \dom \gamma \cdot b$ and $q = a \cdot \cod \gamma
  \cdot b$ for a nontrivial glob $\gamma$ and suitable paths $a$ and $b$. It
  follows that if $q$ and $\dom \phi$ had a common edge, this edge would be an
  edge in the codomain of some nontrivial glob $\gamma$.  This is impossible,
  though, as $\dom \phi \subseteq \dom G$ and no edge of $\dom G$ occurs in the
  codomain of a nontrivial glob.
\end{proof}
\begin{corollary}\label{cor:g2-simplices}
  Consider an $n$-simplex\/ $\sigma = (p_0 \leq \dots \leq p_n)$
  of\/~$\nerve(\Pi) $.  If\/~$\dom \phi$ and\/ $p_0$ have no common edge, then\/
  $\sigma \in \nerve(\Pi_2)$.
\end{corollary}
In fact, \autoref{lem:g1-simplices} and \autoref{lem:g2-simplices} tell us a
bit more about the structure of the simplices of $\nerve(\Pi)$:
\begin{corollary}\label{cor:g-simplices}
  Consider an $n$-simplex\/ $\sigma = (p_0 \leq \dots \leq p_n)$
  of\/~$\nerve(\Pi)$ with minimal witnesses\/ $\gamma_i$ for\/ $p_{i-1} \leq
  p_i$.  Then either\/ $\sigma \in \colim \nerve(\Pi_\ph)$
  or\/~$\dom \phi \subseteq p_0$ and\/ $\gamma_i \nsubseteq G_1$ for some\/ $i
  \in \{1,\dots,n\}$.
\end{corollary}
For a given $n$-simplex 
\begin{equation*}
  \sigma = (p_0 \leq \dots \leq p_n)
\end{equation*}
of $\nerve(\Pi)$ with minimal witnesses $\gamma_i$ for $p_{i-1} \leq p_i$, we
let $c(\sigma)$ be the minimal $c \in \{ 1,\dots,n\}$ such that $\gamma_c
\nsubseteq G_1$. If there is no such $c$, we let $c(\sigma) = n+1$ by
convention.  We call the simplex \emph{fillable} if $c(\sigma) = n+1$ or
$\gamma_{c(\sigma)} \cap G_1 \subseteq \del G_1$.  The following example
illustrates the definition of fillable simplices:
\begin{example}
  In order to illustrate the definition of a fillable simplex, let us draw
  simplices $\sigma$ in $\nerve(\Pi)$ as $n$-marked subgraphs of $G$.
  \begin{enumerate}[label=(\alph*)]
    \item 
       An $n$-simplex $\sigma$ is certainly fillable, whenever $c(\sigma) =
       n+1$, i.\,e.\ whenever all $\gamma_i \subseteq G_1$. This does not
       necessarily imply that $P_\sigma = \bigcup_i p_i \subseteq G_1$ but it
       does imply that all the interior faces of $P_\sigma$ lie in the interior
       of $G_1$, see \autoref{subfig:trivial-fillable}.
     \item
       An $n$-simplex $\sigma$ has $\gamma(c) \leq n$ if and only if there is
       some witness $\gamma_i$ such that $\gamma_i \nsubseteq G_1$. The simplex
       shown in \autoref{subfig:nontrivial-fillable} is fillable as the face of
       $P_\sigma$ labeled $1$ intersects $G_1$ only in its boundary.
     \item
       Finally, consider the face $d_1 \sigma$ of the simplex $\sigma$ of the
       previous example. Its $1$-marked subgraph is shown in
       \autoref{subfig:non-fillable}. We obviously have $c(\sigma) = 1$ but the
       minimal glob $\gamma_1$ witnessing $p_0 \leq p_1$ has nontrivial
       intersection with the interior of $G_1$. This simplex is thus non-fillable.
  \end{enumerate}
 \begin{figure}
  \centering
  \begin{subfigure}[b]{0.8\textwidth}
    \centering
    \begin{tikzpicture}[scale=1]
      \coordinate (s) at (0,0);
      \coordinate (t) at (6,0);
      \fill[black!15] 
        (s.center) to[out=45,in=135] (t.center) -- (s.center); 
      \node[left] at (s) {$s$};
      \node[right] at (t) {$t$};
      \draw[thick,dashed] (s.center) edge[out=45,in=135] coordinate (oben) (t.center);
      \draw[thick,dashed] (s.center) edge[out=-45,in=-135] coordinate (mitte) (t.center);
      \draw[thick,dashed] (s.center)  --  (t.center);
      \draw[very thick] (s.center) to[out=-30,in=-150] ($ (s) ! 0.3 ! (t) $);
      \path[draw,very thick] ($ (s) ! 0.3 ! (t) $) to[out=35,in=90] coordinate[pos=0.4] (x) coordinate[pos=0.9] (y) ($ (s) ! 0.7 ! (t) $);
      \path[draw,very thick] (x) to[out=-30,in=-180] coordinate[pos=0.2] (z) ($ (s) !0.7 !(t) $);
      \path[draw,very thick] (z) -- (y);
      \path[draw,very thick] ($ (s) !0.7 ! (t) $) to[out=-90,in=-150] (t);
  \end{tikzpicture}
    \caption{Sketch of the $n$-marked subgraph of a fillable $n$-simplex $\sigma$ with $c(\sigma) = n+1$.}\label{subfig:trivial-fillable}
  \end{subfigure}
  \vskip.5ex
  \begin{subfigure}[b]{0.8\textwidth}
    \centering
    \begin{tikzpicture}[]
      \coordinate (s) at (0,0);
      \coordinate (t) at (6,0);
      \fill[black!15] 
        (s.center) to[out=45,in=135] (t.center) -- (s.center); 
      \node[left] at (s) {$s$};
      \node[right] at (t) {$t$};
      \draw[thick,dashed] (s.center) edge[out=45,in=135] coordinate (oben) (t.center);
      \draw[thick,dashed] (s.center) edge[out=-45,in=-135] coordinate (mitte) (t.center);
      \draw[thick,dashed] (s.center)  --  (t.center);
      \draw[very thick] (s.center) to[out=-45,in=-135] ($ (s) ! 0.2 ! (t) $);
      \draw[very thick] ($ (s) ! 0.2 ! (t) $) to[out=45,in=135] ($ (s) ! 0.6 ! (t) $);
      \draw[very thick] ($ (s) ! 0.6 ! (t) $) to[out=-45,in=-135] (t);
      \draw[very thick] ($ (s) ! 0.2 ! (t) $) -- (t);
      \node[above,font=\scriptsize] at ($ (s) ! 0.4 ! (t) $) {$2$};
      \node[below,font=\scriptsize] at ($ (s) ! 0.8 ! (t) $) {$1$};
  \end{tikzpicture}
    \caption{The $2$-marked subgraph of a fillable $2$-simplex $\sigma$ with $c(\sigma) = 1$.}\label{subfig:nontrivial-fillable}
  \end{subfigure}
  \vskip.5ex
  \begin{subfigure}[b]{0.8\textwidth}
    \centering
    \begin{tikzpicture}[]
      \coordinate (s) at (0,0);
      \coordinate (t) at (6,0);
      \fill[black!15] 
        (s.center) to[out=45,in=135] (t.center) -- (s.center); 
      \node[left] at (s) {$s$};
      \node[right] at (t) {$t$};
      \draw[thick,dashed] (s.center) edge[out=45,in=135] coordinate (oben) (t.center);
      \draw[thick,dashed] (s.center) edge[out=-45,in=-135] coordinate (mitte) (t.center);
      \draw[thick,dashed] (s.center)  --  (t.center);
      \draw[very thick] (s.center) to[out=-45,in=-135] ($ (s) ! 0.2 ! (t) $);
      \draw[very thick] ($ (s) ! 0.2 ! (t) $) to[out=45,in=135] ($ (s) ! 0.6 ! (t) $);
      \draw[very thick] ($ (s) ! 0.6 ! (t) $) to[out=-45,in=-135] (t);
      \draw[very thick] ($ (s) ! 0.2 ! (t) $) -- (t);
      \node[above,font=\scriptsize] at ($ (s) ! 0.4 ! (t) $) {$1$};
      \node[below,font=\scriptsize] at ($ (s) ! 0.8 ! (t) $) {$1$};
  \end{tikzpicture}
    \caption{The $1$-marked subgraph of the non-fillable face $d_1 \sigma$ of the simplex whose $2$-marked subgraph is shown in \autoref{subfig:nontrivial-fillable}.}\label{subfig:non-fillable}
  \end{subfigure}
  \caption{Illustration of the definition of a fillable simplex. The shaded area of the shown globular graph represents $G_1$.}\label{fig:fillable}
\end{figure}
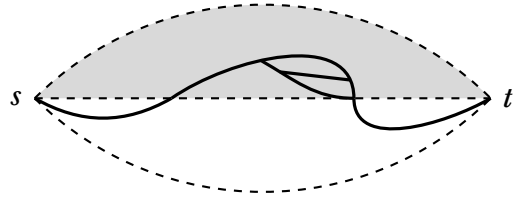
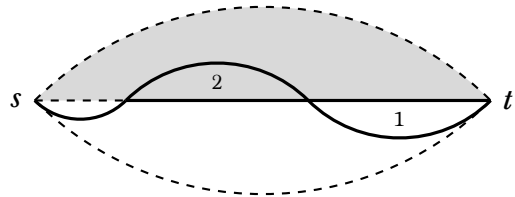
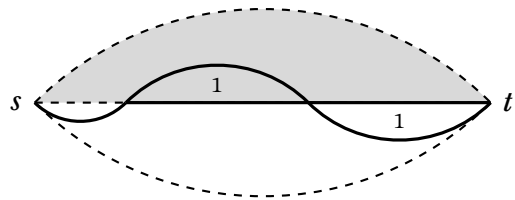

\end{example}

With this definition at hand,
\autoref{cor:g-simplices} admits the following reformulation:
\begin{lemma}\label{cor:fillabel-g-simplices}
  For a fillable simplex\/ $\sigma$ of\/ $\nerve(\Pi)$ of dimension\/~$n$,
  we either have\/~$\sigma \in \colim \nerve(\Pi_\ph)$ or\/~$c(\sigma) \in
  \{2,\dots,n\}$.
\end{lemma}
\begin{proof}
  If $\sigma \notin \colim \nerve(\Pi_\ph)$, then $\sigma \notin
  \nerve(\Pi_2) $ and thus $\dom \phi \subseteq p_0$.  Assume $c(\sigma) = 1$,
  i.\,e.\ $\gamma_1 \nsubseteq G_1$.  Since $\sigma$ is fillable, $\gamma_1$
  intersects $G_1$ only in its boundary. Let us consider the order of the
  vertices $s(\gamma_1), t(\gamma_1), s(\phi), t(\phi)$ on the path $p_0$.
  Observe that the orders
  \begin{equation*}
    s(\phi) \leadsto t(\phi) \leadsto s(\gamma_1) \leadsto t(\gamma_1)
    \qquad\text{and}\qquad
     s(\gamma_1) \leadsto t(\gamma_1) \leadsto s(\phi) \leadsto t(\phi) 
  \end{equation*}
  are impossible as they would imply $\gamma_1 \subseteq G_1$ by our definition
  of $G_1$.
  Moreover, the orders
  \begin{align*}
    &s(\phi)  \leadsto s(\gamma_1) \leadsto  t(\gamma_1) \leadsto t(\phi), &
    t(\gamma_1)  \leadsto t(\phi) \leadsto s(\phi)\leadsto s(\gamma_1)  \\
    &s(\phi) \leadsto s(\gamma_1) \leadsto t(\phi) \leadsto t(\gamma_1), &
     s(\gamma_1) \leadsto s(\phi) \leadsto t(\gamma_1) \leadsto t(\phi)
  \end{align*}
  are also impossible, for $\gamma_1$ is minimal, intersects $G_1$ only in its
  boundary but is not contained in $G_1$ and we have the path $\dom \phi$
  between $s(\phi)$ and $t(\phi)$. However, these are all possible orders and
  we conclude $c(\sigma) \geq 2$.  Moreover, $c(\sigma) \leq n$ as $\sigma
  \notin \nerve(\Pi_1)$.
\end{proof}
In order to finish the proof of \autoref{thm:local}
we want to express the map \eqref{eq:last-map} as a composition of
pushouts of inner horn inclusions $\Lambda^n_i \to \Delta^n$ in which
each $\Delta^n$ maps to a fillable $n$-simplex of $\nerve(\Pi) $. We
therefore need to understand the faces of fillable simplices in
$\nerve(\Pi)$. To this end, we will employ the following remark.
\begin{remark}\label{rem:cutting-simplices-in}
  Consider a glob $\gamma \subseteq G$ with $\gamma \nsubseteq G_1$
  and $\gamma \nsubseteq G_2$. Cutting the glob $\gamma$ along the
  boundary $\del G_1$ supplies us with two globs $\gamma_1$ and
  $\gamma_2$ such that
  \begin{enumerate*}[label=(\roman*)] 
    \item $\gamma = \del( \gamma_1 \cup \gamma_2)$, 
    \item $\gamma_1 \subseteq G_1$ and 
    \item $\gamma_2 \cap G_1 \subseteq \del G_1$. 
  \end{enumerate*}
  This is also illustrated in \autoref{fig:cutting-glob}. Observe
  that $\gamma_1$ corresponds to the part of $\gamma$ that lies left
  of $\del G_1$ while $\gamma_2$ lies right of $\del G_1$.  Given any
  relation $p < r$ witnessed by $\gamma$, we thus find a path $q$
  such that $\gamma_1$ and $\gamma_2$ witness $p < q < r$.
\end{remark}
\begin{remark}\label{cor:any-is-in-boundary-of-fillable}
  The above \autoref{rem:cutting-simplices-in} immediately implies that each
  simplex of $\nerve(\Pi)$ is an inner face of some fillable simplex, for $\Pi
  = (G,\sT)$ is closed under taking subdivisions and the globs $\gamma_1$ and
  $\gamma_2$ constructed above therefore appear in $\sT$ whenever $\gamma \in
  \sT$, see \autoref{lem:nerves-of-subdivision-closed}.
\end{remark}
\begin{figure}
  \centering
    \begin{tikzpicture}
      \begin{scope}
        \clip (0.25,-0.4) to[out=45,in=135] (2.75,-0.4);
        \fill[black!30] (-0.5,0.5) to[out=-45,in=-135] (3.5,0.5);
      \end{scope}
      \draw[name path=delg] (-0.5,0.5) to[out=-45,in=-135] node[pos=0.1,above=1ex] {$\del G_1$} (3.5,0.5);
      \draw[name path=gamma] (0.25,-0.4) to[out=45,in=135] (2.75,-0.4);
      \draw (0.25,-0.4) to[out=-45,in=-135] node[below] (gammalabel) {$\gamma$} (2.75,-0.4);
      \node (gamma1) at (1.5,1) {$\gamma_1$};
      \draw[name intersections={of=delg and gamma}] (gamma1) -- ($ (intersection-1) ! 0.5 ! (intersection-2) $);
      \node (gamma2) at (0,-1) {$\gamma_2$};
      \draw (gamma2) -- ($ (gammalabel) + (0,0.6) $);
    \end{tikzpicture}
  \caption{Cutting a glob $\gamma$ along $\del G_1$.}\label{fig:cutting-glob}
\end{figure}

\newpage
\begin{lemma}\label{lem:fillable-boundaries}
  Let\/ $\sigma$ be a nondegenerate, fillable\/ $n$-simplex of\/
  $\nerve(\Pi)$ and suppose\/ $c = c(\sigma) \in \{2,\dots,n\}$.
  \begin{enumerate}[label=(\alph*)]
    \item The faces $d_i \sigma$ are fillable for all\/ $i \notin \{c-1,c\}$.
    \item The face\/ $d_{c-1} \sigma$ is not fillable and there is no
      fillable $n$-simplex\/ $\tau \neq \sigma$ of\/~$\nerve(\Pi)$ with\/ $c(\tau) \geq c$
      such that\/ $d_{c-1}\sigma \subseteq \del\tau$.
    \item If the face\/ $d_c\sigma$ is not fillable, then there exists a fillable
      $n$-simplex\/ $\tau$ of\/~$\nerve(\Pi)$ such that\/ $d_c\sigma \subseteq \del
      \tau$ and\/ $c(\tau) > c$.
  \end{enumerate}
\end{lemma}
\begin{proof}
  Let us write $\sigma= (p_0 < \dots < p_n)$ and choose minimal witnesses
  $\gamma_j$ of $p_{j-1} < p_j$.  Then $d_i \sigma = ( p_0 < \dots <
  \widehat{p_i} < \dots < p_n)$, where the circumflex signals omission of
  $p_i$. The minimal witnesses for the relations $p_{j-1} < p_j$ that occur in
  $d_i\sigma$ are thus all the $\gamma_j$ with $j \notin \{i,i+1\}$ and
  a minimal glob containing $\gamma_{i} \cup \gamma_{i+1}$.  
  
  Part~(a) is now obvious and so is the assertion that $d_{c-1}\sigma$ is not
  fillable since $\gamma_{c-1} \cup \gamma_c$ is neither contained in $G_1$ nor
  in $G_2$ by the definition of fillable simplices.  Towards the second claim
  in part~(b), it suffices to observe that the decomposition $\gamma_i \cup
  \gamma_{i+1}$ is unique with the property that $\gamma_i \subseteq G_1$ and
  $\gamma_{i+1} \cap G_1 \subseteq \del G_1$. There is thus no fillable simplex
  $\tau \neq \sigma$ with $c(\tau) \geq c$ and $d_{c-1}\sigma \subseteq \del
  \tau$.

  Finally, suppose that $d_c \sigma$ is not fillable. This is equivalent to the
  condition that $(\gamma_c \cup \gamma_{c+1}) \cap G_1 \nsubseteq \del G_1$.
  We may thus decompose $\gamma_c \cup \gamma_{c+1}$ into $\gamma_c' \cup
  \gamma_{c+1}'$ with $\gamma_c' \subseteq G_1$ and $\gamma_{c+1}' \cap G_1
  \subseteq \del G_1$. There then exists a path $q$ with $\gamma_c'$ and $\gamma_{c+1}'$ witnessing $p_{c-1} < q < p_{c+1}$ and the simplex
  \begin{equation*}
    \tau = ( p_0 < \dots p_{c-1} < q < p_{c+1} < \dots < p_n )
  \end{equation*}
  is a fillable simplex with $c(\tau) > c(\sigma)$ and $d_c \sigma \subseteq
  \del \tau$.
\end{proof}

We have finally gathered all technical details to finish the proof of
\autoref{thm:local}.
\begin{lemma}\label{lem:finishing-composition-lemma}
  The map \eqref{eq:last-map} admits a filtration by mid anodyne maps and
  is hence mid anodyne itself.
\end{lemma}
\begin{proof}
  Let us denote the domain of \eqref{eq:last-map} by $X_0$ and let $Y_{n,c}$ be
  the simplicial subset of $\nerve(\Pi)$ that is generated by the fillable
  nondegenerate simplices $\sigma$ of dimension $n$ with $c(\sigma) \geq c$
  together with the fillable nondegenerate simplices of dimension less than
  $n$. 
  Note that \begin{equation*}
    Y_{n-1,2} \cup X_0 = Y_{n-1,1} \cup X_0 = Y_{n,n+1} \cup X_0 
  \end{equation*} by
  \autoref{cor:fillabel-g-simplices}. The first equality follows from the
  fact that any fillable simplex $\sigma$ of $\nerve(\Pi)$ with
  $c(\sigma) = 1$ is already contained in $\colim \nerve(\Pi_\ph)
  \subseteq X_0$ and the second equality follows from the fact that any
  fillable $n$-simplex with $c(\sigma) = n+1$ is contained in $\colim
  \nerve(\Pi_\ph) \subseteq X_0$.
  We thus have a filtration
  \begin{equation*}
    X_0 =  Y_{2,3} \cup X_0 \subseteq Y_{2,2} \cup X_0 = Y_{2,1} \cup X_0 = Y_{3,4} \cup X_0 \subseteq Y_{3,3} \cup X_0 \subseteq \dots
  \end{equation*}
  of the inclusion $X_0 \to \nerve(\Pi)$
  that is exhaustive by \autoref{cor:any-is-in-boundary-of-fillable}. It therefore 
  suffices to show that each of the inclusions $Y_{n,c+1} \cup X_0 \subseteq
  Y_{n,c} \cup X_0$ with $2 \leq c \leq n$ is mid anodyne.
  To this end, it suffices to prove the following statements: 
  \begin{enumerate}[label=(\roman*)]
    \item For any nondegenerate simplex $\sigma \in Y_{n,c}$ that is
      not contained in $Y_{n,c+1} \cup X_0$ there is an inner horn
      $\Lambda^n_{c-1} \subseteq Y_{n,c+1} \cup X_0$ such that $\sigma$ is a
      filler for this horn.
    \item For any simplex $\sigma$ as in (i), the face $d_{c-1} \sigma$ is not
      contained in $Y_{n,c+1} \cup X_0$.
    \item If $\sigma \neq \tau$ are two simplices as in (i), then $d_{c-1}
      \sigma \neq d_{c-1} \tau$. 
  \end{enumerate}

  Consider a nondegenerate simplex $\sigma \in Y_{n,c} \smallsetminus (
  Y_{n,c+1} \cup X_0)$ for some $2 \leq c \leq n$. Note that $\sigma$ is fillable with $c(\sigma) = c$.
  Part (a) of \autoref{lem:fillable-boundaries} implies that all the faces
  $d_i\sigma$ of $\sigma$ with $i \notin \{c-1,c\}$ are fillable of dimension
  strictly less than $n$, i.\,e.\ contained in $Y_{n,c+1}$.  Moreover, by part
  (c) of \autoref{lem:fillable-boundaries}, $d_c\sigma$ is either fillable or a
  face of some $\tau \in Y_{n,c+1}$. In either case, it is contained in
  $Y_{n,c+1}$. The faces $d_i\sigma$ of $\sigma$ with $i \neq c-1$ are therefore
  an inner horn $\Lambda^n_{c-1} \subseteq Y_{n,c+1}$ of the desired kind since $0 < c-1 < n$ by assumption.
  
  Moreover, according to part (b) of \autoref{lem:fillable-boundaries},
  $d_{c-1}\sigma$ is not contained in $Y_{n,c+1}$. Let us assume for the sake of
  contradiction that $d_{c-1}\sigma$ is contained in 
  \begin{equation*}
    X_0 =  \nerve(\Sigma) \  \coprod_{ \mathclap{ \colim \nerve(\Sigma_\ph)  } }\ \colim
  \nerve(\Pi_\ph ),
  \end{equation*}
  i.\,e.\ $d_{c-1} \sigma \in \nerve(\Sigma)$ or $d_{c-1} \sigma
  \in \nerve( \Pi_i)$. As $\Sigma$ and $\Pi_i$ are closed under taking
  subdivisions, this implies $\sigma \in X_0$ by
  \autoref{lem:nerves-of-subdivision-closed}.

  Finally, if $\sigma$ and $\tau$ are two distinct fillable simplices with $c(\sigma) = c(\tau) = c$, \autoref{lem:fillable-boundaries} (b) implies $d_{c-1}\sigma \neq d_{c-1} \tau$.

\end{proof}

\chapter{A Pasting Theorem}\label{chap:thm}

In this final chapter, we collect all the results obtained so far and deduce
our pasting theorem.
Let us commence by combining \autoref{thm:global} and \autoref{thm:local}:
\begin{lastproposition}\label{prop:pasting-0}
  Let $\Sigma \to \Pi$ be an inclusion of complete pasting diagrams such that
  $\Sigma$ and $\Pi$ are both closed under taking subdivisions and contain all
  the interior faces of the underlying graph.
  The functor
  \begin{equation*}
    \dC[\Sigma] \to \dC[\Pi] 
  \end{equation*}
  has the left lifting property against all local mid fibrations.
\end{lastproposition}
\begin{proof}
  Recall that if $\Sigma \to \Pi$ satisfies the hypothesis of
  \autoref{thm:local}, then so do all the restrictions $\Sigma_{x,y} \to
  \Pi_{x,y}$ and all the inclusions $\Sigma_{x,y} \hc \Pi_{x,y} \to \Pi_{x,y}$.
  The proposition thus follows from \autoref{thm:local} and
  \autoref{thm:global}.
\end{proof}

Recall from
\cite{riehl-verity;homotopy-coherent-adjunctions-and-the-formal-theory-of-monads}
that the category of small simplicial categories is simplicially enriched with
mapping spaces given by $\icon(\dA,\dB)_n = \scat( \dA, \dB^{\Delta^n})$, where
$\dB^{\Delta^n}$ denotes the simplicial category with $\dB^{\Delta^n}(a,b) =
\dB(a,b)^{\Delta^n}$. The acronym $\icon$ stands for \enquote{{i}dentity
{c}omponent {o}plax {n}atural transformations} and is due to
Lack in the case of $2$-categories, see \cite{lack;icons}.
\begin{lastproposition}\label{prop:pasting-1}
  Consider an inclusion $\Sigma \to \Pi$ of complete pasting diagrams such that
  both $\Sigma$ and\/ $\Pi$ are closed under taking subdivisions and contain all
  the interior faces of the underlying graph. Further let\/ $\dB \to \dA$ be a
  local mid fibration of simplicial categories.
  The canonical map
  \begin{equation*}
    \icon\bigl( \dC[\Pi], \dB \bigr) \to \icon\bigl( \dC[\Sigma], \dB \bigr)\ \limtimes_{\substack{\\[.5ex] \mathclap{\icon(\dC[\Sigma],\dA )}}}\ \icon\bigl( \dC[\Pi], \dA \bigr)
  \end{equation*}
  is a trivial Kan fibration. 
\end{lastproposition}
\begin{proof}
  Any lifting problem
  \begin{equation*}
    \begin{tikzpicture}[diagram]
      \matrix[objects] {
        |(del)| \del \Delta^n \& |(pi-b)| \icon\bigl( \dC[\Pi], \dB \bigr) \\[-2ex]
        |(simp)| \Delta^n  \& |(pb)| \displaystyle \icon\bigl( \dC[\Sigma], \dB \bigr)\ \limtimes_{\substack{\\[.5ex] \mathclap{\icon(\dC[\Sigma],\dA )}}}\ \icon\bigl( \dC[\Pi], \dA \bigr) \\
      };
      \path[maps,->] 
        (del) edge (simp)
        (del) edge (pi-b)
        (pi-b) edge (pb)
        (simp) edge (pb)
      ;
    \end{tikzpicture}
  \end{equation*}
  transposes to a lifting problem
  \begin{equation*}
    \begin{tikzpicture}[diagram]
      \matrix[objects] {
        |(sigma)| \dC[\Sigma] \& |(b)| \dB^{\Delta^n} \\[-2ex]
        |(pi)| \dC[\Pi] \& |(a)| \displaystyle\dB^{\del \Delta^n}\ \limtimes_{ \substack{\\[.5ex] \mathclap{ \dA^{\del \Delta^n} } } }\ \dA^{\Delta^n} \smash{.}\\
      };
      \path[maps,->] 
        (b) edge (a)
        (sigma) edge (pi)
        (sigma) edge (b)
        (pi) edge (a)
      ;
    \end{tikzpicture}
  \end{equation*}
  Similar to the proof of
  \cite[Lemma~4.4.2]{riehl-verity;homotopy-coherent-adjunctions-and-the-formal-theory-of-monads},
  we observe that the map on the right hand side in this latter diagram is
  locally a mid fibration and the claim thus follows from
  \autoref{prop:pasting-0}.
\end{proof}

Recall from \autoref{thm:labeling} that any labeling $\Lambda$ of a globular
graph $G$ in some simplicial category $\dA$ determines a simplicial functor
$u_\Lambda \colon \dC[\Sigma_{\min}^c] \to \dA$. 
\begin{lastdefinition}\label{def:space-of-comp}
  The space $C(\Lambda)$ of compositions of a labeling $\Lambda$ of a globular
  graph $G$ in some simplicial category $\dA$ is obtained as the pullback in
  \begin{equation*}
    \begin{tikzpicture}[diagram]
      \matrix[objects] {
        |(c)| C(\Lambda) \& |(icon1)| \icon\bigl( \dC[\Pi_{\max}], \dA\bigr) \\
        |(p)| \Delta^0 \& |(icon2)| \icon\bigl( \dC[\Sigma_{\min}^c], \dA \bigr)\smash{.} \\
      };
      \path[maps,->] 
        (c) edge (icon1)
        (c) edge (p)
        (p) edge node[below] {$ u_\Lambda$} (icon2)
        (icon1) edge (icon2)
      ;
    \end{tikzpicture}
  \end{equation*}
\end{lastdefinition}

\begin{lastremark}
  Using the fact that $u_\Lambda$ and $\Lambda$ determine each other uniquely,
  we see that the $0$-simplices in $C(\Lambda)$ are those functors $v \colon
  \dC[\Pi_{\max}] \to \dA$ such that $\Lambda_v = \Lambda$, i.\,e.\ 
  extensions of $\Lambda$ to the fully coherent $\dC[\Pi_{\max}]$.
\end{lastremark}

\pastingthm

\begin{proof}
  Observe that $\dA \to \pt$ is a mid fibration in the case that $\dA$ is
  enriched over quasi-categories. \autoref{prop:pasting-1} thus implies that the map
  \begin{equation*}
    \icon\bigl(\dC[\Pi_{\max}],\dA \bigr) \to \icon\bigl(\dC[\Sigma_{\min}^c],\dA \bigr) 
  \end{equation*}
  is a trivial Kan fibration. The space $C(\Lambda)$ of compositions of
  $\Lambda$ is hence a contractible Kan complex.
\end{proof}

\begin{lastremark}
  Note that we can deduce Power's original result from
  \cite{power;a-2-categorical-pasting-theorem} from \autoref{thm:pasting}.
  Indeed, given any labeling $\Lambda$ of a globular graph $G$ in some
  $2$-category $A$, we obtain a labeling of $G$ in the simplicial category
  $\dA$ obtained by applying the nerve functor to all the categories $A(a,b)$.
  As $C(\Lambda)$ is nonempty there is at least one extension $v \colon
  \dC[\Pi_{\max}] \to \dA$ of $\Lambda$ and $v$ restricts to a map
  \begin{equation*}
    v\colon \nerve(G) = \dC[\Pi_{\max}](s,t) \to \dA(s,t)
  \end{equation*}
  of simplicial sets. We now apply the left adjoint $\tau_1$ of $\nerve$ to the inclusion 
  \begin{equation*}
    \Delta^1 \xto{\del G} \nerve(G) \to \dA(s,t)
  \end{equation*}
  and thus obtain a composite $2$-cell $\phi \colon f \to g$ in $A(s,t) =
  \tau_1 \nerve A(s,t)$.

  Now consider two such extensions $v,v' \in \icon( \dC[\Pi_{\max}],\dA)$ with
  associated composite $2$-cells $\phi$ and $\psi$.  We then find a lift $h$ as
  in the diagram
  \begin{equation}\label{dia:h-construction}
    \begin{tikzpicture}[diagram]
      \matrix[objects] {
        |(del)| \del \Delta^1 \& \& |(c)| \icon\bigl( \dC[\Pi_{\max}],\dA\bigr) \\
        |(1)| \Delta^1 \& |(0)| \Delta^0 \& |(u)| \icon\bigl( \dC[\Sigma_{\min}^c],\dA \bigr) \smash{.}\\
      };
      \path[maps,->] 
        (del) edge (1)
        (del) edge node[above] {$v \amalg v'$} (c)
        (1) edge (0)
        (0) edge node[below] {$u_\Lambda$} (u)
        (c) edge (u)
        (1) edge node {$h$} (c)
      ;
    \end{tikzpicture}
  \end{equation}
  As above, this gives rise to a map 
    $h \colon \Delta^1 \times \nerve(G) \to \dA(s,t)$
    of simplicial sets and hence to a square
    \begin{equation*}
    \begin{tikzpicture}[diagram]
      \matrix[objects] {
        |(1)| \bullet \& |(2)| \bullet \\
        |(3)| \bullet \& |(4)| \bullet \\
      };
      \path[maps,->] 
        (1) edge node[left] {$\phi$} (3)
        (2) edge node[left] {$\psi$} (4)
        (1) edge node[above] {$u$} (2)
        (3) edge node[below] {$v$} (4)
      ;
    \end{tikzpicture}
    \end{equation*}
    in $A(s,t)$. Commutativity of the lower triangle in \eqref{dia:h-construction} 
    now implies that both $u$ and $v$ are the identity, i.\,e.\ $\phi = \psi$.
\end{lastremark}

\printbibliography
\end{document}